\pdfoutput=1
\RequirePackage{silence}
\documentclass[12pt,a4paper,svgnames]{amsart}
\usepackage[hmarginratio={1:1},vmarginratio={1:1},lmargin=60.0pt,tmargin=60.0pt]{geometry}

\usepackage{booktabs}
\allowdisplaybreaks

\usepackage[T1]{fontenc}
\usepackage{latexsym,exscale,amsfonts,amssymb,mathtools}
\usepackage{amsmath,amsthm,amsfonts,amssymb,amscd,textcomp,bbm}
\usepackage{mathrsfs,stackrel}
\usepackage{etoolbox}
\usepackage{tikz,tikz-cd}
\usetikzlibrary{matrix,arrows.meta,decorations.markings,shapes.multipart,decorations.pathreplacing,decorations.pathmorphing}
\usetikzlibrary{shadings}

\pgfdeclareverticalshading{rainbow}{100bp}
{color(0bp)=(red); color(25bp)=(red); color(35bp)=(yellow);
color(45bp)=(green); color(55bp)=(cyan); color(65bp)=(blue);
color(75bp)=(violet); color(100bp)=(violet)}

\tikzset{anchorbase/.style={baseline={([yshift=-0.5ex]current bounding box.center)}},}
\usepackage{enumitem}
\setlist[enumerate]{itemsep=0.15cm,label=\emph{\upshape(\alph*)}}
\setlist[enumerate,2]{itemsep=0.15cm,label=\emph{\upshape(\roman*)}}
\usepackage{aliascnt}
\usepackage{ytableau}
\usepackage{xparse}
\usepackage{cite}

\usepackage{dynkin-diagrams}

\usepackage{array}
\newcolumntype{C}{>{$}c<{$}}

\usepackage{xcolor,colortbl}

\definecolor{mygray}{gray}{0.6}
\definecolor{mygraydark}{gray}{0.4}
\definecolor{mygraylight}{gray}{0.85}
\definecolor{spinach}{RGB}{46,139,87}
\definecolor{tomato}{RGB}{255,99,71}
\definecolor{orchid}{RGB}{143,40,194}
\definecolor{neon}{RGB}{77,77,255}
\definecolor{pumpkin}{RGB}{224,180,80}
\definecolor{citron}{RGB}{190,180,90}

\definecolor{lava}{RGB}{207,16,32}
\definecolor{cream}{RGB}{255,253,208}
\definecolor{verdigris}{RGB}{67,179,174}
\definecolor{Black}{RGB}{0,0,0}
\definecolor{mydarkblue}{RGB}{10,10,170}
\definecolor{darkspinach}{RGB}{20,70,20}
\definecolor{darktomato}{RGB}{155,40,30}
\definecolor{darkorchid}{RGB}{50,10,100}
\definecolor{darklava}{RGB}{150,8,16}

\usepackage{enumitem}
\setlist[enumerate]{itemsep=0.15cm,label=\emph{\upshape(\alph*)}}
\setlist[enumerate,2]{itemsep=0.15cm,label=\emph{\upshape(\roman*)}}
\setlist[enumerate,3]{itemsep=0.15cm,label=\emph{\upshape(\Alph*)}}

\let\emph\relax
\DeclareTextFontCommand{\emph}{\bfseries\em}

\newcommand{\acts}{\centerdot}
\newcommand{\actsleft}{\mathop{\,\;\raisebox{1.7ex}{\rotatebox{-90}{$\circlearrowright$}}\;\,}}
\newcommand{\actsright}{\mathop{\,\;\raisebox{0ex}{\rotatebox{90}{$\circlearrowleft$}}\;\,}}
\renewcommand{\dots}{\text{...}}

\newcommand{\placeholder}{{}_{-}}
\newcommand{\mystrut}{\rule[-0.2\baselineskip]{0pt}{0.9\baselineskip}}

\usepackage{todonotes}

\def\changed#1{{#1}}
\def\ochanged#1{{#1}}

\renewcommand{\to}{\rightarrow}

\def\modules{{\text{-mod}}}

\newcommand{\refequal}[1]{\xy {\ar@{=}^{#1}
(-1,0)*{};(1,0)*{}};
\endxy}

\newcommand{\C}{\mathbb{C}}
\newcommand{\Z}{\mathbb{Z}}

\newcommand{\Q}{\mathbb{Q}}
\newcommand{\F}{\mathbb{F}}
\newcommand{\N}{\mathbb{Z}_{\geq 0}}

\newcommand{\intform}{\mathbb{A}}
\newcommand{\field}{\mathbb{F}}
\newcommand{\ring}{\mathbb{K}}

\DeclareMathOperator{\Hom}{{\rm Hom}}

\DeclareMathOperator{\End}{{\rm End}}

\DeclareMathOperator{\Ext}{{\rm Ext}}

\DeclareMathOperator{\Res}{{\rm Res}}

\DeclareMathOperator{\id}{{\rm id}}

\newcommand{\xsym}{u}

\newcommand{\qbinn}[2]{\genfrac{[}{]}{0pt}{}{#1}{#2}}

\DeclareMathOperator{\dv}{X}
\DeclareMathOperator{\dd}{\partial}

\DeclareMathOperator{\pr}{\textbf{1}}

\DeclareMathOperator{\wt}{\rm{wt}}

\DeclareMathOperator{\ee}{\tilde{{e}}}

\DeclareMathOperator{\ff}{\tilde{{f}}}

\DeclareMathOperator{\canB}{\textbf{B}}
\DeclareMathOperator{\dcanB}{\dot{\textbf{B}}}

\def\NewTheorem#1{%
\newaliascnt{#1}{equation}%
\newtheorem{#1}[#1]{#1}%
\aliascntresetthe{#1}%
\expandafter\def\csname #1autorefname\endcsname{#1}%
}
\def\equationautorefname~#1\null{(#1)\null}

\numberwithin{equation}{subsection}

\NewTheorem{Proposition}
\NewTheorem{Theorem}
\NewTheorem{Conjecture}
\NewTheorem{Corollary}
\AtEndEnvironment{Corollary}{\null\hfill$\square$}%
\NewTheorem{Lemma}
\theoremstyle{definition}
\NewTheorem{Definition}
\AtEndEnvironment{Definition}{\null\hfill$\Diamond$}%
\NewTheorem{Notation}
\AtEndEnvironment{Notation}{\null\hfill$\Diamond$}%
\NewTheorem{Example}
\AtEndEnvironment{Example}{\null\hfill$\Diamond$}%
\NewTheorem{Examples}
\AtEndEnvironment{Examples}{\vskip-10mm\null\hfill$\Diamond$}%

\theoremstyle{remark}
\NewTheorem{Remark}
\AtEndEnvironment{Remark}{\null\hfill$\Diamond$}%
\NewTheorem{Question}
\NewTheorem{Assumption}

\usepackage[hypertexnames=false]{hyperref}
\usepackage{bookmark}
\hypersetup{
pdftoolbar=true,
pdfmenubar=true,
pdffitwindow=false,
pdfstartview={FitH},
pdftitle={On a symplectic quantum Howe duality},
pdfauthor={Elijah Bodish and Daniel Tubbenhauer},
pdfsubject={},
pdfcreator={Elijah Bodish and Daniel Tubbenhauer},
pdfproducer={Elijah Bodish and Daniel Tubbenhauer},
pdfkeywords={},
pdfnewwindow=true,
colorlinks=true,
linkcolor=mydarkblue,
citecolor=teal,
filecolor=magenta,
urlcolor=orchid,
linkbordercolor=lava,
citebordercolor=teal,
urlbordercolor=orchid,
linktocpage=true
}

\newcommand{\nnfootnote}[1]{%
\begin{NoHyper}
\renewcommand\thefootnote{}\footnote{#1}%
\addtocounter{footnote}{-1}%
\end{NoHyper}
}

\def\makeautorefname#1#2{\csdef{#1autorefname}{#2}}
\makeautorefname{section}{Section}%
\makeautorefname{subsection}{Section}%
\makeautorefname{subsubsection}{Section}%

\begin{document}
\title[On a symplectic quantum Howe duality]{On a symplectic quantum Howe duality}
\author[E. Bodish and D. Tubbenhauer]{Elijah Bodish$^{\ast}$ and Daniel Tubbenhauer}
\thanks{$^\ast$Corresponding author.}

\address{E.B.: MIT, Department of Mathematics, Building 2, Office 2-178, Cambridge, MA 02139, https://orcid.org/0000-0003-1499-5136}
\email{ebodish@mit.edu}

\address{D.T.: The University of Sydney, School of Mathematics and Statistics F07, Office Carslaw 827, NSW 2006, Australia, \href{http://www.dtubbenhauer.com}{www.dtubbenhauer.com}, https://orcid.org/0000-0001-7265-5047}
\email{daniel.tubbenhauer@sydney.edu.au}

\begin{abstract} 
We prove a nonsemisimple quantum version of 
Howe's duality with the rank 2n symplectic and the rank 2 special linear group acting on the 
exterior algebra of type C. 
We also discuss the first steps towards the symplectic analog of harmonic analysis on quantum spheres, give character formulas for various fundamental modules, and construct canonical bases of the exterior algebra.
\end{abstract}

\nnfootnote{\textit{Mathematics Subject Classification 2020.} Primary: 
17B37, 20G42; Secondary: 05E10, 20G05, 22E46.}
\nnfootnote{\textit{Keywords.} Quantum Howe duality, tilting modules, Weyl characters, canonical bases, symplectic groups.}

\addtocontents{toc}{\protect\setcounter{tocdepth}{1}}

\maketitle

\tableofcontents

\section{Introduction}\label{S:Intro}

In this paper we prove a nonsemisimple and quantum version of 
one of Howe's dualities. We work mostly over $\Z[q,q^{-1}]$ 
and specializations to any field and any quantum parameter are allowed.

\subsection{Howe dualities}

A cornerstone of modern invariant theory are \emph{Schur--Weyl--Brauer dualities}, which allow 
one to play two commuting actions on tensor space against one another. See for example 
\cite[Section 3]{AnStTu-semisimple-tilting} for a summary, using the language of this paper.

Building on earlier work as e.g. in \cite{Ho-remarks-invariant-theory}, 
these dualities were generalized far beyond tensor space by Howe in the legendary 1995 
Schur lectures \cite{Ho-perspectives-invariant-theory}. Howe's goal, \ochanged{which} was achieved
magnificently, was to present a new approach to classical invariant theory.

As Howe points out, classical invariant theory, historically speaking, 
was mostly about group actions on symmetric algebras 
or, less \ochanged{classically}, on exterior algebras. Howe reformulated these actions in the 
\emph{double centralizer}
(or double commutant) approach of Schur--Weyl--Brauer dualities, and the \ochanged{result} 
is what we call \emph{Howe dualities}. 

An example of a Howe duality 
(here formulated for the Lie algebra instead of the Lie group as 
Howe did in \cite{Ho-perspectives-invariant-theory}) reads as follows.
Consider the $\C$-vector space $\Lambda_{\C}:=\Lambda_{\C}(\C^{m}\otimes\C^{n})$, 
i.e. the exterior 
algebra of $\C^{m}\otimes\C^{n}$. Then:
\begin{enumerate}

\item There are commuting actions
\begin{gather*}
U_{\C}(\mathfrak{gl}_{m})\actsleft
\Lambda_{\C}
\actsright
U_{\C}(\mathfrak{gl}_{n}).
\end{gather*}	

\item Let $\phi_{\C}$ and $\psi_{\C}$ be the 
$\C$-algebra homomorphisms 
induced by the two actions from (a). Then:
\begin{gather*}
\phi_{\C}\colon
U_{\C}(\mathfrak{gl}_{m})\twoheadrightarrow\End_{U_{\C}(\mathfrak{gl}_{n})}(\Lambda_{\intform})
,\quad
\psi_{\C}\colon
U_{\C}(\mathfrak{gl}_{n})\twoheadrightarrow\End_{U_{\C}(\mathfrak{gl}_{m})}(\Lambda_{\intform})
.
\end{gather*}
That is, the two actions generate the others centralizer.

\item There is an explicit $U_{\C}(\mathfrak{gl}_{m})\otimes U_{\C}(\mathfrak{gl}_{n})^{op}$-module decomposition of $\Lambda_{\C}$ pairing a Weyl $U_{\C}(\mathfrak{gl}_{m})$-module
and a dual Weyl $U_{\C}(\mathfrak{gl}_{n})$-module.

\end{enumerate}
This is known as \emph{exterior type A Howe duality} and appears 
in \cite[Theorem 4.1.1]{Ho-perspectives-invariant-theory} as
one of many examples considered by Howe.

\begin{Remark}
The more classical version is \emph{symmetric type A Howe duality}, but the 
above is closer to the results in this paper, so we decide to prefer exterior 
over symmetric type A Howe duality.
\end{Remark}

As Howe explains, when working over $\C$, the various forms of Howe dualities are 
equivalent to the respective Schur--Weyl--Brauer dualities.
In particular, Howe-type dualities have 
been of crucial importance for invariant theory every since, 
but are also pervasive in other fields. Most notably, 
representation theory, low dimensional topology, and categorification.

Let us recall two generalizations of Howe's dualities.

\begin{Remark}
We are interested in nonsemisimple quantum versions of Howe dualities 
and, in order to not lose focus, we will ignore other possible generalizations like super versions.
\end{Remark}

A first possible generalization is to study \emph{nonsemisimple versions} 
of Howe's dualities\changed{, for example, working over the
ground ring $\Z$}. Note hereby that \cite{Ho-perspectives-invariant-theory} 
stays over $\C$ which makes the story semisimple and character-type 
arguments apply as explained masterfully by Howe.
In the nonsemisimple cases such 
character arguments are not as powerful as over $\C$, and one needs some form 
of replacement. While (a) is straightforward and 
characteristic independent, replacements for Howe's arguments showing (b) 
and (c) are needed.

As explained for Schur--Weyl(--Brauer) dualities 
in \cite{DuPaSc-schur-weyl-tilting} and for Howe dualities in 
\cite{AdRy-tilting-howe-positive-char}
the theory of tilting modules in the spirit of 
\cite{Ri-good-filtrations} and \cite{Do-tilting-alg-groups} 
usually plays a crucial role. 
In particular, \cite{AdRy-tilting-howe-positive-char} studies exterior versions 
of Howe dualities using the tilting module approach. 
Whenever tilting theory is not available, like for the symmetric invariants, 
the nonsemisimple study of invariant theory tends to be much trickier, see e.g. 
\cite{CoPr-invariant-theory}, and is less 
accessible via Howe's approach.

\begin{Remark}
As stated above, over $\C$ or, more generally, in the semisimple case,
Schur--Weyl--Brauer dualities and Howe dualities are equivalent. This is no longer true 
in the nonsemisimple situation, and \cite{DuPaSc-schur-weyl-tilting} and 
\cite{AdRy-tilting-howe-positive-char} seem to have appeared independently.
\end{Remark}

A second possible generalization are \emph{quantum versions} involving 
some form of quantum groups instead of Lie groups. For starters, let us 
assume \ochanged{that} we work over $\C(q)$ for a formal parameter $q$. In this case 
quantum Howe dualities are semisimple, so Howe's character-type arguments apply.
In other words, (b) and (c) are rather easy to prove. However, (a) can be quite difficult.

Part of the problem is that quantum groups, and quantizations in general, often crucially 
depend on the involved choices of quantization. Sometimes one needs to work 
hard to get (a) for quantum groups instead of Lie groups, sometimes
one even needs to give up on quantum groups altogether.

First examples of quantum versions of Howe's dualities appeared in the 90s, see e.g. 
the work \cite{NoUmWa-sl2-son-duality}, and indicated that nonstandard quantizations 
often play a role. See also \cite{EhSt-nw-algebras-howe} and \cite{SaTu-bcd-webs} 
for the same type of phenomena. Sometimes, however, Howe's dualities quantize 
nicely (meaning with the standard Drinfeld--Jimbo quantum groups). In type A, see for example 
\cite{LeZhZh-q-first-fundamental-theorem} or \cite{CaKaMo-webs-skew-howe} 
for exterior, \cite{RoTu-symmetric-howe} 
for symmetric, \cite{TuVaWe-super-howe} 
for exterior-symmetric and \cite{LaTuVa-verma-howe} for Verma versions of nicely quantized Howe dualities. But again, 
proving (a) might still be tricky, in particular outside of type A, as e.g. in \cite{Su-harmonic-analysis-son}.

Combining both, the nonsemisimple and quantum setting, \ochanged{then needs} new arguments
(compared to \cite{Ho-perspectives-invariant-theory}) for (a), (b) and (c).
Modulo this paper and our ignorance, the only 
known nonsemisimple quantum Howe duality is the exterior type A Howe duality 
from above. In the exterior type A Howe duality case it turns out that (a), (b) and (c) are rather easy to prove, but that is a coincidence since the exterior powers of the vector representation in type A 
are tilting regardless of the characteristic.
This is implicit already in \cite{LeZhZh-q-first-fundamental-theorem}
and \cite{CaKaMo-webs-skew-howe}, and has then been pointed out explicitly,
and masterfully used, in \cite{El-ladders-clasps}.

\subsection{What this paper does}

The starting point of this paper is to take the ``easiest'' of Howe's dualities outside of 
type A and to prove it in the nonsemisimple and quantum case. 

In \cite{Ho-perspectives-invariant-theory} there are four essentially equivalent dualities listed:
\begin{enumerate}[itemsep=0.15cm,label=\emph{\upshape(\Alph*)}]

\item We have $SP_{2n}$ acting on $\Lambda\big((\C^{2n})^{\otimes m}\big)$
or $O_{n}$ acting on $\mathrm{Sym}\big((\C^{n})^{\otimes m}\big)$ with the dual action by $SP_{2m}$.

\item We have $O_{n}$ acting on $\Lambda\big((\C^{n})^{\otimes m}\big)$ 
or $SP_{2n}$ acting on $\mathrm{Sym}\big((\C^{2n})^{\otimes m}\big)$ with the dual action by $SO_{2m}$.

\end{enumerate}
Since we want nonsemisimple versions, thus tilting theory should come in handy, 
the exterior cases are certainly easier than the symmetric ones, so we 
stay with these. Moreover, \cite{SaTu-bcd-webs}
suggests that quantization might be an issue in general, but we have one special case where 
quantization is nice due to the small number coincidence that $SP_{2}\cong SL_{2}$, see \autoref{R:HowePartOneCoideal} for a detailed discussion.

It turns out that one is lucky and \emph{one case remains}: $SP_{2n}$ acting on $\Lambda((\C^{2n})^{\otimes 1})$ with the dual action by $SP_{2\cdot 1}\cong SL_{2}$. This is the content of the present paper.

That is, our main result, (a version of) \emph{symplectic quantum Howe duality}, \ochanged{is} as follows.
Let $\intform=\Z[q,q^{-1}]$ be the integral form where $q$ is a variable, and let $\Lambda_{\intform}$ the exterior algebra of the symplectic quantum vector representation. Then:

\begin{Theorem}\label{T:IntroMain}
Consider the free $\intform$-module $\Lambda_{\intform}$.

\begin{enumerate}

\item There are commuting actions
\begin{gather*}
U_{\intform}(\mathfrak{sp}_{2n})\actsleft
\Lambda_{\intform}
\actsright
U_{\intform}(\mathfrak{sl}_{2})^{op}.
\end{gather*}	

\item Let $\phi_{\intform}$ and $\psi_{\intform}$ be the 
$\intform$-algebra homomorphisms 
induced by the two actions from (a). Then:
\begin{gather*}
\phi_{\intform}\colon
U_{\intform}(\mathfrak{sp}_{2n})\twoheadrightarrow\End_{U^{\intform}(\mathfrak{sl}_{2})^{op}}(\Lambda_{\intform})
,\quad
\psi_{\intform}\colon
\changed{U_{\intform}(\mathfrak{sl}_{2})}\twoheadrightarrow\End_{U_{\intform}(\mathfrak{sp}_{2n})}(\Lambda_{\intform})
.
\end{gather*}
That is, the two actions generate the others centralizer.

\item Let $\textbf{U}_{\intform}:=U_{\intform}(\mathfrak{sp}_{2n})\otimes U_{\intform}(\mathfrak{sl}_{2})$. Then $\Lambda_{\intform}$ is \ochanged{a 
Howe} tilting $\textbf{U}_{\intform}$-module
(see \autoref{D:ETilting} for the definition), that is, 
$\Lambda_{\intform}$ is tilting for $U_{\intform}(\mathfrak{sp}_{2n})$ and $U_{\intform}(\mathfrak{sl}_{2})$, as 
well as full tilting for $\End_{U_{\intform}(\mathfrak{sp}_{2n})}(\Lambda_{\intform})$ and $\End_{U_{\intform}(\mathfrak{sl}_{2})}(\Lambda_{\intform})$, both
separately, and satisfies
\begin{gather}\label{Eq:MainOne}
\Hom_{U_{\intform}(\mathfrak{sp}_{2n})}\big(\Delta_{\intform}(\varpi_{k}),\Lambda_{\intform}\big)\cong
\nabla_{\intform}(n-k)
,\quad
\Hom_{U_{\intform}(\mathfrak{sl}_{2})}\big(\Delta_{\intform}(n-k),\Lambda_{\intform}\big)\cong
\nabla_{\intform}(\varpi_{k})
.
\end{gather}
Moreover, in the Grothendieck ring
\begin{gather}\label{Eq:MainTwo}
[\Lambda_{\intform}] 
= 
\sum_{k=0}^{n} 
[\Delta_{\intform}(\varpi_{k})\otimes\nabla_{\intform}(n-k)].
\end{gather}
Here $\Delta_{\intform}(\varpi_{k})$ and $\nabla_{\intform}(n-k)$ denote the 
Weyl $U_{\intform}(\mathfrak{sp}_{2n})$-module and dual Weyl 
$U_{\intform}(\mathfrak{sl}_{2})$-module of the indicated highest weights.

\end{enumerate}
\end{Theorem}

Note that $U_{\intform}(\mathfrak{sl}_{2})$ is 
isomorphic to $U_{\intform}(\mathfrak{sp}_{2})$. Thus, 
\autoref{T:IntroMain} is the nonsemisimple and quantum version of 
symplectic Howe duality as in \cite[Theorem 3.8.9.3]{Ho-perspectives-invariant-theory} 
in the special case $m=1$.

\begin{Remark}
In the semisimple case \autoref{Eq:MainOne} implies the $\textbf{U}_{\intform}$-module decomposition
\begin{gather*}
\Lambda_{\intform}\cong\bigoplus_{k=0}^{n}\Delta_{\intform}(\varpi_{k})\otimes\nabla_{\intform}(n-k),
\end{gather*}
and \autoref{Eq:MainTwo} its combinatorial shadow. \ochanged{Thus, we can} view \autoref{T:IntroMain}.(c) as a nonsemisimple analog of 
(c) of the respective Howe duality.
\end{Remark}

\begin{Remark}\label{R:IntroMain}
\leavevmode	

\begin{enumerate}

\item The ``op'' makes an appearance in \autoref{T:IntroMain} since we 
have two left actions, one of $U_{\intform}(\mathfrak{sp}_{2n})$ and one of 
$U_{\intform}(\mathfrak{sl}_{2})$, instead of a left and a right action.

\item As we will elaborate in \autoref{R:HowePartOneCoideal}, it is rather 
surprising that \cite[Theorem 3.8.9.3, $m=1$]{Ho-perspectives-invariant-theory} can be quantized 
as in \autoref{T:IntroMain}, namely involving two standard quantum groups.

\item The nonsemisimple but non-quantum versions of the exterior cases in (A) and (B) above 
are discussed in \cite{AdRy-tilting-howe-positive-char}, see e.g. 
\cite[Theorem 2.1]{AdRy-tilting-howe-positive-char}.

\item Although formulated differently, 
the result \cite[Theorem 3.4]{Su-harmonic-analysis-son} can be interpreted as the
orthogonal version of \autoref{T:IntroMain} corresponding to (A) for $m=1$ above. Note that 
\cite{Su-harmonic-analysis-son} works over $\Q(q^{1/2})$ for a variable $q$, i.e. the semisimple case,
which \ochanged{is} somewhat expected due to the appearance of the symmetric algebra. 
(The fraction power of $q$ is related to the usage of the braiding in \cite{Su-harmonic-analysis-son} 
and does not play any essential role.) Note also that 
\cite{Su-harmonic-analysis-son} has a nice quantization, again due to the 
small number coincidence $SP_{2}\cong SL_{2}$.

\end{enumerate}
We will elaborate on these points in the main text.
\end{Remark}

Our proof of \autoref{T:IntroMain} splits into two parts: we first 
prove the semisimple version and then go to the nonsemisimple case.
Along the way we sketch the beginning of the symplectic analog of harmonic analysis on quantum spheres, prove some consequences of our Howe duality for the characters and various multiplicities 
of fundamental simple, Weyl or tilting modules in type C, and we 
construct canonical bases of $\Lambda_{\intform}$.

\subsection{Further directions}

Here are a few ideas to continue the story:

\begin{enumerate}

\item As we explain in more details in \autoref{R:HowePartOneCoideal}, our 
quantization results can be seen as a ``dual'' to \cite{SaTu-bcd-webs} for $m=1$.
A version for $m>1$ could give a Howe duality approach to quantum
BCD web categories as in e.g. \cite{BoElRoTa-c-webs} or their various 
applications to tilting modules as e.g. in \cite{Bo-c2-tilting}.
These webs could then be compared \ochanged{with} the $q=1$ version of \cite{SaTu-bcd-webs}.

\item Some of the most important modules in this paper are simple, Weyl and tilting modules for fundamental weights of the symplectic group. Not much is known about these modules for the nonsemisimple quantum case. Our results, \ochanged{particularly} those in \autoref{SS:Tilting}, 
already give the Weyl character of the 
fundamental tilting modules; in 
particular, we can tell when a Weyl module 
is simple. These results might 
help to generalize some of the known facts 
about these \ochanged{modules} such as those that can be found in, for example, 
\cite{PrSu-fund-weyl-sp}, \cite{McNi-charp-howe}, \cite{Fo-symplectic-simple} or \cite{GiMa-modular-dim-symplectic}.

\item The results of \autoref{SS:Canonical} suggest 
the existence of a categorification of our story.

\end{enumerate}

\begin{Remark}
Some calculations in this paper were done using code. 
If the reader wants to run that code themselves, then 
they can find it here:

\url{https://github.com/dtubbenhauer/sp2n-fundamental-tilting-charaters}

\url{https://github.com/dtubbenhauer/symplectic-diff-operators}. 
\end{Remark}


\noindent\textbf{Acknowledgments.}
We like to thank Ben Elias, Pavel Etingof, Sasha Kleshchev, and Joanna Meinel for
freely sharing ideas, pointers to the
literature and stimulating discussions. \ochanged{Special thanks to the referee for their thorough review of the document, providing many helpful comments, and spotting an incredible amount of nonsense on our part.} 
Special thanks to Sasha Kleshchev for (implicitly) suggesting this project, and to Magma, Mathematica and SageMath for computational help.
D.T. thanks the rainbow flag for supporting and inspiring them during the years.

\section{The module and the two commuting actions}\label{S:ExtAlg}

Throughout, and unless stated otherwise, e.g. in 
specific examples, we fix $n\in\Z_{\geq 1}$.

\begin{Remark}
All colors in this paper are a visual aid only.
\end{Remark}

\subsection{Quantum exterior algebras of type C}

We start by stating our quantum notation:

\begin{Notation}
For a parameter $q$ let $\intform:=\mathbb{Z}[q,q^{-1}]$ and $\field:=\text{frac}(\intform)=\mathbb{Q}(q)$. Working with $\intform$ 
as the ground ring means working \emph{integrally}, working with $\field$ makes the story \emph{semisimple}.
The specializations $\intform\to\C$ given by $q\mapsto 1$ or 
$q\mapsto -1$ are called the \emph{classical case}.

In general:
\begin{enumerate}

\item We use the subscript $\intform$ to denote integral versions of modules, divided power versions of algebras, or, more generally, any notion that involves $\intform$.

\item We use the subscript $q$ to denote the $\field$ versions of modules and quantum groups,
or, more generally, any notion that involves $\field$.

\item Similarly, we use the subscript $1$ to denote the classical counterparts. 
(The difference between $1$ and $-1$ does not play an essential role for us.)

\item We will later use arbitrary specializations which are denoted by $\ring$.

\item In general, we often only define $\intform$-versions and their variants over 
other rings can be obtained by scalar extension, or scalar extension and specialization.

\end{enumerate}
For $k\in\N$ let
\begin{gather*}
[k]_{\intform}:=(q^{k}-q^{-k})/(q-q^{-1})=q^{-k+1}+q^{-k+3}+\dots+q^{k-3}+q^{k-1}
,
\end{gather*}
denote the usual quantum numbers. For $k\in\Z_{<0}$ let 
$[k]_{\intform}=-[-k]_{\intform}$.
For $k\in\Z$ and $l\in\N$ also define
\begin{gather*}
[l]_{\intform}!:=[l]_{\intform}[l-1]_{\intform}\dots[1]_{\intform}
,\quad
\begin{bmatrix}k\\l\end{bmatrix}_{\intform}
:=
\frac{[k]_{\intform}[k-1]_{\intform}\dots[k-l+1]_{\intform}}{[l]_{\intform}!}
.
\end{gather*}
We let $q_{1}=\dots=q_{n-1}=q$ and $q_{n}=q^{2}$. As usual, we can then define the quantum numbers (and factorials and binomials etc.) 
in $q_{i}$ instead of $q$ and indicate this using a subscript $i$.
\end{Notation}

\begin{Notation}
Using interval notation, we denote by
\begin{gather*}
[1,-1]:=\{1<2<\dots<n<-n<\dots<-2<-1\}
\end{gather*}
the set that we order as indicated. We will also use e.g. $\leq$ on this set, having the evident meaning. Moreover, for $i,j\in [1,-1]$ with $i\leq j$ we set 
\begin{gather*}
[i,j]:=\{z\in[1,-1]\,|\,i\leq z\leq j\}.
\end{gather*}
Note that we will often write $i\in[1,n]$ and then index elements by $-i$ 
which, up to the order, is the same as $i\in[-n,-1]$ and using the index $i$.
\end{Notation}

The following is the main object under study in this paper:

\begin{Definition}\label{D:QExt}
Let $\Lambda_{\intform}=\Lambda_{\intform}(V_{\intform})$ 
(the notation $\Lambda_{\intform}(V_{\intform})$ will be clear later) be the $\intform$-algebra generated by elements
\begin{gather*}
v_{1},v_{2},\dots,v_{n},v_{-n},\dots,v_{-2},v_{-1},
\end{gather*}
modulo the relations
\begin{gather*}
v_{i}^{2}=0,i\in[1,-1]
,\quad
v_{j}v_{i}=-q\cdot v_{i}v_{j}\text{ and }
v_{-i}v_{-j}=-q\cdot v_{-j}v_{-i}\text{ for }i,j\in[1,n],i<j
,
\\
v_{-i}v_{j}=-q\cdot v_{j}v_{-i}\text{ for }i,j\in[1,n],i\neq j
,
\\
v_{-i}v_{i}=-q^{2}\cdot v_{i}v_{-i}
+(q-q^{-1})
\sum_{k\in[1,n-i]}(-q)^{k+1}\cdot v_{i+k}v_{-(i+k)}\text{ for }i\in[1,n]
.
\end{gather*}
We call $\Lambda_{\intform}$ the \emph{quantum exterior algebra of type C}.
\end{Definition}

\begin{Remark}
Note that the final relation in 
\autoref{D:QExt} is very different for the classical counterpart of $\Lambda_{\intform}$ 
due to the factor $(q-q^{-1})$.
\end{Remark}

For $k\in\N$ let further $\Lambda_{\intform}^{k}$ denote 
the $\intform$-span of all monomials in $k$ symbols.

\begin{Lemma}
\changed{The $\intform$-algebra $\Lambda_{\intform}$ is $\N$-graded, and the} set $\Lambda_{\intform}^{k}$ is the 
degree $k$ part of $\Lambda_{\intform}$.
\end{Lemma}

\begin{proof}
The relations are homogeneous with respect to the number of appearing monomials.
\end{proof}

Note that the first three relations in \autoref{D:QExt} are saying that 
$\{v_{i}|i\in[1,n]\}$ and $\{v_{i}|i\in[-n,-1]\}$
span \emph{quantum exterior algebras of type A} defined as in e.g. \cite{BeZw-braided-ext-sym}.

\begin{Remark}
The $\intform$-algebra $\Lambda_{\intform}$ 
morally is a special case of the quantum exterior algebras considered in
\cite{BeZw-braided-ext-sym} where the underlying module is 
the vector representation of $\mathfrak{sp}_{2n}$. By the exterior version of 
\cite[Corollary 4.26]{Zw-r-matrix-poisson} the 
$\intform$-algebra $\Lambda_{\intform}$ is therefore flat, 
i.e. $\Lambda_{\intform}$ is free and has the same rank as 
the dimension of its classical counterpart. We will prove this result in
our setting in \autoref{P:ExtFlat} below.
The symmetric version of the 
$\intform$-algebra $\Lambda_{\intform}$ also appears in 
\cite[Section 5]{LeZhZh-q-first-fundamental-theorem}, and this is where we got the idea 
for the generators and relations of $\Lambda_{\intform}$.
\end{Remark}

The following is an $\intform$-basis of $\Lambda_{\intform}$:

\begin{Definition}\label{D:Standard}
For a(n ordered) subset $S\subset[1,-1]$ with $S=\{x_{1}<x_{2}<\dots<x_{|S|}\}$ 
we define \emph{standard basis vectors} as
\begin{gather*}
v_{S}:=v_{x_{1}}v_{x_{2}}\cdots v_{x_{|S|}}\in\Lambda_{\intform}.
\end{gather*}
For $i\in\{1,\dots,n\}$ let $\epsilon_{i}$ be the standard basis vector of 
the $n$ dimensional Euclidean space.
We define the \emph{weight} $\wt$ on these
$\wt v_{i}=\epsilon_{i}$, $\wt v_{-i}=-\epsilon_{i}$, and
$\wt v_{S}=\sum_{x\in S}\wt v_{x}$.
\end{Definition}

\begin{Example}
\leavevmode
\begin{enumerate}

\item For $n=1$ the generating relations of $\Lambda_{\intform}$ are
\begin{gather*}
v_{1}^{2}=v_{-1}^{2}=0,
\quad
v_{-1}v_{1}=q^{2}\cdot v_{1}v_{-1}.
\end{gather*}
We have $[1,-1]=1<-1$, and the ordered subsets $S\subset[1,-1]$ 
are $S_{\emptyset}=\emptyset$, $S_{1}=\{1\}$, $S_{-1}=\{-1\}$ and $S_{1,-1}=\{1<-1\}$.
The $\intform$-algebra $\Lambda_{\intform}$ is free of rank four and
\begin{gather*}
\{
v_{S_{\emptyset}}=1,
v_{S_{\{1\}}}=v_{1},
v_{S_{\{-1\}}}=v_{-1},
v_{S_{\{1,-1\}}}=v_{1}v_{-1},
\}
\end{gather*}
is a $\intform$-basis of $\Lambda_{\intform}$. The total rank, 
ordered by degree, is $1+2+1=4=2^{2\cdot 1}$.

\item For $n=2$ the most complicated relation is
\begin{gather*}
v_{-1}v_{1}=q^{2}\cdot v_{1}v_{-1}
+(q-q^{-1})(-q)^{2}\cdot v_{2}v_{-2}.
\end{gather*}
We have sixteen subsets $S\subset[1,-1]$ and the standard basis of 
$\Lambda_{\intform}$ is
\begin{gather*}
\{v_{S_{\emptyset}}=1,
v_{S_{\{1\}}}=v_{1},
\dots,
v_{S_{\{1,2,-1\}}}=v_{1}v_{2}v_{-1},
v_{S_{\{1,2,-2,-1\}}}=v_{1}v_{2}v_{-2}v_{-1}\}.
\end{gather*}

\end{enumerate}
Note that for $n=1$ we have that $\mathfrak{sp}_{2\cdot 1}$ and
$\mathfrak{sl}_{2}$ are isomorphic, and this is reflected in $\Lambda_{\intform}$ 
as the most complicated relations do not play any role.
\end{Example}

The standard basis vectors form what we call the \emph{standard basis} of $\Lambda_{\intform}$ given by:

\begin{Proposition}\label{P:ExtFlat}
We have the following.
\begin{enumerate}

\item For all $k\in\{0,1,\dots,2n\}$, the sets
\begin{gather*}
\{
v_{S}|S\subset[1,-1]
\}
,\quad
\{
v_{S}|S\subset[1,-1],|S|=k
\}
\end{gather*}
are $\intform$-basis for $\Lambda_{\intform}$ and $\Lambda_{\intform}^{k}$, respectively.

\item $\Lambda_{\intform}$ is flat, i.e. $\mathrm{rk}_{\intform}\Lambda_{\intform}^{k}=\dim_{\C}\Lambda_{1}^{k}=\binom{2n}{k}$ and thus, $\mathrm{rk}_{\intform}\Lambda_{\intform}=\dim_{\C}\Lambda_{1}=2^{2n}$.

\end{enumerate}

\end{Proposition}

\begin{proof}
\textit{(a).} As it happens quite often in this context, see e.g. 
\cite[Proof of Proposition 2.2]{Su-harmonic-analysis-son} for a closely related usage,
the result is a consequence of Bergman's diamond lemma \cite[Theorem 1.2]{Be-diamond-lemma}. 

Precisely, one easily checks that, if $v_{1}<\dots<v_{n}<v_{-n}<\dots< v_{-1}$, 
then the lexicographic order on monomials satisfies the 
hypotheses of the diamond lemma, and the irreducible monomials 
are exactly the $v_{S}$ such that $S\subset[1,-1]$. 

It remains to verify that one can resolve all the 
following ambiguities coming from the defining relations:
\begin{gather*}
v_{x}v_{x}v_{x},\quad 
v_{-x}v_{-x}v_{-x},\quad 
v_{j}v_{j}v_{i},\quad 
v_{-x}v_{-x}v_{y},\quad 
v_{-i}v_{-i}v_{-j},\quad 
v_{-x}v_{-x}v_{x},\quad
v_{j}v_{i}v_{i},\quad 
v_{j}v_{i}v_{k},
\\
v_{-x}v_{y}v_{y},\quad 
v_{-x}v_{y}v_{z},\quad 
v_{-z}v_{y}v_{z},\quad
v_{-k}v_{-i}v_{-j},\quad
v_{-z}v_{-y}v_{y},\quad
v_{-x}v_{x}v_{x},\quad 
v_{-j}v_{j}v_{i},
\end{gather*}
where $x,y,z\in[1,n]$, such that $z<y$ and $y\neq x \neq z$, and $i,j,k\in [1,n]$, such that $k<i<j$. 
We give one example, showing that $(v_{-z}v_{y})v_{z}=v_{-z}(v_{y}v_{z})$, and leave the rest of the verification as an exercise.

If $X$ denotes a replacement of $v_{-z}v_{y}$ and $X^{\prime}$
one of $v_{y}v_{z}$, then we need to show that $Y$ in
\begin{gather*}
\begin{tikzcd}[ampersand replacement=\&,column sep=0.75em,row sep=0.5em]
\&
v_{-z}v_{y}v_{z}
\ar[dl,very thick,-]
\ar[dr,very thick,-]
\&
\\
Xv_{z}
\ar[dr,very thick,-]
\&
\&
v_{-z}X^{\prime}
\ar[dl,very thick,-]
\\
\&
Y
\&
\end{tikzcd}
\end{gather*}
is the same for both ways to go around the diamond.
To this end, the first calculation is:
\begin{gather*}
(v_{-z}v_{y})v_{z}
=-q\cdot v_{y}v_{-z}v_{z}
=q^{3}\cdot v_{y}v_{z}v_{-z}
+(q-q^{-1})\sum_{r=1}^{n-z}(-q)^{r+2}\cdot v_{y}v_{z+r}v_{-z-r}
\\
=-q^{4}\cdot v_{z}v_{y}v_{-z}
+(q-q^{-1})\sum_{r=1}^{y-z-1}(-q)^{r+3}\cdot v_{z+r}v_{y}v_{-z-r}
+(q-q^{-1})\hspace*{-0.3cm}\sum_{r=y-z+1}^{n-z}(-q)^{r+2}\cdot v_{y}v_{z+r}v_{-z-r}.
\end{gather*}
The second is:
\begin{gather*}
v_{-z}(v_{y}v_{z})=-q\cdot v_{-z}v_{z}v_{y}
=q^{3}\cdot v_{z}v_{-z}v_{y}+(q-q^{-1})\sum_{r=1}^{n-z}(-q)^{r+2}\cdot v_{z+r}v_{-z-r}v_{y}
\\
=-q^{4}\cdot v_{z}v_{y}v_{-z}+
(q-q^{-1})\big(\sum_{\substack{1\leq r\leq n-z \\ z+r\neq y}}(-q)^{r+3}\cdot v_{z+r}v_{y}v_{-z-r} 
+(-q)^{y-z+2}\cdot v_{y}v_{-y}v_{y}\big)
\\
=-q^{4}\cdot v_{z}v_{y}v_{-z}+
(q-q^{-1})\big(\sum_{r=1}^{y-z-1}(-q)^{r+3}\cdot v_{z+r}v_{y}v_{-z-r}
+\sum_{r=y-z+1}^{n-z}(-q)^{r+4}\cdot v_{y}v_{z+r}v_{-z-r}\big)
\\
+(-q)^{y-z+2}(q-q^{-1})^{2}\sum_{s=1}^{n-y}(-q)^{s+1}\cdot v_{y}v_{y+s}v_{-y-s}
\\
=-q^{4}\cdot v_{z}v_{y}v_{-z}
+(q-q^{-1})\big(\sum_{r=1}^{y-z-1}(-q)^{r+3}\cdot v_{z+r}v_{y}v_{-z-r}
+\sum_{r=y-z+1}^{n-z}(-q)^{r+4}\cdot v_{y}v_{z+r}v_{-z-r}\big)
\\
+(q-q^{-1})^{2}\sum_{r=y-z+1}^{n-z}(-q)^{z+r-y+1}(-q)^{y-z+2}\cdot v_{y}v_{z+r}v_{-z-r}
=-q^4v_{z}v_{y}v_{-z}
\\
+(q-q^{-1})\Big(\sum_{r=1}^{y-z-1}(-q)^{r+3}\cdot v_{z+r}v_{y}v_{-z-r}
+\hspace*{-0.2cm}\sum_{r=y-z+1}^{n-z}\big((-q)^{r+4}+(-q)^{r+3}(q-q^{-1})\big)\cdot v_{y}v_{z+r}v_{-z-r}\Big)
\\
=-q^4v_{z}v_{y}v_{-z}
+(q-q^{-1})\big(\sum_{r=1}^{y-z-1}(-q)^{r+3}\cdot v_{z+r}v_{y}v_{-z-r}
+\sum_{r=y-z+1}^{n-z}(-q)^{r+2}\cdot v_{y}v_{z+r}v_{-z-r}\big).
\end{gather*}
Hence, we find that the ambiguity $v_{-z}v_{y}v_{z}$ is resolvable.

\textit{(b).} Directly from (a) since the classical $q=1$ version has
the claimed dimensions.
\end{proof}

Note that, by \autoref{P:ExtFlat}, 
$\Lambda_{\intform}^{k}$ is zero unless $k\in\{0,\dots,2n\}$, so we will only consider those 
degrees.

\subsection{The symplectic action}\label{SS:SymplecticAction}

The $\intform$-algebra $\Lambda_{\intform}$ is actually a module 
of the associated quantum group (we explain this first in the semisimple case only). Before 
we come to the definition we need the following root conventions:

\begin{Notation}\label{N:Uqsp-notation}
All the below are the conventions from \cite[Plate III]{Bo-chapters-4-6}.

The Dynkin diagram that we will use is
$\dynkin C{}$ where we have $n$ nodes labeled, from left to right, $1$ to $n$, and the symmetrizer is 
$(1,\dots,1,2)$.

Let $X:=\oplus_{i=1}^{n}\Z\epsilon_{i}$ (the $\epsilon_{i}$ are as in \autoref{D:Standard}) be the weight lattice for $\mathfrak{sp}_{2n}$, and write $\Phi_{C_{n}}\subset X$ for the root system. We fix the following choice of simple roots:
\begin{gather*}
\alpha_{i}=
\epsilon_{i}-\epsilon_{i+1}, 
\text{ for }i\in\{1,\dots,n-1\},
\text{ and }\alpha_{n}=2\epsilon_{n}.
\end{gather*}
The set of reflections through simple roots $s_{i}:=s_{\alpha_{i}}$ generates the Weyl group $W_{C_{n}}$, which acts on $X$. The form on $X$ is determined by $(\epsilon_{i},\epsilon_{j})= \delta_{i,j}$ and is $W$-invariant. 

Furthermore, we have associated coroots using $\alpha^{\vee}:=2\alpha/(\alpha,\alpha)$, Cartan integers using $a_{ij}:=(\alpha_{i}^{\vee},\alpha_{j})$, and fundamental weights 
by $\varpi_{i}:=\epsilon_{1}+\dots+\epsilon_{i}$  for $i\in\{1,\dots,n\}$.
\end{Notation}

We use the following conventions for quantum groups over $\field$. Integral versions 
will be recalled later in \autoref{S:Weyl}.

The \emph{quantum group} $U_{q}(\mathfrak{sp}_{2n})$ is defined as the $\field$-algebra generated by elements $e_{i}$, $f_{i}$, and $k_{i}^{\pm 1}$, for $i\in\{1,\dots,n\}$, called the \emph{($\mathfrak{sp}_{2n}$) Chevalley generators}. We take everything modulo the relations that $k_{i}$ is the two-sided inverse of $k_{i}^{-1}$, the $k_{i}$ commute among each other and:
\begin{gather*}
k_{i}e_{j}=q^{(\alpha_{i},\alpha_{j})}\cdot e_{j}k_{i}, 
\quad
k_{i}f_{j}=q^{-(\alpha_{i},\alpha_{j})}\cdot f_{j}k_{i},
\quad
e_{i}f_{j}-f_{j}e_{i}=
\delta_{i,j}\cdot\frac{k_{i}-k_{i}^{-1}}{q_{i}-q_{i}^{-1}},
\\
\changed{\sum_{s=0}^{1-a_{ij}}(-1)^{s}
\begin{bmatrix}1-a_{ij} \\ s\end{bmatrix}_{q,i}\cdot
e^{1-a_{ij}-s}e_{j}e_{i}^{s}=0,
\quad
\sum_{s=0}^{1-a_{ij}}(-1)^{s}
\begin{bmatrix}1-a_{ij}\\ s\end{bmatrix}_{q,i}\cdot
f^{1-a_{ij}-s}f_{j}f_{i}^{s}=0.}
\end{gather*}

We choose a Hopf structure on $U_{q}(\mathfrak{sp}_{2n})$ by demanding that the $k_{i}$ are group-like and:
\begin{gather}\label{Eq:CoProductSP}
\begin{gathered}
\Delta(e_{i})=e_{i}\otimes 1+k_{i}\otimes e_{i},
\quad
\Delta(f_{i})=f_{i}\otimes k_{i}^{-1}+1\otimes f_{i},
\\
S(e_{i})=-k_{i}^{-1}e_{i},
\quad
S(f_{i})=-f_{i}k_{i},
\quad 
\epsilon(e_{i})=0,
\quad 
\epsilon(f_{i})=0.
\end{gathered}
\end{gather}
One can check that \autoref{Eq:CoProductSP} defines the structure of a Hopf algebra on $U_{q}(\mathfrak{sp}_{2n})$. 
In fact, this convention is in line with \cite{Ja-lectures-qgroups}.

Recall that there is a $\field$-algebra automorphism, the \emph{Chevalley involution}, determined by
\begin{gather*}
\omega\colon U_{q}(\mathfrak{sp}_{2n})\to U_{q}(\mathfrak{sp}_{2n}),
\quad\omega(e_{i})=f_{i},
\omega(f_{i})=e_{i},
\omega(k_{i}^{\pm 1})=k_{i}^{\mp 1}.  
\end{gather*}
This gives:

\begin{Definition}\label{D:Twist}
Given a $U_{q}(\mathfrak{sp}_{2n})$-module $(V,\rho_{V})$ 
with the indicated action $\rho_{V}$, we define the \emph{$\omega$-twist} 
of $V$, denoted ${^{\omega}V}$, as the $U_{q}(\mathfrak{sp}_{2n})$-module $(V,\rho_{V}\circ\omega)$.
\end{Definition}

\begin{Definition}\label{D:VecRep}
The \emph{vector representation} $V_{q}$ of $U_{q}(\mathfrak{sp}_{2n})$ is the $\field$-vector space
with basis $\{v_{1},v_{2},\dots,v_{n},v_{-n},\dots,v_{-2},v_{-1}\}$ and 
$U_{q}(\mathfrak{sp}_{2n})$-action given by
\begin{gather}\label{Eq:VecRep}
\begin{gathered}
k_{i}\acts v_{\pm j}= 
q_{i}^{\pm\delta_{i,\pm j}}\cdot v_{\pm j},
j\in[-n,-1],
\\
\changed{i\neq n\colon
\begin{cases}
e_{i}\acts v_{i+1}=v_{i}, 
\\
e_{i}\acts v_{-i}=v_{-(i+1)}, 
\\
e_{i}\acts v_{j}=0\text{ if }j\notin\{-i,i+1\},
\end{cases}
\quad
\begin{cases}
f_{i}\acts v_{i}=v_{i+1},
\\
f_{i}\acts v_{-(i+1)}=v_{-i},
\\
f_{i}\acts v_{j}=0\text{ if }j\notin\{i,-(i+1)\},
\end{cases}}
\\
\changed{
\begin{cases}
e_{n}\acts v_{-n}=v_{n},
\\
e_{n}\acts v_{j}=0\text{ else},
\end{cases}
\quad
\begin{cases}
f_{n}\acts v_{n}=v_{-n},
\\
f_{n}\acts v_{j}=0\text{ else}.
\end{cases}}
\end{gathered}
\end{gather}
We call the $\field$-span of $\{v_{1},v_{2},\dots,v_{n}\}$ 
and of $\{v_{-n},\dots,v_{-2},v_{-1}\}$ \emph{vector representations of type A} and its \emph{dual}.
\end{Definition}

\begin{Remark}
The $U_{q}(\mathfrak{sp}_{2n})$-modules such that $K_{i}$ acts on $\mu$ weight vectors as $q^{(\alpha_{i}^{\vee},\mu)}$ are called \emph{type-$1$ modules}.
We will only consider $U_{q}(\mathfrak{sp}_{2n})$-modules of this form, see 
\cite[Section 5.2]{Ja-lectures-qgroups} for background and why this restriction is harmless.
\end{Remark}

\begin{Example}\label{E:VecRepSP}
Let us as usual denote the action of $e_{i}$, $f_{i}$ and $k_{i}$ as labeled edges, with 
the convention that not displayed arrows mean trivial actions.
\begin{gather*}
\begin{tikzcd}[ampersand replacement=\&,column sep=1cm,row sep=1cm]
\colorbox{blue!35}{\mystrut$v_{-1}$}
\ar[r,"e_{1}",yshift=0.075cm,blue]
\ar[loop above,"k_{1}{\colon}q^{-1}",blue]
\&
\colorbox{blue!35}{\mystrut$v_{-2}$}
\ar[r,"e_{2}",yshift=0.075cm,blue]
\ar[l,"f_{1}",yshift=-0.075cm,blue]
\ar[loop above,"k_{2}{\colon}q^{-1};k_{1}{\colon}q",blue]
\&
\colorbox{blue!35}{\mystrut$v_{-3}$}
\ar[r,"e_{3}",yshift=0.075cm,blue]
\ar[l,"f_{2}",yshift=-0.075cm,blue]
\ar[loop above,"k_{3}{\colon}q^{-1};k_{2}{\colon}q",blue]
\&
\colorbox{blue!35}{\mystrut$v_{-4}$}
\ar[r,"e_{4}",yshift=0.075cm,black]
\ar[l,"f_{3}",yshift=-0.075cm,blue]
\ar[loop above,"k_{4}{\colon}q^{-2};k_{3}{\colon}q",blue]
\&
\colorbox{spinach!35}{\mystrut$\,v_{4}\,$}
\ar[r,"e_{3}",yshift=0.075cm,spinach]
\ar[l,"f_{4}",yshift=-0.075cm,black]
\ar[loop above,"k_{4}{\colon}q^{2};k_{3}{\colon}q^{-1}",spinach]
\&
\colorbox{spinach!35}{\mystrut$\,v_{3}\,$}
\ar[r,"e_{2}",yshift=0.075cm,spinach]
\ar[l,"f_{3}",yshift=-0.075cm,spinach]
\ar[loop above,"k_{3}{\colon}q;k_{2}{\colon}q^{-1}",spinach]
\&
\colorbox{spinach!35}{\mystrut$\,v_{2}\,$}
\ar[r,"e_{1}",yshift=0.075cm,spinach]
\ar[l,"f_{2}",yshift=-0.075cm,spinach]
\ar[loop above,"k_{2}{\colon}q;k_{1}{\colon}q^{-1}",spinach]
\&
\colorbox{spinach!35}{\mystrut$\,v_{1}\,$}
\ar[l,"f_{1}",yshift=-0.075cm,spinach]
\ar[loop above,"k_{1}{\colon}q",spinach]
\end{tikzcd}
.
\end{gather*}
The colors illustrate the underlying vector representations of type A.
The Chevalley generators $e_{i}$ and $f_{i}$ always act by one, while the scalar for 
$k_{i}$ is illustrated.
\end{Example}

\begin{Remark}\label{R:ScalingVecRep}
Scaling the standard basis vectors by $v_{-n+1}\mapsto -u_{-n+1}$, $v_{-n+3}\mapsto -u_{-n+3}$ 
all the way to $v_{-1}\mapsto -u_{-1}$ or $v_{-2}\mapsto -u_{-2}$, depending on parity, and $v_{i}\mapsto u_{i}$ otherwise, we see that $V_{q}$ is isomorphic 
to the $U_{q}(\mathfrak{sp}_{2n})$-module illustrated by:
\begin{gather*}
\begin{tikzcd}[ampersand replacement=\&,column sep=1cm,row sep=1cm]
\colorbox{blue!35}{\mystrut$u_{-1}$}
\ar[r,"e_{1}:-1",yshift=0.075cm,blue]
\ar[loop above,"k_{1}{\colon}q^{-1}",blue]
\&
\colorbox{blue!35}{\mystrut$u_{-2}$}
\ar[r,"e_{2}:-1",yshift=0.075cm,blue]
\ar[l,"f_{1}:-1",yshift=-0.075cm,blue]
\ar[loop above,"k_{2}{\colon}q^{-1};k_{1}{\colon}q",blue]
\&
\colorbox{blue!35}{\mystrut$u_{-3}$}
\ar[r,"e_{3}:-1",yshift=0.075cm,blue]
\ar[l,"f_{2}:-1",yshift=-0.075cm,blue]
\ar[loop above,"k_{3}{\colon}q^{-1};k_{2}{\colon}q",blue]
\&
\colorbox{blue!35}{\mystrut$u_{-4}$}
\ar[r,"e_{4}:1",yshift=0.075cm,black]
\ar[l,"f_{3}:-1",yshift=-0.075cm,blue]
\ar[loop above,"k_{4}{\colon}q^{-2};k_{3}{\colon}q",blue]
\&
\colorbox{spinach!35}{\mystrut$\,u_{4}\,$}
\ar[r,"e_{3}:1",yshift=0.075cm,spinach]
\ar[l,"f_{4}:1",yshift=-0.075cm,black]
\ar[loop above,"k_{4}{\colon}q^{2};k_{3}{\colon}q^{-1}",spinach]
\&
\colorbox{spinach!35}{\mystrut$\,u_{3}\,$}
\ar[r,"e_{2}:1",yshift=0.075cm,spinach]
\ar[l,"f_{3}:1",yshift=-0.075cm,spinach]
\ar[loop above,"k_{3}{\colon}q;k_{2}{\colon}q^{-1}",spinach]
\&
\colorbox{spinach!35}{\mystrut$\,u_{2}\,$}
\ar[r,"e_{1}:1",yshift=0.075cm,spinach]
\ar[l,"f_{2}:1",yshift=-0.075cm,spinach]
\ar[loop above,"k_{2}{\colon}q;k_{1}{\colon}q^{-1}",spinach]
\&
\colorbox{spinach!35}{\mystrut$\,u_{1}\,$}
\ar[l,"f_{1}:1",yshift=-0.075cm,spinach]
\ar[loop above,"k_{1}{\colon}q",spinach]
\end{tikzcd}
.
\end{gather*}
This matches the conventions 
one often sees in the literature, 
as for example in \cite[Section 5]{LeZhZh-q-first-fundamental-theorem}, 
and is be a bit better to identify with the classical case.
\end{Remark}

\begin{Lemma}\label{L:ActionExt}
There is an $U_{q}(\mathfrak{sp}_{2n})$-action on $\Lambda_{q}$ determined by \autoref{Eq:VecRep} 
and:
\begin{gather*}
k_{i}\acts v_{I}= 
q^{(\alpha_{i},\wt I)}\cdot v_{i},
\quad
e_{i}\acts v_{x_{1}}v_{x_{2}}\dots v_{x_{m}}= 
\sum_{j=1}^{m}
q^{(\alpha_{i},\wt v_{1}\dots v_{j-1})}\cdot v_{x_{1}}\dots(e_{i}\acts v_{x_{j}})\dots v_{x_{m}},
\\
f_{i}\acts v_{x_{1}}v_{x_{2}}\dots v_{x_{m}}=\sum_{j=1}^{m}q^{(\alpha_{i},\wt v_{x_{j+1}}\dots v_{x_{m}})}\cdot v_{x_{1}}\dots(f_{i}\acts v_{x_{j}})\dots v_{x_{m}}.
\end{gather*}
\end{Lemma}

\begin{proof}
We have that $\Lambda_{q}$ is the tensor algebra in $V_{q}$ modulo the symmetric algebra, that is
$\Lambda_{q}=\Lambda_{q}(V_{q})$.
Thus, the formulas above can be obtained from the coproduct and the claim follows.
\end{proof}

For $S\subset[1,-1]$ we let $S_{i,i+1}:=S\cap\{\pm i,\pm (i+1)\}$ and 
$S_{n}:=S\cap\{\pm n\}$. We need these as these are the only places where the $i$th 
Chevalley generators act nontrivially, see \autoref{Eq:VecRep}.

We can now explicitly describe the $U_{q}(\mathfrak{sp}_{2n})$-action 
from \autoref{L:ActionExt} on the standard basis:

\begin{Lemma}\label{L:ExtAction}
For $i\in\{1,\dots,n-1\}$ and $i=n$ we have:
\begin{align*}
f_{i}\acts v_{S}&= 
\begin{cases}
v_{(S\setminus\{i\})\cup\{i+1\}} & \text{if }S_{i,i+1}=\{i\},
\\
v_{(S\setminus\{-(i+1)\})\cup\{-i\}} & \text{if }S_{i,i+1}=\{-(i+1)\},
\\
v_{(S\setminus\{-(i+1)\})\cup\{-i \}}
+q^{-1}\cdot v_{(S\setminus\{i\})\cup\{i+1\}} & \text{if }S_{i,i+1}=\{i,-(i+1)\},
\\
q\cdot v_{(S\setminus\{i\})\cup\{i+1\}} & \text{if }S_{i,i+1}=\{i,-i\},
\\
v_{(S\setminus\{-(i+1)\})\cup\{-i\}} & \text{if }S_{i,i+1}=\{i+1,-(i+1)\},
\\
v_{(S\setminus\{i\})\cup\{i+1\}} & \text{if }S_{i, i+1}=\{i,-(i+1),-i\},
\\
v_{(S\setminus\{-(i+1)\})\cup\{-i\}} & \text{if }S_{i,i+1}=\{i,(i+1),-(i+1)\},
\\
0 & \text{otherwise},
\end{cases}
\\
f_{n}\acts v_{S}&= 
\begin{cases} 
v_{(S\setminus\{n\})\cup\{-n\}}& \text{if }S_{n}=\{n\},
\\
0 & \text{otherwise}.
\end{cases}
\\
e_{i}\acts v_{S}&= 
\begin{cases}
v_{(S\setminus\{-i)\})\cup\{-(i+1)\}} & \text{if }S_{i,i+1}=\{-i\},
\\
v_{(S\setminus\{i+1\})\cup\{i\}} & \text{if }S_{i,i+1}=\{i+1\},
\\
v_{(S\setminus\{i+1\})\cup\{i\}}
+q^{-1}\cdot v_{S\setminus\{-i\})\cup\{-(i+1)\}} & \text{if }S_{i,i+1}=\{i+1,-i\},
\\
q\cdot v_{(S\setminus\{-i\})\cup\{-(i+1)\}} & \text{if }S_{i,i+1}=\{i,-i\},
\\
v_{(S\setminus\{i+1\})\cup\{i\}} & \text{if }S_{i,i+1}=\{i+1,-(i+1)\},
\\
v_{(S\setminus\{-i\})\cup\{-(i+1)\}} & \text{if }S_{i, i+1}=\{i,i+1,-i\},
\\
v_{(S\setminus\{i+1\})\cup\{i\}} & \text{if }S_{i,i+1}=\{i+1,-(i+1), -i\},
\\
0 & \text{otherwise},
\end{cases}
\\
e_{n}\acts v_{S}&= 
\begin{cases} 
v_{(S\setminus\{-n\})\cup\{n\}}& \text{if }S_{n}=\{-n\},
\\
0 & \text{otherwise}.
\end{cases}
\end{align*}  
\end{Lemma}

\begin{proof}
An annoying but straightforward calculation, so we only give one 
prototypical example. Say $S_{i,i+1}=\{i,-i\}$. Since $f_{i}\acts v_{j}$ is zero unless 
$j=i$ or $j=-(i+1)$, we immediately see that only one summand of $f_{i}\acts v_{x_{1}}v_{x_{2}}\dots v_{x_{m}}$ can survive, namely the one where $f_{i}$ hits $v_{i}$. The vector 
$v_{i}$ is changed to $v_{i+1}$ and thus, $S$ changes to $S\setminus\{i\}\cup\{i+1\}$. 
Moreover, due to our choice of the coproduct $k_{i}^{-1}$ is involved as well and gives the $q$-power.
\end{proof}

\begin{Remark}
Note that the formula for the action of $e_{i}$ is related 
to the formula for the action of $f_{i}$ by ``multiplying $\{i,i+1,-(i+1),-i\}$ by $-1$''. 
In what follows we will often just present the formulas 
for the action of $f_{i}$, noting that the formula for $e_{i}$ is similar.
\end{Remark}

\begin{Lemma}\label{L:Lambda-iso-to-omega-twist}
The map $v_{S}\mapsto v_{-S}$ defines a $U_{q}(\mathfrak{sp}_{2n})$-module isomorphisms $\Lambda_{q}^{k}\to{^{\omega}\Lambda_{q}^{k}}$. 
\end{Lemma}

\begin{proof}
Note that $|S|=|-S|$, so the $\field$-linear map 
determined by $v_{S}\mapsto v_{-S}$ is an
isomorphism $\Lambda^{k}\rightarrow\Lambda^{k}$. Since 
$\wt v_{-S}=\wt v_{S}$, the map clearly 
commutes with the actions of $k_{i}^{\pm 1}$, 
for $i\in\{1,\dots,n\}$. The lemma then follows 
from analyzing the descriptions of actions of 
$e_{i}$ and $f_{i}$, for $i\in\{1,\dots,n\}$ in \autoref{L:ExtAction}.
\end{proof}

\begin{Remark}\label{R:twist-intform-iso}
By \autoref{P:ExtFlat}, the map $v_{S}\mapsto v_{-S}$ also determines an isomorphism of $\intform$-modules $\Lambda_{\intform}^{k}\rightarrow\Lambda_{\intform}^{k}$.
\end{Remark}

The same as in \autoref{L:Lambda-iso-to-omega-twist}, of course, also holds for the whole space
$\Lambda_{q}$.

\subsection{Dot diagrams}

We now rephrase the previous results using certain decorated Young diagrams that we 
call \emph{dot diagrams}.
In this paper all Young diagrams are illustrated in the English convention, and we have 
two rows and $n$ columns:

\begin{Definition}
A \emph{dot diagram} is a two row rectangular Young diagram 
with $2n$ nodes and these nodes can be decorated with zero or one dot.

We view the nodes in the first row as numbered, 
from left to right, $1$ to $n$, and in the second row $-1$ to $-n$. We number the 
columns $1$ to $n$ reading left to right, so that the nodes in the $i$th column 
contain the numbers $i$ and $-i$.
\end{Definition}

\begin{Example}
Here are a few examples of dot diagrams. The leftmost diagram illustrates the numbering convention:
\begin{gather*}
\ytableausetup{centertableaux}
\begin{ytableau}
1 & 2 & 3 \\
{-}1 & {-}2 & {-}3
\end{ytableau}
\,,\quad
\begin{ytableau}
\phantom{a} & *(spinach!50)\bullet & \phantom{a} \\
\phantom{a} & \phantom{a} & \phantom{a}
\end{ytableau}
\,,\quad
\begin{ytableau}
\phantom{a} & *(spinach!50)\bullet & \phantom{a} \\
\phantom{a} & *(spinach!50)\bullet & \phantom{a}
\end{ytableau}
\,,\quad
\begin{ytableau}
\phantom{a} & *(spinach!50)\bullet & \phantom{a} \\
\phantom{a} & *(spinach!50)\bullet & *(spinach!50)\bullet
\end{ytableau}
\,,\quad
\begin{ytableau}
*(spinach!50)\bullet & *(spinach!50)\bullet & \phantom{a} \\
\phantom{a} & *(spinach!50)\bullet & *(spinach!50)\bullet
\end{ytableau}
\,.
\end{gather*}
In these examples $n=3$, which is the number of columns. The color shading is just a visual aid and has no other relevance.
\end{Example}

\begin{Remark}
Various versions of dot diagrams have been used in several papers dealing with symplectic groups.
We got the idea to use dot diagrams for illustrations from \cite{Do-sp-fund}.
\end{Remark}

\begin{Definition}\label{L:DotDiagramsTwo}
We define a \emph{dot map} 
from the standard basis $\{v_{S}|S\subset[1,-1]\}$ to the set of dot diagrams 
with $2n$ nodes by assigning $S=\{x_{1}<x_{2}<\dots<x_{|S|}\}$ to the dot diagram with one dot in each $x_{j}$ position. 
\end{Definition}

\begin{Example}
The dot map is exemplified by:
\begin{gather*}
\ytableausetup{centertableaux}
S=\{1,2,-3,-2\}\mapsto
\begin{ytableau}
*(spinach!50)\bullet & *(spinach!50)\bullet & \phantom{a} \\
\phantom{a} & *(spinach!50)\bullet & *(spinach!50)\bullet
\end{ytableau}
\,.
\end{gather*}
Here $n=3$.
\end{Example}

The following lemma will be used to silently identify dot diagrams with $S\subset[1,-1]$ 
or with $v_{S}$ for $S\subset[1,-1]$.

\begin{Lemma}\label{L:DotDiagrams}
The dot map is a bijection that restricts also to a bijection on the degree $k$ part where the image are dot diagrams with $k$ dots.
\end{Lemma}

\begin{proof}
Immediate.
\end{proof}

Using \autoref{L:DotDiagrams}, we can rephrase the $U_{q}(\mathfrak{sp}_{2n})$-action in \autoref{L:ExtAction} as follows (in all illustrations we have $n=3$):

\begin{enumerate}[itemsep=0.15cm,label=\emph{\upshape(\roman*)}]

\item In general, as soon as two dots end in the same node in the process of the action, then the result is zero. We will therefore only describe the nonzero actions.

\item The Chevalley operators $f_{i}$ for $i\in\{1,\dots,n-1\}$ move dots along the first row rightwards and along the second row leftwards. Here the operator $f_{i}$ only moves dots in the 
$i$th top column and the $(i+1)$th bottom column. For example:
\begin{gather*}
f_{2}\acts
\begin{ytableau}
\phantom{a} & *(spinach!50)\bullet & \phantom{a} \\
\phantom{a} & \phantom{a} & \phantom{a}
\end{ytableau}
=
\begin{ytableau}
\phantom{a} & \phantom{a} & *(spinach!50)\bullet \\
\phantom{a} & \phantom{a} & \phantom{a}
\end{ytableau}
\,,\quad
f_{2}\acts
\begin{ytableau}
\phantom{a} & \phantom{a} & \phantom{a} \\
\phantom{a} & \phantom{a} & *(spinach!50)\bullet
\end{ytableau}
=
\begin{ytableau}
\phantom{a} & \phantom{a} & \phantom{a} \\
\phantom{a} & *(spinach!50)\bullet & \phantom{a}
\end{ytableau}
\,.
\end{gather*}
No other potential dots are moved.

\item The Chevalley operators $e_{i}$ for $i\in\{1,\dots,n-1\}$ do the opposite and move dots leftwards in the first row and rightwards in the second row. 

\item The operators $f_{n}$ and $e_{n}$ act, for example, by
\begin{gather*}
f_{n}\acts
\begin{ytableau}
\phantom{a} & \phantom{a} & *(spinach!50)\bullet \\
\phantom{a} & \phantom{a} & \phantom{a}
\end{ytableau}
=
\begin{ytableau}
\phantom{a} & \phantom{a} & \phantom{a} \\
\phantom{a} & \phantom{a} & *(spinach!50)\bullet
\end{ytableau}
\,,\quad
e_{n}\acts
\begin{ytableau}
\phantom{a} & \phantom{a} & \phantom{a} \\
\phantom{a} & \phantom{a} & *(spinach!50)\bullet
\end{ytableau}
=
\begin{ytableau}
\phantom{a} & \phantom{a} & *(spinach!50)\bullet \\
\phantom{a} & \phantom{a} & \phantom{a}
\end{ytableau}
\,.
\end{gather*}
Again, no other potential dots are moved. The general picture is similar.

\item The $q$ comes in for the following local configurations 
\begin{gather*}
f_{i}\acts
\begin{ytableau}
*(spinach!50)\bullet & \phantom{a} \\
*(spinach!50)\bullet & \phantom{a}
\end{ytableau}
=
q\cdot
\begin{ytableau}
\phantom{a} & *(spinach!50)\bullet \\
*(spinach!50)\bullet & \phantom{a}
\end{ytableau}
\,,\quad
f_{i}\acts
\begin{ytableau}
*(spinach!50)\bullet & \phantom{a} \\
\phantom{a} & *(spinach!50)\bullet
\end{ytableau}
=
\begin{ytableau}
*(spinach!50)\bullet & \phantom{a}\\
*(spinach!50)\bullet & \phantom{a}
\end{ytableau}
+q^{-1}\cdot 
\begin{ytableau}
\phantom{a} & *(spinach!50)\bullet \\
\phantom{a} & *(spinach!50)\bullet
\end{ytableau}
\,.
\end{gather*}
Here we assume that the leftmost displayed column is the $i$th column for $i\in\{1,\dots,n-1\}$.
Similarly for the action of $e_{i}$.

\end{enumerate}

Let us point out that the Chevalley operators 
$e_{i}$ and $f_{i}$ \textit{do not change the number of dots}.
We leave it to the reader to rewrite \autoref{L:ExtAction} in terms of dot diagrams.

\begin{Remark}
We later need to keep track of additional data and to this end we will have \emph{rainbow diagrams}, certain decorated dot diagrams, in \autoref{SS:Rainbow} below.
\end{Remark}

\subsection{Self duality of the exterior algebra}

Let $(\placeholder)^{*}$ denote the $\field$-vector space dual.
Recall that the antipode allows to extend $U_{q}(\mathfrak{sp}_{2n})$-actions to 
these duals.

It is well-known that the $U_{q}(\mathfrak{sp}_{2n})$-module $V_{q}$ and its dual
$V_{q}^{*}$ are isomorphic. The following is an explicit isomorphism:

\begin{Lemma}
There is an isomorphism of $U_{q}(\mathfrak{sp}_{2n})$-modules
\begin{gather*}
\varphi\colon V_{q}\to V_{q}^{*}
,\quad v_{i}\mapsto(-q)^{i-1}\cdot v_{-i}^{*},i\in[1,n]\quad
v_{-i}\mapsto -q^{2n}(-q)^{-(i-1)}\cdot v_{i}^{*},i\in[1,n].
\end{gather*}
\end{Lemma}

\begin{proof}
Essentially by definition, the map $\varphi$ is a $\field$-linear bijection.
Moreover, $U_{q}(\mathfrak{sp}_{2n})$-equivariance is checked by direct calculation. 
For example, if $i\in[1,n-1]$, then
\begin{gather*}
\varphi(f_{i}\acts v_{i})=
\varphi(v_{i+1})=(-q)^{i}\cdot v_{-(i+1)}^{*}
=(-q)^{i-1}\cdot f_{i}\acts v_{-i}^{*}
=f_{i}\acts\varphi(v_{i}),
\\
\varphi(f_{n}\acts v_{n})=\varphi(v_{-n})=-q^{2n}(-q)^{-(n-1)}\cdot v_{-n}^{*}
=q(-q)^{n}\cdot v_{-n}^{*}
\\
=-(-q)^{n-1}q^{2}\cdot v_{-n}^{*}
=(-q)^{n-1}\cdot f_{n}\acts v_{n}^{*}
=f_{n}\acts\varphi(v_{n}),
\end{gather*}
where we used that $S(f_{i})=-f_{i}k_{i}$ (this holds also 
for $i=n$) and $k_{i}\acts v_{i}=q\cdot v_{i}$, as well 
as $k_{n}\acts v_{n}=q^{2}\cdot v_{n}$.
\end{proof}

For $i\in[1,n]$ and $S=\{x_{1},x_{2},\dots,x_{m}\}\subset[1,-1]$ such that 
$x_{1}<x_{2}<\dots<x_{m}$, we define:
\begin{gather*}
d_{i}:=\varphi^{-1}(v_{-i}^{*})=(-q)^{-(i-1)}\cdot v_{i},
\quad
d_{-i}:=\varphi^{-1}(v_{i}^{*})=q^{-2n}(-q)^{i-1}\cdot v_{-i},
\\
d_{S}:=d_{x_{1}}d_{x_{2}}\dots d_{x_{m}}. 
\end{gather*}
These elements will be used in \autoref{P:SelfDual} below
to describe the \emph{dual quantum exterior algebra of type C}.

\begin{Lemma}\label{L:dualext}
The $\intform$-algebra 
$\Lambda_{\intform}$ is generated by the elements 
\begin{gather*}
d_{1},\dots,d_{n},d_{-n},\dots,d_{-1},
\end{gather*}
subject only to the following relations:
\begin{gather*}
d_{i}^{2}=0,i\in[1,-1],
\quad
d_{j}d_{i}=-q^{-1}\cdot d_{i}d_{j}\text{ and }
d_{-i}d_{-j}=-q^{-1}\cdot d_{-j}d_{-i}\text{ for }i,j\in[1,n],i<j
\\
d_{-i}d_{j}=-q^{-1}d_{j}d_{-i}\text{ for }i,\in[1,n],i\neq j
\\
d_{-i}d_{i}=-q^{-2}\cdot d_{i}d_{-i}+ 
(q^{-1}-q)\sum_{k\in [1,n-i]}(-q^{-1})^{k+1}\cdot d_{i+k}d_{-(i+k)}\text{ for }i\in[1,n].
\end{gather*}

Moreover, for all $k\in\{0,1,\dots,2n\}$, the sets
\begin{gather*}
\{
d_{S}|S\subset[1,-1]
\}
,\quad
\{
d_{S}|S\subset[1,-1],|S|=k
\}
\end{gather*}
are $\intform$-bases for $\Lambda_{\intform}$ and $\Lambda_{\intform}^{k}$, respectively.
\end{Lemma}

\begin{proof}
Up to units in $\intform$, the set of generators $\{v_{i}\mid i\in [1,-1]\}$ is equal to the set of elements $\{d_{i}\mid i\in [1,-1]\}$. So the latter generates $\Lambda_{\intform}$. That the given sets are $\intform$-bases can be proven as in \autoref{P:ExtFlat}.
\end{proof}

\begin{Lemma}\label{L:DualOne}
\leavevmode

\begin{enumerate}

\item For $i\in\{1,\dots,n-1\}$ and $i=n$ we have:
\begin{align*}
f_{i}\acts d_{S}&= 
\begin{cases}
(-q)\cdot d_{(S\setminus\{i\})\cup\{i+1\}} & \text{if }S_{i,i+1}=\{i\},
\\
(-q)\cdot d_{(S\setminus\{-(i+1)\})\cup\{-i\}}& \text{if }S_{i,i+1}=\{-(i+1)\}, 
\\
(-q)\cdot\big(d_{(S\setminus\{-(i+1)\})\cup\{-i\}} 
+q^{-1}\cdot d_{(S\setminus\{i\})\cup\{i+1\}}\big) & \text{if }S_{i,i+1}=\{i,-(i+1)\},
\\
(-q)\cdot\big(q\cdot d_{(S\setminus\{i\})\cup\{i+1\}}\big) & \text{if }S_{i,i+1}=\{i,-i\},
\\
(-q)\cdot d_{(S\setminus\{-(i+1)\})\cup\{-i\}} & \text{if }S_{i,i+1}=\{i+1,-(i+1)\},
\\
(-q)\cdot d_{(S\setminus\{i\})\cup\{i+1\}} & \text{if }S_{i,i+1}=\{i,-(i+1),-i\},
\\
(-q)\cdot d_{(S\setminus\{-(i+1)\})\cup\{-i\}} & \text{if }S_{i,i+1}=\{i,(i+1),-(i+1)\},
\\
0 & \text{otherwise}.
\end{cases}
\\
f_{n}\acts d_{S}&= 
\begin{cases} 
(-q^{2})\cdot d_{(S\setminus\{n\})\cup\{-n\}}& \text{if }S_{n}=\{ n\},
\\
0 & \text{otherwise}.
\end{cases}
\end{align*}
Similarly, for $e_{i}\acts d_{S}$ and $e_{n}\acts d_{S}$. 

\item For $i\in\{1,\dots,n-1\}$ and $i=n$ we have:
\begin{align*}
f_{i}\acts v_{-S}^{*}&= 
\begin{cases}
(-q)\cdot v^{*}_{-((S\setminus\{i\})\cup\{ i+1\})} & \text{if }S_{i,i+1}=\{i\}, 
\\
(-q)\cdot v^{*}_{-((S\setminus\{-(i+1)\})\cup\{-i\})} & \text{if }S_{i,i+1}=\{-(i+1)\},
\\
\scalebox{0.93}{$(-q)\cdot\big(v^{*}_{-((S\setminus\{-(i+1)\})\cup\{-i\})}
+q^{-1}\cdot v^{*}_{-((S\setminus\{i\})\cup\{i+1\})}\big)$} & \text{if }S_{i, i+1}=\{i,-(i+1)\},
\\
(-q)\cdot\big(q\cdot v^{*}_{-((S\setminus\{i\})\cup\{i+1\})}\big) 
& \text{if }S_{i,i+1}=\{i,-i\},
\\
(-q)\cdot v^{*}_{-((S\setminus\{-(i+1)\})\cup\{ -i \})} & \text{if }S_{i,i+1}=\{i+1,-(i+1)\},
\\
(-q)\cdot v^{*}_{-((S\setminus\{i\})\cup\{i+1\})} & \text{if }S_{i,i+1}= \{i,-(i+1),-i\},
\\
(-q)\cdot v^{*}_{-((S\setminus\{-(i+1)\})\cup\{-i\})} & \text{if }S_{i,i+1}=\{i,(i+1),-(i+1)\},
\\
0 & \text{otherwise}.
\end{cases}
\\
f_{n}\acts v_{-S}^{*}&= 
\begin{cases} 
(-q^{2})\cdot v^{*}_{-((S\setminus\{n\})\cup\{-n\})} & \text{if }S_{n}=\{n\},
\\
0 & \text{otherwise}.
\end{cases}
\end{align*}
Similarly, for $e_{i}\acts v_{-S}^{*}$ and $e_{n}\acts v_{-S}^{*}$.

\end{enumerate}
See also \autoref{L:ExtAction} for the action of $e_{i}$: 
the comparison between the actions of $f_{i}$ 
and $e_{i}$ therein is similar to the ones in these formulas.
\end{Lemma}

\begin{proof}
As in \autoref{L:ExtAction}, and omitted.
\end{proof}

We leave it to the reader to formulate the dot diagram analog of \autoref{L:DualOne}.

\begin{Proposition}\label{P:SelfDual}
The map
\begin{align*}
\varphi^{k}_{q}\colon\Lambda_{q}^{k}\to(\Lambda_{q}^{k})^{*},
\quad d_{S}\mapsto v_{-S}^{*}
\end{align*}
is an isomorphism of $U_{q}(\mathfrak{sp}_{2n})$-modules.
\end{Proposition}

\begin{proof}
\autoref{P:ExtFlat} shows that $\Lambda^{k}_{\intform}$ is a free $\intform$-module with 
basis $\{v_{S}|S\subset[1,-1],|S|=k\}$. Therefore, $\Hom_{\intform}(\Lambda^{k}_{\intform},\intform)$ is also free with basis $\{v_{S}^{*}|S\subset[1,-1],|S|=k\}$. \autoref{L:dualext} implies that the map $\varphi^{k}_{\intform}\colon\Lambda^{k}_{\intform}\to\Hom_{\intform}(\Lambda^{k}_{\intform}, \intform)$, defined by $d_{S}\mapsto v_{-S}^{*}$, is an $\intform$-module isomorphism. Tensoring with $\field$ yields a $\field$-vector space isomorphism $\varphi_{q}^{k}\colon\Lambda_{q}^{k}\to(\Lambda_{q}^{k})^{*}$.

To see $U_{q}(\mathfrak{sp}_{2n})$-equivariance one uses 
(a) and (b) of \autoref{L:DualOne} to compare 
the actions of the Chevalley generators on the bases 
given by $d_{S}$ and $v_{-S}^{*}$, respectively.
\end{proof}

\subsection{The special linear action}\label{SS:SLAction}

We now discuss the dual action. 
In the classical case and for 
$\C$ as the ground field, the formula $2^{2n}=4^{n}$ shows that, 
as $\C$-vector spaces, we have
\begin{gather*}
\Lambda_{1}(V_{1})\cong\Lambda_{1}(\C^{2n})\cong\big(\Lambda_{1}(\C^{2})\big)^{\otimes n}
\cong(\C\oplus\C^{2}\oplus\C)^{\otimes n}.
\end{gather*}
There is thus a rather straightforward $U_{1}(\mathfrak{sl}_{2})$-action on $\Lambda_{1}(V_{1})$, 
simply by acting on $\C\oplus\C^{2}\oplus\C$ (we read this 
as the trivial plus vector plus 
trivial $U_{1}(\mathfrak{sl}_{2})$-module) and then taking the $n$th tensor power.
The same works in the classical case in general as pointed out in 
\cite[Section 2.2]{Fo-symplectic-simple}.
We now describe the quantum version of this action.

\begin{Remark}
The classical version of the below is not the same as to Howe's approach from 
\cite[Theorem 3.8.9.3]{Ho-perspectives-invariant-theory}
where differential operators play a key role.
We will describe the quantum analog of differential operators later in \autoref{SS:DualAction}.
\end{Remark}

Recall that 
the \emph{quantum group} $U_{q}(\mathfrak{sl}_{2})$ 
is the $\field$-algebra generated by $E$, $F$ and $K^{\pm 1}$, 
the \emph{($\mathfrak{sl}_{2}$) Chevalley generators}, such that 
$K$ and $K^{-1}$ are two-sided inverses and
\begin{gather*}
KE=q^{2}EK,\quad KF=q^{-2}FK,\quad EF-FE=\frac{K-K^{-1}}{q-q^{-1}}.
\end{gather*}

We will additionally choose the following coproduct $\Delta$, antipode 
$S$ and counit $\epsilon$ on $U_{q}(\mathfrak{sl}_{2})$.
First, $K$ and $K^{-1}$ are group like, and second:
\begin{gather}\label{Eq:CoProductSL}
\begin{gathered}
\Delta(E)=1\otimes E+E\otimes K,
\quad
\Delta(F)=F\otimes 1+K^{-1}\otimes F,
\\
S(F)=-K^{-1}F,
\quad
S(E)=EK,
\quad
\epsilon(E)=0,
\quad
\epsilon(F)=0.
\end{gathered}
\end{gather}
One easily checks that this defines the structure of a Hopf algebra on 
$U_{q}(\mathfrak{sl}_{2})$. We will use this structure throughout.

\begin{Remark}
Comparing \autoref{Eq:CoProductSP} and \autoref{Eq:CoProductSL},
note that we use the opposite coproduct conventions 
for $U_{q}(\mathfrak{sp}_{2n})$ and for $U_{q}(\mathfrak{sl}_{2})$.
\end{Remark}

A $U_{q}(\mathfrak{sl}_{2})$-module is a \emph{weight module}, 
if the action of $K$ is diagonalizable.  
Given a finite dimensional weight $U_{q}(\mathfrak{sl}_{2})$-module, the eigenvalues of $K$ are contained in the set $\{\pm q^{n}\}_{n\in\mathbb{Z}}$, see e.g. 
\cite[Proposition 2.3]{Ja-lectures-qgroups}. We will only consider 
weight modules for $U_{q}(\mathfrak{sl}_{2})$ and will drop that adjective.

\begin{Definition}
Following for example \cite{Ja-lectures-qgroups}, if all of 
the eigenvalues of a $U_{q}(\mathfrak{sl}_{2})$-module 
are \changed{positive powers of $q$}, then we say the module 
is \emph{type-1} and if all of the eigenvalues are 
\changed{negative powers of $q$}, then we say the module is \emph{type-(-1)}.
\end{Definition}

\begin{Definition}
Let $\sigma$ be the $\field$-algebra automorphism of $U_{q}(\mathfrak{sl}_2)$ determined by
\begin{gather*}
\sigma(E)=-E,
\quad 
\sigma(F)=F,
\quad
\sigma(K^{\pm 1})=-K^{\pm 1}. 
\end{gather*}
Given a $U_{q}(\mathfrak{sl}_2)$-module $(V,\rho_{V})$ with action $\rho_{V}$, we define the \emph{$\sigma$-twist of $V$}, denoted by ${^{\sigma}V}$, as the $U_{q}(\mathfrak{sl}_2)$-module
$(V,\rho_{V}\circ\sigma)$.
\end{Definition}

As one easily checks, if $V$ is a type-$\epsilon$ $U_{q}(\mathfrak{sl}_2)$-module, then ${^{\sigma}V}$ is a type-$(-\epsilon)$ $U_{q}(\mathfrak{sl}_2)$-module.

\begin{Lemma}\label{L:twist-sigma-equivalence}
The category of finite dimensional $U_{q}(\mathfrak{sl}_2)$-modules is a direct sum of the category of type-1 $U_{q}(\mathfrak{sl}_2)$-modules and the category of type-$(-1)$ $U_{q}(\mathfrak{sl}_2)$-modules. Moreover, twisting by $\sigma$ induces an equivalence of categories between the two summands. 
\end{Lemma}

\begin{proof}
See \cite[Section 5.2]{Ja-lectures-qgroups}.
\end{proof}

We will use the following two $U_{q}(\mathfrak{sl}_{2})$-modules, 
one will be used as a type-1 module, the other as a type-(-1) module:

\begin{Definition}\label{D:VecRepSL2}
Let $\epsilon\in\{\pm 1\}$. The \emph{type-$\epsilon$ 
trivial representation} ${^{\epsilon}\F}_{\ast}$
of $U_{q}(\mathfrak{sl}_{2})$ is the $\field$-vector space
with basis $\{w_{0}\}$ and 
$U_{q}(\mathfrak{sl}_{2})$-action given by
\begin{gather*}
E\acts w_{0}=F\acts w_{0}=0,\quad
K\acts w_{0}=\epsilon\cdot w_{0}.
\end{gather*}
Moreover, the \emph{type-$\epsilon$
vector representation} 
${^{\epsilon}W}_{q}$ of $U_{q}(\mathfrak{sl}_{2})$ is the $\field$-vector space
with basis $\{w_{1},w_{-1}\}$ and 
$U_{q}(\mathfrak{sl}_{2})$-action given by
\begin{gather*}
\begin{tikzcd}[ampersand replacement=\&,column sep=1cm,row sep=1cm]
\colorbox{blue!35}{\mystrut$w_{-1}$}
\ar[r,"E",yshift=0.075cm,black]
\ar[loop above,"K{\colon}\epsilon q^{-1}",blue]
\&
\colorbox{spinach!35}{\mystrut$\,w_{1}\,$}
\ar[l,"F",yshift=-0.075cm,black]
\ar[loop above,"K{\colon}\epsilon q",spinach]
\end{tikzcd}
\end{gather*}
The picture is to be read as in \autoref{E:VecRepSP}.
\end{Definition} 

Consider the $U_{q}(\mathfrak{sl}_{2})$-module
\begin{gather*}
X_{q}={^{1}\mathbb{F}}_{\ast}\oplus{^{-1}W}_{q}\oplus{^{1}\mathbb{F}}_{\ast}
=\field\{w_{0},w_{-1},w_{1},w_{0}^{\prime}\},
\end{gather*}
with the indicated $\field$-basis.
The $n$th tensor power of $X_{q}$ has an $\field$-basis
\begin{gather*}
\changed{\big\{\vec{\xsym}=\xsym_{1}\otimes\xsym_{2}\otimes\dots\otimes\xsym_{n}|
\xsym_{i}\in\{w_{0},w_{-1},w_{1},w_{0}^{\prime}\}\big\}.}
\end{gather*}

\begin{Definition}\label{D:TensorToDot}
\changed{For each $\vec{\xsym}$ we associate a set 
$S=S_{\vec{\xsym}}\subset[1,-1]$} as follows:
\begin{gather*}
S\cap\{\pm i\}=
\left\{
\changed{\begin{aligned}
i & \text{ if $\changed{\xsym_{i}}=w_{0}$},
\\
-i & \text{ if $\changed{\xsym_{i}}=w_{0}^{\prime}$}, 
\\
\pm i & \text{ if $\changed{\xsym_{i}}=w_{1}$},
\\
\emptyset & \text{ if $\changed{\xsym_{i}}=w_{-1}$}.
\end{aligned}}
\right.
\end{gather*}
(Note that this determines $S$.)
We use this to define an $\field$-linear map
$\Psi\colon
X_{q}^{\otimes n}\to\Lambda_{q}$ by 
\changed{$\vec{\xsym}\mapsto v_{S_{\vec{\xsym}}}$}.
\end{Definition}

In dot diagrams \autoref{D:TensorToDot} means that the $i$th column 
of $S$ is as follows:
\begin{gather}\label{Eq:TensorToDot}
\changed{\xsym_{i}}=w_{0}
\leftrightsquigarrow
\begin{ytableau}
*(spinach!50)\bullet \\
\phantom{a}
\end{ytableau}
\,,\quad
\changed{\xsym_{i}}=w_{0}^{\prime}
\leftrightsquigarrow
\begin{ytableau}
\phantom{a} \\
*(spinach!50)\bullet
\end{ytableau}
\,,\quad
\changed{\xsym_{i}}=w_{1}
\leftrightsquigarrow
\begin{ytableau}
*(spinach!50)\bullet \\
*(spinach!50)\bullet
\end{ytableau}
\,,\quad
\changed{\xsym_{i}}=w_{-1}
\leftrightsquigarrow
\begin{ytableau}
\phantom{a} \\
\phantom{a}
\end{ytableau}
\,.
\end{gather}
These are local pictures showing only the $i$th column.
The two column configurations on the 
right are called \emph{fully dotted} and \emph{undotted}.

\begin{Example}
For $n=4$ we have
\begin{gather*}
w_{1}\otimes w_{-1}\otimes w_{0}\otimes w_{0}^{\prime}
\leftrightsquigarrow
S=\{1,3,-4,-1\}
\leftrightsquigarrow
\begin{ytableau}
*(spinach!50)\bullet & \phantom{a} & *(spinach!50)\bullet & \phantom{a}\\
*(spinach!50)\bullet & \phantom{a} & \phantom{a} & *(spinach!50)\bullet
\end{ytableau}
\,.
\end{gather*}
The leftmost column is fully dotted, the second column from the left is undotted.
\end{Example}

We use the following lemma to identify the tensor product basis 
coming from $X_{q}^{\otimes n}$ with dot diagrams.

\begin{Lemma}
The map $\Psi$ is an isomorphism 
of $\field$-vector spaces with inverse $\Psi^{-1}(v_{S})=\changed{\vec{\xsym}}$, where 
\begin{gather*}
\changed{\xsym_{i}}:=
\left\{
\changed{\begin{aligned}
w_{0} & \text{ if $S\cap\{\pm i\}=i$},
\\
w_{0}^{\prime} & \text{ if $S\cap\{\pm i\}=-i$}, 
\\
w_{-1} & \text{ if $S\cap\{\pm i\}=\emptyset$},
\\
w_{1} & \text{ if $S\cap\{\pm i\}=\pm i$}.
\end{aligned}}
\right.
\end{gather*}

\end{Lemma}

\begin{proof}
A direct verification.
\end{proof}

The previous lemma endows $\Lambda_{q}$ with an action of $U_{q}(\mathfrak{sl}_{2})$:
\begin{gather*}
X\acts v_{S}:=\Psi\big(X\acts\Psi^{-1}(v_{S})\big)
\text{ for all $X\in U_{q}(\mathfrak{sl}_{2})$}. 
\end{gather*}
We will use this action throughout.

\begin{Notation}\label{N:Subsets}
Let $S\subset[1,-1]$. Define $S_{0}:=S\cap -S$, i.e. the subset of $S$ consisting of all fully dotted 
columns. Also, write $S^{c}=[1,-1]\setminus S$ and $S_{0}^{c}:=S^{c}\cap -S^{c}$, i.e. the subset of $[1,-1]$ consisting of all undotted columns. 

For $k\in\{1,\dots,n\}$, and $T\subset[1,-1]$, write $T_{>k}:=T\cap\{\pm (k+1),\pm (k+2),\dots,\pm n\}$. Similarly define $T_{\geq k}$, $T_{<k}$, and $T_{\leq k}$. We set $\textbf{w}(T):=1/2\big(|T_{0}|-|T_{0}^{c}|\big)$,
$\textbf{w}_{>i}(T):=1/2\big(|(T_{0})_{>i}|-|(T_{0}^{c})_{>i}|\big)$
and $\textbf{w}_{<i}(T):=1/2\big(-|(T_{0})_{<i}|+|(T_{0}^{c})_{<i}|\big)$.
\end{Notation}

\begin{Example}\label{E:DotStatistic}
Let $n=6$ and $S=\{1,2,5,-5,-4\}$. Then $-S=\{4,5,-5,-2,-1\}$, $S_{0}=\{5,-5\}$,
$S^{c}=\{3,4,6,-6,-3,-2,-1\}$, $-S^{c}=\{1,2,6,-6,-4,-3\}$, and $S_{0}^{c}=\{6,-6\}$.
In dot diagrams:
\begin{gather*}
\begin{array}{r@{\ }c@{\ }l}
S
\leftrightsquigarrow
\begin{ytableau}
*(spinach!50)\bullet & *(spinach!50)\bullet & \phantom{a} & \phantom{a} & *(spinach!50)\bullet & \phantom{a} \\
\phantom{a} & \phantom{a} & \phantom{a} & *(spinach!50)\bullet & *(spinach!50)\bullet & \phantom{a}
\end{ytableau}
\,,&\quad
S^{c}
&\leftrightsquigarrow
\begin{ytableau}
\phantom{a} & \phantom{a} & *(spinach!50)\bullet & *(spinach!50)\bullet & \phantom{a} & *(spinach!50)\bullet \\
*(spinach!50)\bullet & *(spinach!50)\bullet & *(spinach!50)\bullet & \phantom{a} & \phantom{a} & *(spinach!50)\bullet
\end{ytableau}
\,,
\\[0.5cm]
-S
\leftrightsquigarrow
\begin{ytableau}
\phantom{a} & \phantom{a} & \phantom{a} & *(spinach!50)\bullet & *(spinach!50)\bullet & \phantom{a} \\
*(spinach!50)\bullet & *(spinach!50)\bullet & \phantom{a} & \phantom{a} & *(spinach!50)\bullet & \phantom{a}
\end{ytableau}
\,,&\quad
-S^{c}
&\leftrightsquigarrow
\begin{ytableau}
*(spinach!50)\bullet & *(spinach!50)\bullet & *(spinach!50)\bullet & \phantom{a} & \phantom{a} & *(spinach!50)\bullet \\
\phantom{a} & \phantom{a} & *(spinach!50)\bullet & *(spinach!50)\bullet & \phantom{a} & *(spinach!50)\bullet
\end{ytableau}
\,,
\\[0.5cm]
S_{0}
\leftrightsquigarrow
\begin{ytableau}
\phantom{a} & \phantom{a} & \phantom{a} & \phantom{a} & *(spinach!50)\bullet & \phantom{a} \\
\phantom{a} & \phantom{a} & \phantom{a} & \phantom{a} & *(spinach!50)\bullet & \phantom{a}
\end{ytableau}
\,,&\quad
S_{0}^{c}
&\leftrightsquigarrow
\begin{ytableau}
\phantom{a} & \phantom{a} & *(spinach!50)\bullet & \phantom{a} & \phantom{a} & *(spinach!50)\bullet \\
\phantom{a} & \phantom{a} & *(spinach!50)\bullet & \phantom{a} & \phantom{a} & *(spinach!50)\bullet
\end{ytableau}
\,.
\end{array}
\end{gather*}
In words, the minus operation flips the dot diagram upside down, the complement fills in the open nodes with dots and removes all other dots, and the subscript zero is the part stable under mirroring upside down.
\end{Example}

\begin{Example}
Let $n=6$ and let 
\begin{gather*}
S=\{1,3,-3,-5\}
\leftrightsquigarrow
\begin{ytableau}
*(spinach!50)\bullet & *(blue!50) & *(spinach!50)\bullet & *(blue!50) & \phantom{a} & *(blue!50) \\
\phantom{a} & *(blue!50) & *(spinach!50)\bullet & *(blue!50) & *(spinach!50)\bullet & *(blue!50)
\end{ytableau}
\,,
\end{gather*}
where we shaded the undotted columns of 
the dot diagram for $S$. Then
\begin{gather*}
(S_{0}^{c})_{>3}=\{4,-4,6,-6\}
\leftrightsquigarrow
\begin{ytableau}
\phantom{a} & \phantom{a} & \phantom{a} & *(spinach!50)\bullet & \phantom{a} & *(spinach!50)\bullet \\
\phantom{a} & \phantom{a} & \phantom{a} & *(spinach!50)\bullet & \phantom{a} & *(spinach!50)\bullet
\end{ytableau}
\,.
\end{gather*}
Note that this fills in exactly the shaded parts above that are to the right of column $3$.
\end{Example}

We now describe the action of $U_{q}(\mathfrak{sl}_{2})$ on $\Lambda_{q}$.
By \autoref{Eq:TensorToDot} and up to $q$-powers, $E$ acts on 
a dot diagram by adding two dots into undotted columns while $F$ acts by removing 
the two dots in fully dotted columns. In particular, $E$ and $F$ do not fix the number of 
dots, but only \textit{preserve the number of dots modulo two}.

\begin{Example}\label{E:TensorAction}
As an example for the action of $F$ for $n=4$ we have 
\begin{gather*}
\begin{aligned}
F\acts w_{1}\otimes w_{-1}\otimes w_{1}\otimes w_{0}^{\prime}
&=
(F\acts w_{1})\otimes w_{-1}\otimes w_{1}\otimes w_{0}^{\prime}
+
\scalebox{0.92}{$(K^{-1}\acts w_{1})\otimes (K^{-1}\acts w_{-1})\otimes (F\acts w_{1})\otimes w_{0}^{\prime}$}
\\
&=
w_{-1}\otimes w_{-1}\otimes w_{1}\otimes w_{0}^{\prime}
+
(-q^{-1})(-q)\cdot
w_{1}\otimes w_{-1}\otimes w_{-1}\otimes w_{0}^{\prime}
\end{aligned}
\\
F\acts
\begin{ytableau}
*(spinach!50)\bullet & \phantom{a} & *(spinach!50)\bullet & \phantom{a}\\
*(spinach!50)\bullet & \phantom{a} & *(spinach!50)\bullet & *(spinach!50)\bullet
\end{ytableau}
=
\begin{ytableau}
\phantom{a} & \phantom{a} & *(spinach!50)\bullet & \phantom{a}\\
\phantom{a} & \phantom{a} & *(spinach!50)\bullet & *(spinach!50)\bullet
\end{ytableau}
+
\begin{ytableau}
*(spinach!50)\bullet & \phantom{a} & \phantom{a} & \phantom{a}\\
*(spinach!50)\bullet & \phantom{a} & \phantom{a} & *(spinach!50)\bullet
\end{ytableau}
\,,
\end{gather*}
These two calculations are the same, but using vectors and dot diagrams, respectively.
\end{Example}

\autoref{E:TensorAction} formalized gives:

\begin{Lemma}\label{L:EFops}
The action of $U_{q}(\mathfrak{sl}_{2})$ on $\Lambda_{q}$ has the 
following explicit description. 
\begin{gather*}
\changed{K\acts v_{S}=
(-q)^{\mathbf{w}(S)}\cdot v_{S},}
\\
\changed{E\acts v_{S}
=\sum_{\{\pm i\}\subset S_{0}^{c}}(-q)^{\mathbf{w}_{> i}(S)}
\cdot
v_{S\cup\{\pm i\}},
\quad
F\acts v_{S}=\sum_{\{\pm i\}\subset S_{0}}(-q)^{\mathbf{w}_{< i}(S)}
\cdot
v_{S\setminus\{\pm i\}}.}
\end{gather*}
\end{Lemma}

\begin{proof}
Directly from \autoref{Eq:CoProductSL} and the translation of $\vec{w}$ to dot diagrams.
\end{proof}

\begin{Lemma}
The $(-q)^{k}$ weight space of the $U_{q}(\mathfrak{sl}_{2})$-module
$\Lambda_{q}$ is exactly $\Lambda_{q}^{n+k}$.
\end{Lemma}

\begin{proof}
Note that $|S|=|S_{0}|+x$, where $x$ is the number 
of pairs $\pm i$ such that $|S\cap\{\pm i\}|=1$. 
Similarly, $|S^{c}|=|S_{0}^{c}|+x$. Therefore, 
$2n=|S|+|S^{c}|=|S_{0}|+2x+|S_{0}^{c}|$,
so $|S_{0}|-|S_{0}^{c}|/2
=|S_{0}|/2-(n-x-|S_{0}|/2)=-n+|S|$.
\end{proof}

\begin{Proposition}\label{P:Commute}
The two actions, the one of $U_{q}(\mathfrak{sp}_{2n})$ and the one of ${U}_{q}(\mathfrak{sl}_{2})$, on $\Lambda_{q}$ commute.
\end{Proposition}

\begin{proof}
It suffices to show that if $S\subset[1,-1]$ and 
$k\in\{1,\dots,n\}$, then $x_{i}\acts(X\acts v_{S})=X\acts(x_{i}\acts v_{S})$, 
for $x_{i}\in\{e_{i},f_{i}\}$ and $X\in\{E,F\}$. Using 
\autoref{L:ExtAction} and \autoref{L:EFops} this 
is a series of elementary calculations. We give for 
following example, and leave the rest to the reader. 

Suppose $S_{k,k+1}=\{k,-(k+1)\}$ and write 
$T:=(S\setminus\{-(k+1)\})\cup\{-k\}$ and $U:=(S\setminus\{k\})\cup\{k+1\}$. Then
\begin{align*}
E\acts (f_{k}\acts v_{S}) &= E\acts (v_{T} + q^{-1}\cdot v_{U}) \\
&=\sum_{\{\pm i\}\subset T_{0}^{c}}(-q)^{\textbf{w}_{>i}(T)}
\cdot
v_{T\cup\{\pm i\}} + q^{-1}\sum_{\{\pm i\}\subset U_{0}^{c}}(-q)^{\textbf{w}_{>i}(U)}
\cdot
v_{U\cup\{\pm i\}} \\
&=\sum_{\substack{\{\pm i\}\subset T_{0}^{c} \\ i \neq k+1}}(-q)^{\textbf{w}_{>i}(T)}
\cdot
v_{T\cup\{\pm i\}} + (-q)^{\textbf{w}_{>k+1}(T)}
\cdot
v_{T\cup\{\pm (k+1)\}}\\
&+q^{-1}(-q)^{\textbf{w}_{>k}(U)}
\cdot
v_{U\cup\{\pm k\}}+ q^{-1}\sum_{\substack{\{\pm i\}\subset U_{0}^{c} \\ i\neq k}}(-q)^{\textbf{w}_{>i}(U)}
\cdot
v_{U\cup\{\pm i\}} \\
&=\sum_{\substack{\{\pm i\}\subset T_{0}^{c} \\ i \neq k+1}}(-q)^{\textbf{w}_{>i}(T)}
\cdot
v_{T\cup\{\pm i\}} + 
q^{-1}\sum_{\substack{\{\pm i\}\subset U_{0}^{c} \\ i\neq k}}(-q)^{\textbf{w}_{>i}(U)}
\cdot
v_{U\cup\{\pm i\}} \\
&=\sum_{\{\pm i\}\subset S_{0}^{c}}(-q)^{\textbf{w}_{>i}(S)}
\cdot
v_{T\cup\{\pm i\}} + 
q^{-1}\sum_{\{\pm i\}\subset S_{0}^{c}}(-q)^{\textbf{w}_{>i}(S)}
\cdot
v_{U\cup\{\pm i\}} \\
&=\sum_{\{\pm i\}\subset S_{0}^{c}}(-q)^{\textbf{w}_{>i}(S)}
\cdot
(v_{T\cup\{\pm i\}} + q^{-1} \cdot v_{U\cup\{\pm i\}}) \\
&=\sum_{\{\pm i\}\subset S_{0}^{c}}(-q)^{\textbf{w}_{>i}(S)}
\cdot
(v_{(S\cup\{\pm i\})\setminus\{-(k+1)\}\cup \{-k\}} + q^{-1} \cdot v_{(S\cup\{\pm i\})\setminus \{k\}\cup \{-k\}}) \\
&= f_{k}\acts\big(\sum_{\{\pm i\}\subset S_{0}^{c}}(-q)^{\textbf{w}_{>i}(S)}
\cdot
v_{S\cup\{\pm i\}}\big) \\
&=f_{k}\acts(E\acts v_{S}).
\end{align*}
All other cases can be calculated similarly.
\end{proof}

\section{Howe duality -- part 1}\label{S:HowePartOne}

We now state and prove a quantum version of 
\cite[Theorem 3.8.9.3]{Ho-perspectives-invariant-theory} for $m=1$.

\subsection{Semisimple classical Howe duality}\label{SS:HowePartOne}

We recall the classical case.

Let $V_{1}=\C\{u_{1},\dots,u_{n},u_{-n},\dots,u_{-1}\}$.
One can check that
\begin{gather*}
\ochanged{\left\{
\begin{gathered}
(u_{i},u_{-i})=1,\quad(u_{-i},u_{i})=-1,\quad
\text{ and }\quad(u_{i},u_{j})=0\text{ otherwise},
\\
\text{pairing matrix: }
J=
\begin{pmatrix}
0 & -A
\\
A & 0
\end{pmatrix}
,\quad
A:=\text{anti-diagonal matrix with $1$ on the anti-diagonal}
,
\end{gathered}
\right.}
\end{gather*}
equips $V_{1}$ with a symplectic form, which we are going to use on $\Lambda_{1}$. 
Let
\begin{gather*}
\mathfrak{sp}_{2n}=
\left\{
M\middle|JM=-M^{T}J
\right\}
=\left\{
\begin{pmatrix}
A & B
\\
C & -A^{T}
\end{pmatrix}
\middle|
B=B^{T},C=C^{T}
\right\}
\subset\mathfrak{gl}_{2n}
.
\end{gather*}
Then $\mathfrak{sp}_{2n}$ acts on $V_{1}$ as the linear endomorphisms 
$X\colon V_{1}\rightarrow V_{1}$ such that 
$(Xv,w)+(v,Xw)=0$ for all $v,w\in V$.
This is the $q=1$ version of $V_{q}$ from \autoref{D:VecRep}, i.e. $V_{1}$ is the 
\emph{vector representation} of $\mathfrak{sp}_{2n}$. This can be 
seen by using \autoref{R:ScalingVecRep} since, in elementary matrices, we have
\begin{gather*}
e_{i}\leftrightsquigarrow
E_{i,i+1}-E_{-i-1,-i}
,\quad
e_{n}\leftrightsquigarrow E_{n,-n},
\quad\text{etc.},
\end{gather*}
so that, for example, $e_{1}\acts u_{2}=u_{1}$ and $e_{1}\acts u_{-1}=-u_{-2}$.

There is an action of $\mathfrak{sp}_{2n}$ on the exterior algebra of $V$, denoted $\Lambda_{1}$, such that $X\in\mathfrak{sp}_{2n}$ acts on $\Lambda^{1}_{1}\cong V_{1}$ as before and on $\Lambda^{k}_{1}$ by derivations. Note that $\Lambda_{1}$ is the $q=1$ version of $\Lambda_{\intform}$ 
from \autoref{D:QExt} (when using the $u_{i}$ 
basis vectors). In other words, $\Lambda_{1}$ is the exterior algebra.

The $\C$-vector space 
$\Lambda_{1}$ has the $\C$-basis
\begin{gather}\label{Eq:ExtBasis}
\{u_{A,B,C}:=u_{a_{1}}\dots u_{a_{|A|}}u_{-b_{1}}\dots u_{-b_{|B|}}u_{c_{1}}u_{-c_{1}}\dots u_{c_{|C|}}u_{-c_{|C|}}\mid A\sqcup B\sqcup C\subset[1,n]\}.
\end{gather}
Note that this basis is not the classical version of the basis given by the $v_{S}$.

Let $E=\begin{psmallmatrix}0 & 1\\0 & 0\end{psmallmatrix}$, 
$F=\begin{psmallmatrix}0 & 0\\1 & 0\end{psmallmatrix}$ and 
$H=\begin{psmallmatrix}1 & 0\\0 & -1\end{psmallmatrix}$ be the usual generators of 
$\mathfrak{sl}_{2}$. The $q=1$ version of the dual action from \autoref{SS:SLAction} arises as follows 
(a precise comparison is given in \autoref{R:Dequantum}).
$\mathfrak{sl}_{2}$ acts on $\Lambda_{1}$ such that the $E$, $F$, and $H$ in $\mathfrak{sl}_{2}$ act in the basis from \autoref{Eq:ExtBasis} as
\begin{gather*}
E\acts u_{A,B,C}= 
\sum_{x\in\{ 1,\dots,n\}\setminus C} 
u_{A,B,C\cup\{x\}},
\quad
F\acts u_{A,B,C}= 
\sum_{x\in C}u_{A,B,C\setminus\{x\}},
\\
H\acts u_{A,B,C}=(-n+|A|+|B|+2|C|)\cdot u_{A,B,C}.
\end{gather*}

We now use the enveloping $\C$-algebras $U_{1}(\mathfrak{sp}_{2n})$
and $U_{1}(\mathfrak{sl}_{2})$. Let $\textbf{U}_{1}:=U_{1}(\mathfrak{sp}_{2n})\otimes U_{1}(\mathfrak{sl}_{2})$. Note that a $U_{1}(\mathfrak{sp}_{2n})$-$U_{1}(\mathfrak{sl}_{2})^{op}$-bimodule is the same as a $\textbf{U}_{1}$-module.

\begin{Theorem}\label{T:ClassicalHowe}
Consider the $\C$-vector space $\Lambda_{1}$.

\begin{enumerate}

\item There are commuting actions
\begin{gather*}
U_{1}(\mathfrak{sp}_{2n})\actsleft
\Lambda_{1}
\actsright
U_{1}(\mathfrak{sl}_{2})^{op}.
\end{gather*}
as described above.

\item Let $\phi_{1}$ and $\psi_{1}$ be the $\C$-algebra homomorphisms 
induced by the two actions from (a). Then:
\begin{gather*}
\phi_{1}\colon
U_{1}(\mathfrak{sp}_{2n})\twoheadrightarrow\End_{U_{1}(\mathfrak{sl}_{2})^{op}}(\Lambda_{1})
,\quad
\psi_{1}\colon
\changed{U_{1}(\mathfrak{sl}_{2})}\twoheadrightarrow\End_{U_{1}(\mathfrak{sp}_{2n})}(\Lambda_{1})
.
\end{gather*}
That is, the two actions generate the others centralizer.

\item The $\textbf{U}_{1}$-module
$\Lambda_{1}$ decomposes as
\begin{gather*}
\Lambda_{1}\cong\bigoplus_{k=0}^{n}L_{1}(\varpi_{k})\otimes L_{1}(n-k).
\end{gather*}
Here $L_{1}(\varpi_{k})$ and $L_{1}(n-k)$ denote the simple 
$U_{1}(\mathfrak{sp}_{2n})$-module and $U_{1}(\mathfrak{sl}_{2})^{op}$-module 
of the indicated highest weights.

\end{enumerate}

\end{Theorem}

\begin{proof}
This is a special case of \cite[Theorem 3.8.9.3]{Ho-perspectives-invariant-theory}, 
but in a different convention, so let us sketch a proof.

First, as one easily checks, the actions of $\mathfrak{sp}_{2n}$ and $\mathfrak{sl}_{2}$ 
from above commute. This is seen by identifying $E$ and $F$ with 
multiplication and contraction by the element 
$\omega:=\sum_{i=1}^{n}u_{i}u_{-i}\in\Lambda^{2}_{1}$ preserved by $\mathfrak{sp}_{2n}$. This shows (a).

Then (b) follows from (c), so we prove (c). Let $x$ be a 
highest weight vector for $\mathfrak{sp}_{2n}$ and a lowest 
weight vector for $\mathfrak{sl}_{2}$. Then using that $x$ is 
a lowest weight vector for $\mathfrak{sl}_{2}$ we get 
$x\in\Lambda^{m}_{1}$ for $0\leq m\leq n$. Since $x$ is a highest weight vector for $\mathfrak{sp}_{2n}$ in $\Lambda_{1}$, it follows \changed{from the explicit $\C$-basis in \autoref{Eq:ExtBasis} and the $\mathfrak{sp}_{2n}$-action on it by derivations} that 
$x=\sum_{C\subset\{k+1,\dots,n\}}\xi_{C}\cdot u_{\{1,\dots,k\},\emptyset,C}$ 
\changed{for some $\xi_{C}\in\C$}. 
Suppose that $C\neq\emptyset$. Then from \ochanged{$F\acts x=0$} it is possible to find a contradiction. Thus, up to scalar all the joint highest-lowest weight vectors for 
$\mathfrak{sp}_{2n}$-$\mathfrak{sl}_{2}$ are of the form 
$u_{1}u_{2}\cdots u_{k}$ for $k\in\{0,1,\dots,n\}$, and (c) follows.
\end{proof}

\subsection{Semisimple quantum Howe duality}

Using \autoref{P:QBimodule}, we will consider $\Lambda_{q}$ as a $U_{q}(\mathfrak{sp}_{2n})$-$U_{q}(\mathfrak{sl}_{2})^{op}$-bimodule or, alternatively, 
as a $\textbf{U}_{q}:=U_{q}(\mathfrak{sp}_{2n})\otimes U_{q}(\mathfrak{sl}_{2})$-module.
We use the convention that $\varpi_{0}:=0$, so the irreducible with highest weight $\varpi_0$ is the trivial module.

\begin{Lemma}\label{L:Uhwt}
Let $\epsilon\in\{\pm 1\}$.

\begin{enumerate}

\item The finite dimensional simple $\textbf{U}_{q}$-modules 
are all of the form $L_{q}(\lambda)\otimes{^{\epsilon}L}_{q}(k)$ where $\lambda$ is an integral dominant weight 
for $\mathfrak{sp}_{2n}$ and $k\in\N$, thus, $k$ is a dominant 
integral weight for $\mathfrak{sl}_{2}$. 
Their characters are given by Weyl's quantum character formula per component.

\item Let $V$ be a finite dimensional module over $\textbf{U}_{q}$. 
Then $V$ is completely reducible. Moreover, assume that $v\in V$ is such that:
\begin{enumerate}

\item $K_{i}\acts v=q^{(\alpha_{i}^{\vee},\lambda)}\cdot v$ and $(\alpha_{i}^{\vee},\lambda)\in\N$ 
for $i\in\{1,\dots,n\}$.

\item $K\acts v=\epsilon q^{k}\cdot v$ for some $k\in\N$.

\item $e_{i}\acts v=0$ for $i\in\{1,\dots,n\}$.

\item $F\acts v=0$.

\end{enumerate}
Then the $\textbf{U}_{q}$ module 
$L_{q}(\lambda)\otimes{^{\epsilon}L}_{q}(-k)$ is a direct summand of $V$.

\end{enumerate}  
\end{Lemma}

\begin{proof}
\textit{(a).} By, for example, \cite[Theorems 5.15 and 5.17]{Ja-lectures-qgroups} 
or \cite[Section 6]{AnPoWe-representation-qalgebras}
the category of finite dimensional $\textbf{U}_{q}$ modules is semisimple
and has the same character combinatorics as the classical 
case, so the claims follow from classical theory.

\textit{(b).} As in the proof of (a) this statement follows by classical theory 
and keeping track of quantum weight combinatorics.
\end{proof}

\begin{Definition}
Let $V$ be a finite dimensional $\textbf{U}_{q}$ module. 
We write $V_{\lambda,k}^{+}$ to denote the subspace of 
vectors in $V$ satisfying the conditions in 
\autoref{L:Uhwt}.(b). We refer to a vector in 
$V_{\lambda,k}^{+}$ as a \emph{singular vector} of weight $(\lambda,k)$.
\end{Definition}

\begin{Lemma}\label{L:Singular}
Let $V$ be a finite dimensional $\textbf{U}_{q}$-module. \changed{Then there 
integral dominant weights $\lambda$
for $\mathfrak{sp}_{2n}$ and $k\in\N$}
\begin{gather*}
\changed{V\cong\bigoplus_{\lambda,k}L_{q}(\lambda)\otimes\big({^{\epsilon}L}_{q}(-k)\big)^{\oplus\dim V_{\lambda, k}^{+}}.}
\end{gather*}
\end{Lemma}

\begin{proof}
Again, using \cite[Theorems 5.15 and 5.17]{Ja-lectures-qgroups} 
or \cite[Section 6]{AnPoWe-representation-qalgebras}, this follows from classical theory.
\end{proof}

\begin{Theorem}\label{T:QuantumHowe}
Consider the $\field$-vector space $\Lambda_{q}$.

\begin{enumerate}

\item There are commuting actions
\begin{gather*}
U_{q}(\mathfrak{sp}_{2n})\actsleft
\Lambda_{q}
\actsright
U_{q}(\mathfrak{sl}_{2})^{op}.
\end{gather*}
as in \autoref{S:ExtAlg}.

\item Let $\phi_{q}$ and $\psi_{q}$ be the $\field$-algebra homomorphisms 
induced by the two actions from (a). Then:
\begin{gather*}
\phi_{q}\colon
U_{q}(\mathfrak{sp}_{2n})\twoheadrightarrow\End_{U_{q}(\mathfrak{sl}_{2})^{op}}(\Lambda_{q})
,\quad
\psi_{q}\colon
\changed{U_{q}(\mathfrak{sl}_{2})}\twoheadrightarrow\End_{U_{q}(\mathfrak{sp}_{2n})}(\Lambda_{q})
.
\end{gather*}
That is, the two actions generate the others centralizer.

\item The $\textbf{U}_{q}$-module
$\Lambda_{q}$ decomposes as
\begin{gather*}
\Lambda_{q}\cong\bigoplus_{k=0}^{n}L_{q}(\varpi_{k})\otimes{^{\epsilon}L}_{q}(n-k).
\end{gather*}
Here $L_{q}(\varpi_{k})$ and ${^{\epsilon}L}_{q}(n-k)$ denote the simple 
$U_{q}(\mathfrak{sp}_{2n})$-module and $U_{q}(\mathfrak{sl}_{2})^{op}$-module 
of the indicated highest weights.
Here $\epsilon=1$ for even degrees and $\epsilon=-1$ for odd degrees.

\end{enumerate}

\end{Theorem}

\begin{proof}
\text{(a)+(b).} The first statement is \autoref{P:QBimodule}, while (b) is a consequence of (c).

\text{(c).} It is easy to see that each vector in the set
\begin{gather*}
\{ 
v_{\emptyset},v_{\{1\}},v_{\{1,2\}},\dots,v_{\{1,2,\dots,n\}}
\}
\end{gather*}
is a singular vector of weight
\begin{gather*}
(\varpi_{0},-n),
(\varpi_{1},-n+1),
(\varpi_{2},-n+2),
\dots,(\varpi_{n},0).
\end{gather*}
Thus, by \autoref{L:Singular} we have the following inclusion as a direct summand:
\begin{gather}\label{summandofqext}
\bigoplus_{k=0}^{n}L_{q}(\varpi_{k})\otimes{^{\epsilon}L}_{q}(n-k)\subset_{\oplus}\Lambda_{q}.
\end{gather}
To prove equality, recall that $\dim_{\field}L_{q}(\lambda)=\dim_{\field}L(\lambda)$ and $\dim_{\field}L_{q}(k)=\dim_{\field}L(k)$, see 
\autoref{L:Uhwt}, 
and $\dim_{\field}\Lambda_{q}=\dim_{\C}\Lambda_{1}$ by \autoref{P:ExtFlat}.
It hence follows from \autoref{T:ClassicalHowe} that 
\begin{gather*}
\dim_{\field}\big(
\bigoplus_{k=0}^{n} 
L_{q}(\varpi_{k})\otimes{^{\epsilon}L}_{q}(n-k)\big) 
=\dim_{\field}\Lambda_{q}. 
\end{gather*}
Thus, the inclusion in \autoref{summandofqext} is an equality. 

Note that the 
$U_{q}(\mathfrak{sp}_{2n})$-action preserves the degree while the $U_{q}(\mathfrak{sl}_{2})$-action preserves the parity of the degree so that we really do have two separate action pairs. 
Since the $U_{q}(\mathfrak{sl}_{2})$-action is on tensor powers of 
$X_{q}={^{1}\mathbb{F}}_{\ast}\oplus{^{-1}W}_{q}\oplus{^{1}\mathbb{F}}_{\ast}$ the even-odd 
claim follows.
\end{proof}

\begin{Remark}\label{R:Dequantum}
The meticulous reader might note that \autoref{T:QuantumHowe} does not quite quantize 
the classical version in \autoref{T:ClassicalHowe} because of the appearance of 
type-1 and type-(-1) $U_{q}(\mathfrak{sl}_{2})$-modules. This can be fixed by observing that $\Lambda_{q}^{even}$ is a type-1 $U_{q}(\mathfrak{sl}_{2})$-module and $\Lambda_{q}^{odd}$ is a type-(-1) $U_{q}(\mathfrak{sl}_{2})$-module, and then specialize these two separately.

It is actually more accurate to compare the action in \autoref{SS:HowePartOne} to the one given 
in \autoref{SS:DualAction} later on. However, the difference will not play a role for us.
\end{Remark}

\begin{Remark}\label{R:HowePartOneCoideal}
The quantum $(\mathfrak{sp}_{2n},\mathfrak{sp}_{2k})$-duality 
that was given in \cite[Theorem A]{SaTu-bcd-webs} suggests that our quantum 
$(\mathfrak{sp}_{2n},\mathfrak{sl}_{2}\cong\mathfrak{sp}_{2})$-duality 
might need a quantum symmetric pair and a quantum group instead of two quantum groups. 
(Note that \cite[Theorem A]{SaTu-bcd-webs} replaces $\mathfrak{sp}_{2n}$ 
by a quantum symmetric pair, but one might suspect the existence of a 
``dual'' of \cite[Theorem A]{SaTu-bcd-webs} with the roles of the quantum group and the 
quantum symmetric pair swapped.)

Indeed, we do have a quantum symmetric pair: According to 
\cite[Theorem A]{SaTu-bcd-webs} one of the quantum objects 
associated to symplectic exterior Howe duality is the quantum 
symmetric pair $U_{q}^{\prime}(\mathfrak{sp}_{2k})$ for the 
Satake-type diagram $\dynkin A{II}$. Algebraically, 
$U_{q}^{\prime}(\mathfrak{sp}_{2k})$ is the subalgebra of $U_{q}(\mathfrak{gl}_{2k})$ generated by
\begin{align*}
E_{i},F_{i},K_{i}^{\pm 1}\quad&\text{for }i\in\{1,3,\dots,2k-1\},
\\
B_{i}=F_{i}-K_{i}^{-1}ad(E_{i-1}E_{i+1})\acts E_{i}
\quad&\text{for }i\in\{2,4,\dots,2k\}.
\end{align*}
In the special case $k=1$ we thus get the ``small number coincidence'' 
\begin{gather*}
U_{q}^{\prime}(\mathfrak{sp}_{2k})
\cong
U_{q}(\mathfrak{sp}_{2})
\cong
U_{q}(\mathfrak{sl}_{2}),
\end{gather*}
and we thus have a quantum symmetric pair appearing, albeit 
in a roundabout way.
The same is true for the orthogonal version of \autoref{T:QuantumHowe}, 
see \cite[Theorem 3.4]{Su-harmonic-analysis-son} 
(note that this theorem is formulated differently). That theorem corresponds to 
a quantum $(\mathfrak{so}_{m},\mathfrak{sl}_{2}\cong\mathfrak{sp}_{2})$-duality
where the quantum exterior is replaced by a quantum symmetric algebra.

The remaining two dualities in \cite[Theorem A]{SaTu-bcd-webs}
come from the Satake-type diagram $\dynkin A{I}$. In this case the quantum 
symmetric pair does not correspond to a quantum group in any \ochanged{known} way, so 
a quantization might be difficult.
\end{Remark}

\begin{Corollary}\label{C:sskernelisisimple}
We have an isomorphism of 
$U_{q}(\mathfrak{sp}_{2n})$-modules
\begin{gather*}
\ker F\cap\Lambda_{q}^{k}\cong L_{q}(\varpi_{k})
\end{gather*}
with distinguished highest weight vector $v_{\{1,\dots,k\}}$ 
for the $U_{q}(\mathfrak{sp}_{2n})$-action.
\end{Corollary}

\begin{proof}
Let $V=\bigoplus_{k\in\mathbb{Z}}V_{k}$ be a finite dimensional representation of 
$U_{q}(\mathfrak{sl}_{2})$ such that $K\acts v=\epsilon q^{k}\cdot v$ for all $v\in V_{k}$. Then 
\begin{gather*}
V\cong\bigoplus_{k\in\N}
{^{\epsilon}L}_{q}(k)^{\oplus\dim_{\field}\ker F\cap V_{-k}}.
\end{gather*}
Observing that $K$ acts on $\Lambda_{q}^{k}$ as $\epsilon q^{-n+k}$ we deduce that the multiplicity of $L_{q}(n-k)$ as a direct summand of $\Lambda_{q}$ is equal to $\dim_{\field}\ker F\cap\Lambda_{q}^{k}$. By \autoref{T:QuantumHowe} we find that this multiplicity is also equal to $\dim L_{q}(\varpi_{k})$. 

Since $v_{1}\dots v_{k}\in\ker F\cap\Lambda_{q}^{k}$ there is a nonzero map 
$L_{q}(\varpi_{k})\to\ker F\cap\Lambda_{q}^{k}$. 
The representation $L_{q}(\varpi_{k})$ is simple, so the map is injective, hence, an isomorphism.
\end{proof}

\section{Integral versions and Weyl filtrations}\label{S:Weyl}

We now extend the previous results to $\intform$.

\subsection{Divided power quantum groups}\label{SS:dividedpower}

Recall from \cite{Lu-qgroups-root-of-1} that the 
\emph{divided power quantum group} $U_{\intform}(\mathfrak{sp}_{2n})$ 
is the $\intform$-subalgebra of $U_{q}(\mathfrak{sp}_{2n})$ 
generated by the \emph{divided power generators} $E_{i}^{(k)}=E_{i}^{k}/[k]_{q}$ 
and $F_{i}^{(k)}=F_{i}^{k}/[k]_{q}$ as well as an adjusted Cartan part. We 
will not need the Cartan part so we only refer to, 
for example, \cite{AnPoWe-representation-qalgebras} for details.

Before we go to the most important,
for the purpose of this paper, $U_{\intform}(\mathfrak{sp}_{2n})$-modules, 
we state a general restriction lemma. 
Let $\Res^{\field}_{\intform}$ denote the scalar restriction from 
$\field$ to $\intform$.

\begin{Lemma}\label{L:kernelispreserved}
Let $V_{q}$ and $W_{q}$ be $U_{q}(\mathfrak{sp}_{2n})$-modules 
and let $\varphi\in\Hom_{U_{q}(\mathfrak{sp}_{2n})}(V_{q},W_{q})$.
Suppose $U_{\intform}(\mathfrak{sp}_{2n})\subset U_{q}(\mathfrak{sp}_{2n})$ 
preserves $V_{\intform}\subset\Res^{\field}_{\intform}V_{q}$ and $W_{\intform}\subset\Res^{\field}_{\intform}W_{q}$, and that 
$\varphi(V_{\intform})\subset W_{\intform}$. Then 
$\ker\varphi\cap V_{\intform}$ is a 
$U_{\intform}(\mathfrak{sp}_{2n})$-submodule of $V_{\intform}$. 
\end{Lemma}

\begin{proof}
Directly from the definitions.
\end{proof}

\begin{Remark}
The theory below about Weyl and dual Weyl modules is standard, 
but mostly used over fields, see e.g. \cite{Do-tilting-alg-groups} 
and \cite{Ri-good-filtrations} for the original references. 
Some papers however are very careful with the ground ring, for example
\cite{Ka-Based-Filt} or \cite{RyHa-q-kempf}, which are helpful and are main sources for the integral statements.
\end{Remark}

As usual in the theory, see again \cite{AnPoWe-representation-qalgebras}, 
there are two types of $U_{\intform}(\mathfrak{sp}_{2n})$-modules indexed 
by dominant integral weights $\lambda\in X^{+}$:
\begin{enumerate}

\item The \emph{dual Weyl modules} $\nabla_{\intform}(\lambda)$
which we, following \cite{AnPoWe-representation-qalgebras}, define 
as induced modules where one induces a free rank one $\intform$-module over the negative Borel subalgebra, see e.g. the arXiv version of \cite{AnStTu-cellular-tilting} for details.

\item The \emph{Weyl modules} $\Delta_{\intform}(\lambda)$ that 
we define, as a free $\intform$-module, 
as the $\intform$-module dual of $\nabla_{\intform}(\lambda)$ whose action is the 
one on $\nabla_{\intform}(\lambda)$ twisted by the antipode and the Chevalley involution $\omega$,
using \autoref{D:Twist}, see 
e.g. \cite[Section 2A]{AnStTu-cellular-tilting} for details.

\end{enumerate}

Although the above definitions are the gold standard in the field, 
the following construction of the Weyl module 
is what we actually use throughout.
The next lemma also summarizes some properties of these $U_{\intform}(\mathfrak{sp}_{2n})$-modules.

\begin{Lemma}\label{L:Weylmodule-omnibuslemma}
\leavevmode

\begin{enumerate}

\item We have $\Delta_{\intform}(\lambda)\cong U_{\intform}(\mathfrak{sp}_{2n})\acts v_{\lambda}\subset L_{q}(\lambda)$, where $L_{q}(\lambda)$ is the simple $U_{q}(\mathfrak{sp}_{2n})$-module of highest weight $\lambda$ with highest weight vector $v_{\lambda}$.

\item The $U_{\intform}(\mathfrak{sp}_{2n})$-modules $\Delta_{\intform}(\lambda)$ and $\nabla_{\intform}(\lambda)$ are $\intform$-free and have the same character as the simple $U_{q}(\mathfrak{sp}_{2n})$-module $L_{q}(\lambda)$ from \autoref{L:Uhwt}.

\item We have \emph{Ext-vanishing}, i.e. $\mathrm{Ext}^{i}\big(\Delta_{\intform}(\lambda),\nabla_{\intform}(\mu)\big)=0$ unless $i=0$ and $\lambda=\mu$, in which case we have
\begin{gather*}
\Hom\big(\Delta_{\intform}(\lambda),\nabla_{\intform}(\lambda)\big)\cong\intform.
\end{gather*}

\end{enumerate}
\end{Lemma}

\begin{proof}
All of this is well-known, see for example \cite{AnPoWe-representation-qalgebras} or \cite{AnStTu-cellular-tilting} for details, and \cite{RyHa-q-kempf} for integral formulations. 
For example, the comparison in (a) can be obtained from \cite[Lemma 4.3 and Theorem 5.4]{RyHa-q-kempf}.
\end{proof}

\begin{Definition}\label{D:Tilting}
A $U_{\intform}(\mathfrak{sp}_{2n})$-module is called (integrally) \emph{tilting} if it has a 
Weyl and a dual Weyl filtration.
\end{Definition}

\begin{Example}
Let $V_{\intform}$ be as in \autoref{D:VecRep}, but over $\intform$ instead of $\field$.
For all $k\in\N$, the $U_{\intform}(\mathfrak{sp}_{2n})$-module 
$\Delta_{\intform}(\varpi_{1})^{\otimes k}\cong V_{\intform}^{\otimes k}$ is tilting.
In general, tensor products of Weyl modules for minuscule weights are examples 
of tilting modules, see e.g. \cite[Proposition 2.3]{AnStTu-semisimple-tilting} 
for an explicit statement with a self-contained proof.
\end{Example}

More generally than the previous example:

\begin{Lemma}\label{L:TensorTilting}
Tensor products of tilting modules are tilting.
\end{Lemma}

\begin{proof}
See \cite[Theorem 3.3]{Pa-tilting-tensor}, which can be adapted to work over $\intform$.
\end{proof}

Recall that \emph{specialization} in this context is:

\begin{Notation}
Let $\ring$ be a field and fix $\xi\in\ring\setminus\{0\}$. Then 
$q\mapsto\xi$ induces an action of $\intform$ on $\ring$. 
We will say that $\ring$ is a \emph{field over $\intform$} 
or a \emph{specialization}, 
leaving the data of $\xi\in\ring\setminus\{0\}$ implicit.
\end{Notation}

\begin{Example}
There are essentially four different types of specializations 
(roughly ordered by complexity):
\begin{enumerate}

\item The semisimple specializations, where $\xi$ is not a root of unity.
As an example consider the fraction field 
$\ring=\field$ with $\xi=q$. Another example is $\ring=\F_{7}(q)$ for generic $q$ 
when we let $\xi=q$.

\item The complex root of unity case. That is, $\ring=\C$ and $\xi$ is a root of unity.

\item The (strictly) mixed cases. For example, $\ring=\F_{7}$ and $\xi=2$.

\item Characteristic $p$, which, for example, could be $\ring=\F_{7}$ and $\xi=1$.

\end{enumerate}
The only simple Lie algebra where the 
representation theory is understood in all 
of the above cases is $\mathfrak{sl}_{2}$; for the combinatorics 
of the modules see e.g. \cite[Section 3.4]{Do-q-schur}.
Other expositions, with a focus on special cases, are e.g.
\cite[Section 2]{AnTu-tilting} or \cite{TuWe-quiver-tilting}.
\end{Example}

We recall the following:

\begin{Lemma}\label{L:AllFour}
For all $\lambda\in X^{+}$ there exists a simple module $L_{\ring}(\lambda)$, 
a Weyl module $\Delta_{\ring}(\lambda)$, a dual Weyl module $\nabla_{\ring}(\lambda)$ 
and an indecomposable tilting module $T_{\ring}(\lambda)$, where $\lambda$ 
is the highest weight. They are pairwise isomorphic as 
$U_{\ring}(\mathfrak{sp}_{2n})$-modules if and only if all of the four types are isomorphic.
\end{Lemma}

\begin{proof}
Again, this is well-known, see the references in the proof of 
\autoref{L:Weylmodule-omnibuslemma}.
\end{proof}

\begin{Example}
Note that integral tilting $U_{\intform}(\mathfrak{sp}_{2n})$-modules will 
be tilting after specialization, but not all tilting modules after specialization come from 
integral tilting modules. An easy example is $\Delta_{\intform}(2\varpi_{1})$ which is not 
integrally tilting but $\Delta_{q}(2\varpi_{1})$ is tilting, or more generally, 
$\Delta_{\ring}(2\varpi_{1})$ is tilting if $\xi+\xi^{-1}\neq 0$.
\end{Example}

Similarly, the \emph{divided power quantum group} $U_{\intform}(\mathfrak{sl}_{2})$ 
is the $\intform$-subalgebra of $U_{q}(\mathfrak{sl}_{2})$ 
generated by the \emph{divided power generators} $E^{(k)}=E^{k}/[k]_{q}$ 
and $F^{(k)}=F^{k}/[k]_{q}$, as well as $K$ and $K^{-1}$.

\begin{Remark}
When considering tilting modules for the quantum group for 
$\mathfrak{sl}_{2}$, one generally 
restricts to type-1 modules, since this subcategory is 
closed under tensor product. However, the $\intform$ 
version of \autoref{L:twist-sigma-equivalence} allows 
us to define type-(-1) Weyl modules, dual Weyl modules, 
and tilting modules, as $\sigma$-twists of the respective 
type-1 analogs, and the theory recalled above 
still works out. One must just take that the tensor product of two type-(-1) 
modules is of type-1 and the tensor product of a 
type-(-1) module with a type-1 module is a type-(-1) module.
(In this paper, modules that are not of type-1 will only be relevant for $\mathfrak{sl}_{2}$.)
\end{Remark}

\begin{Definition}\label{D:IntVecRep}
The \emph{vector representation} $V_{\intform}$ of $U_{\intform}(\mathfrak{sp}_{2n})$ is the free $\intform$-submodule of $V_{q}$ with basis $\{v_{1},v_{2},\dots,v_{n},v_{-n},\dots,v_{-2},v_{-1}\}$ and induced action.
\end{Definition}

Since e.g. $e_{i}^{2}v_{j}=0$, we see that $V_{\intform}$ is indeed a $U_{\intform}(\mathfrak{sp}_{2n})$-module.

\begin{Proposition}\label{P:IntActionExt}
We have the following.
\begin{enumerate}

\item There is a $U_{\intform}(\mathfrak{sp}_{2n})$-action on $\Lambda_{\intform}$, the $\intform$-version of the action from \autoref{SS:SymplecticAction}.

\item There is a $U_{\intform}(\mathfrak{sl}_{2})$-action on $\Lambda_{\intform}$, the $\intform$-version of the action from \autoref{SS:SLAction}.

\item The actions commute.

\end{enumerate}
\end{Proposition}

\begin{proof}
\textit{(a).} By symmetry between the 
actions of $e_{i}$ and $f_{i}$, it suffices to show that $f^{k}_{i}$ acts 
on $v_{S}$ with a scalar $[k]_{\intform}!$ that can be canceled. Using \autoref{L:ExtAction} 
this boils down to a case-by-case check where the hardest 
local configuration to check is
\begin{gather*}
f_{i}^{2}\acts
\begin{ytableau}
*(spinach!50)\bullet & \phantom{a} \\
\phantom{a} & *(spinach!50)\bullet
\end{ytableau}
=
f_{i}\acts\left(\,
\begin{ytableau}
*(spinach!50)\bullet & \phantom{a} \\
*(spinach!50)\bullet & \phantom{a}
\end{ytableau}
+q^{-1}\cdot
\begin{ytableau}
\phantom{a} & *(spinach!50)\bullet \\
\phantom{a} & *(spinach!50)\bullet
\end{ytableau}
\,\right)
=
[2]_{\intform}
\cdot
\begin{ytableau}
\phantom{a} & *(spinach!50)\bullet \\
*(spinach!50)\bullet & \phantom{a}
\end{ytableau}
\,.
\end{gather*}
Here $i\neq n$ and the leftmost column is the $i$th column. 
Any further application of $f_{i}$ annihilates the vectors, so we are done.
All other cases are easier and omitted.

\textit{(b).} Immediate from the construction in \autoref{SS:SLAction}.

\textit{(c).} As in \autoref{P:Commute}.
\end{proof}

Hence, we get a $U_{\intform}
:=U_{\intform}(\mathfrak{sp}_{2n})\otimes U_{\intform}(\mathfrak{sl}_{2})$-action on $\Lambda_{\intform}$. We use this action throughout.

\begin{Remark}
The following definition uses facts from the theory of \emph{Ringel duality}, 
all of which can be found in e.g. \cite[Appendix]{Do-q-schur}. 
(Note that \cite[Appendix]{Do-q-schur} works over a field, 
but the parts of the theory we use work integrally.)
The definition also generalizes immediately beyond the setting of this paper.
\end{Remark}

\begin{Definition}\label{D:ETilting}
Given a $U_{\intform}$-module $M$ we define the \emph{$M$-Schur algebras} as
\begin{gather*}
S_{\intform}^{l}=\End_{U_{\intform}(\mathfrak{sp}_{2n})}(M),
\quad
S_{\intform}^{r}=\End_{U_{\intform}(\mathfrak{sl}_{2})}(M).
\end{gather*}
We call $M$ \emph{Howe tilting} if it is 
tilting for $U_{\intform}(\mathfrak{sp}_{2n})$ and $U_{\intform}(\mathfrak{sl}_{2})$, as 
well as full tilting for $S_{\intform}^{l}$ and $S_{\intform}^{r}$, both
separately, and satisfies
\begin{gather*}
\Hom_{U_{\intform}(\mathfrak{sp}_{2n})}\big(\Delta_{\intform}(\lambda),M\big)\cong
\nabla_{\intform}(\lambda^{\prime})
,\quad
\Hom_{U_{\intform}(\mathfrak{sl}_{2})}\big(\Delta_{\intform}(\lambda^{\prime}),M\big)\cong
\nabla_{\intform}(\lambda)
\end{gather*}
for the bijection $(\placeholder)^{\prime}$ determined by Ringel duality.
\end{Definition}

A main ingredient for the integral version of the Howe duality 
from \autoref{T:QuantumHowe} is to show that $\Lambda_{\intform}$ is 
Howe tilting for the above $U_{\intform}$-action on $\Lambda_{\intform}$.

\subsection{Canonical basis}\label{SS:Rainbow}

Recall the standard basis from \autoref{D:Standard}. We now define a canonical-type basis 
(a justification for that name will be given in \autoref{SS:Canonical} below), 
but we need some notation first. Before that, 
recall from \autoref{N:Subsets} that $S^{c}=[1,-1]\setminus S$,  
$S_{0}=S\cap -S$ 
($S_{0}$ corresponds to fully dotted columns) and $S_{0}^{c}=S^{c}\cap -S^{c}$ 
($S_{0}^{c}$ corresponds to undotted columns).

\begin{Definition}
Fix $S\subset[1,-1]$. To every $x\in S_{0}$ we associate 
$x^{S}$ such that either $x^{S}\in S_{0}^{c}$ or $x^{S}=x$.
\changed{This is done inductively as follows: reading left-to-right, pick the first not yet considered fully dotted column $\pm x$ and pair it with the rightmost available undotted column 
$\pm x^{S}$ if the number of dots to the right of the fully dotted column is 
strictly smaller than the number of 
columns to the right, and set $x^{S}=x$ otherwise.}
Define $L(S):=\{x\in S_{0}|x>0,x\neq x^{S}\}$, 
$R(S):=\{x^{S}|x\in L(S)\}$ and additionally
$\binom{L(S)}{k}:=\big\{X\subset 
L(S)\,|\,|X|=k\big\}$. We also write $\pm L(S):=\{\pm x|x\in L(S)\}$ 
and $\pm R(S):=\{\pm x|x\in R(S)\}$.
\end{Definition}

\begin{Definition}
A \emph{rainbow diagram} is a dot diagram where every fully dotted column $x$
with $x\neq x^{S}$ is paired to 
some undotted column via an arc on top of the diagram connecting $x$ to $x^{S}$.
\end{Definition}

\begin{Notation}
The name rainbow diagrams comes from us thinking of rainbows connecting a fully dotted 
column with an undotted column as for example in:\vspace*{-2cm}
\begin{gather*}
\hspace*{-4cm}
\scalebox{1.5}{$\begin{tikzpicture}
\begin{scope}[xscale=1.9]
\clip (-2,0.65) rectangle (2,2.65);
\shade[shading=rainbow,shading angle=270] [even odd rule]
(0,0.65) circle (0.5) (0,0.65) circle (0.15);
\end{scope}
\node at (0,0) {$\begin{ytableau}
*(spinach!50)\bullet & \phantom{a} & \phantom{a} \\
*(spinach!50)\bullet & *(spinach!50)\bullet & \phantom{a}
\end{ytableau}$};
\node at (2,0) {$\leftrightsquigarrow$};
\node at (4,0.325) {$\begin{ytableau}
\none[1] & \none[\phantom{a}] & \none[1] \\
*(spinach!50)\bullet & \phantom{a} & \phantom{a} \\
*(spinach!50)\bullet & *(spinach!50)\bullet & \phantom{a}
\end{ytableau}$};
\end{tikzpicture}$}
\end{gather*}
However, as drawing rainbows gets a bit demanding, we decided to 
simply put numbers on top of the diagram so that the same numbers indicate the 
start and end of a rainbow; see above for an example.
\end{Notation}

\begin{Lemma}\label{L:Rainbow}
The definition of $x^{S}$ is well-defined such that every left endpoint of 
a rainbow is a fully dotted column, every right endpoint is an undotted column, and 
and rainbows do not intersect.
\end{Lemma}

\begin{proof}
Directly from the definition.
\end{proof}

\begin{Example}\label{E:Rainbow}
\leavevmode
\begin{enumerate}

\item We continue \autoref{E:DotStatistic}.
Since $S_{0}=\{5,-5\}$, so $x_{1}=5$, and $S_{0}^{c}=\{3,6,-6,-3\}$, 
and we get $x_{1}^{S}=3$. We get $L(S)=\{5\}$, $R(S)=\{3\}$, 
$\binom{L(S)}{0}=\{\emptyset\}$ and $\binom{L(S)}{1}=\big\{\{5\}\big\}$.

\item Consider the following example:
\ochanged{\begin{gather*}
S=S_{0}=
\begin{ytableau}
*(spinach!50)\bullet & *(spinach!50)\bullet & \phantom{a} & \phantom{a} \\
*(spinach!50)\bullet & *(spinach!50)\bullet & \phantom{a} & \phantom{a}
\end{ytableau}
\,,\quad
S^{c}=S_{0}^{c}=
\begin{ytableau}
\phantom{a} & \phantom{a} & *(spinach!50)\bullet & *(spinach!50)\bullet \\
\phantom{a} & \phantom{a} & *(spinach!50)\bullet & *(spinach!50)\bullet
\end{ytableau}
\,,\\
\text{rainbow: }
\begin{ytableau}
\none[1] & \none[2]  & \none[2] & \none[1] \\
*(spinach!50)\bullet & *(spinach!50)\bullet & \phantom{a} & \phantom{a} \\
*(spinach!50)\bullet & *(spinach!50)\bullet & \phantom{a} & \phantom{a} \\
\none[\phantom{a}] & \none[\phantom{a}]  & \none[\phantom{a}] & \none[\phantom{a}]
\end{ytableau}
\,.
\end{gather*}}
We can see that $\pm L(S)=\{1,2,-2,-1\}$, the left end 
points of the rainbows, $\pm R(S)=\{3,4,-4,-3\}$, the right 
end points of the rainbows.

\end{enumerate}
More examples are given below.
\end{Example}

\begin{Definition}\label{D:Canonical}
Let $S\subset[1,-1]$. We define \emph{canonical basis vectors} by
\begin{gather*}
b_{S}:= 
\sum_{p=0}^{|L(S)|}\sum_{\{x_{1},\dots,x_{p}\}\in\binom{L(S)}{p}}
q^{-p}\cdot \changed{v_{(S\setminus\{\pm x_{1},\dots,\pm x_{p}\})\cup\{\pm x_{1}^{S},\dots,\pm x_{p}^{S}\}}}
.
\end{gather*}
The set $\{b_{S}|S\subset[1,-1]\}$
is the \emph{canonical basis}.
\end{Definition}

The rainbow diagrams give the canonical basis by \textit{jumping dots along rainbows}:

\begin{Example}\label{E:RainbowCalcs}
\leavevmode
\begin{enumerate}

\item For $n=1$ we have $b_{S}=v_{S}$ for all $S\subset[1,-1]$ since in this case $S_{0}=\emptyset$ 
or $S_{0}^{c}=\emptyset$.

\item For $n=2$ we already see a difference between the standard and the canonical bases. 
Precisely, let $S=\{1,-1\}$ so that $S_{0}=S$ and $S_{0}^{c}=\{2,-2\}$, 
and we have $1^{S}=2$. For $T=\{2,-1\}$ we have $T_{0}=\emptyset$. We get
\begin{gather*}
\ochanged{b_{\{1,-1\}}=
\begin{ytableau}
\none[1] & \none[1] \\
*(spinach!50)\bullet & \phantom{a} \\
*(spinach!50)\bullet & \phantom{a} \\
\none[\phantom{1}] & \none[\phantom{1}]
\end{ytableau}
=
\begin{ytableau}
*(spinach!50)\bullet & \phantom{a} \\
*(spinach!50)\bullet & \phantom{a}
\end{ytableau}
+q^{-1}\cdot
\begin{ytableau}
\phantom{a} & *(spinach!50)\bullet \\
\phantom{a} & *(spinach!50)\bullet
\end{ytableau}
\,,}\\
\ochanged{f_{1}\acts b_{\{1,-1\}}
=[2]_{\intform}\cdot
\begin{ytableau}
\phantom{a} & *(spinach!50)\bullet \\
*(spinach!50)\bullet & \phantom{a}
\end{ytableau}
=[2]_{\intform}\cdot b_{\{2,-1\}},}
\end{gather*}
where we use dot diagrams for the standard basis.

\item Continuing \autoref{E:Rainbow}.(b), we get that $b_{\{1,2,-2,-1\}}$ is equal to
\ochanged{\begin{align*}
\begin{ytableau}
\none[1] & \none[2]  & \none[2] & \none[1] \\
*(spinach!50)\bullet & *(spinach!50)\bullet & \phantom{a} & \phantom{a} \\
*(spinach!50)\bullet & *(spinach!50)\bullet & \phantom{a} & \phantom{a} \\
\none[\phantom{a}] & \none[\phantom{a}]  & \none[\phantom{a}] & \none[\phantom{a}]
\end{ytableau}
=&
\begin{ytableau}
*(spinach!50)\bullet & *(spinach!50)\bullet & \phantom{a} & \phantom{a} \\
*(spinach!50)\bullet & *(spinach!50)\bullet & \phantom{a} & \phantom{a}
\end{ytableau}
+q^{-1}\cdot
\begin{ytableau}
*(spinach!50)\bullet & \phantom{a} & *(spinach!50)\bullet & \phantom{a} \\
*(spinach!50)\bullet & \phantom{a} & *(spinach!50)\bullet & \phantom{a}
\end{ytableau}
\\
&+q^{-1}\cdot
\begin{ytableau}
\phantom{a} & *(spinach!50)\bullet & \phantom{a} & *(spinach!50)\bullet \\
\phantom{a} & *(spinach!50)\bullet & \phantom{a} & *(spinach!50)\bullet
\end{ytableau}
+q^{-2}\cdot
\begin{ytableau}
\phantom{a} & \phantom{a} & *(spinach!50)\bullet & *(spinach!50)\bullet \\
\phantom{a} & \phantom{a} & *(spinach!50)\bullet & *(spinach!50)\bullet
\end{ytableau}
\,.
\end{align*}}
We get that $f_{2}\acts b_{\{1,2,-2,-1\}}=[2]_{\intform}\cdot b_{\{1,3,-2,-1\}}$
since $b_{\{1,3,-2,-1\}}$ is 
\begin{gather*}
\begin{ytableau}
\none[1] & \none[2]  & \none[2] & \none[1] \\
*(spinach!50)\bullet & \phantom{a} &*(spinach!50)\bullet & \phantom{a} \\
*(spinach!50)\bullet & *(spinach!50)\bullet & \phantom{a} & \phantom{a} \\
\none[\phantom{a}] & \none[\phantom{a}]  & \none[\phantom{a}] & \none[\phantom{a}]
\end{ytableau}
=
\begin{ytableau}
*(spinach!50)\bullet & \phantom{a} &*(spinach!50)\bullet & \phantom{a} \\
*(spinach!50)\bullet & *(spinach!50)\bullet & \phantom{a} & \phantom{a}
\end{ytableau}
+q^{-1}\cdot
\begin{ytableau}
\phantom{a} & \phantom{a} & *(spinach!50)\bullet & *(spinach!50)\bullet \\
\phantom{a} & *(spinach!50)\bullet & \phantom{a} & *(spinach!50)\bullet
\end{ytableau}
\,.
\end{gather*}

\end{enumerate}
In general, $b_{S}=v_{S}$ if $S_{0}=\emptyset$ 
or $S_{0}^{c}=\emptyset$.
\end{Example}

As often in the theory of canonical bases, the explicit form of the basis elements
is less important than their properties. The same is true for us. For example, 
the $U_{\intform}(\mathfrak{sp}_{2n})$-action on the canonical basis is, to a large extend, 
coefficient free:

\begin{Lemma}\label{L:IntExtAction}
For $i\in\{1,\dots,n-1\}$ and $i=n$ we have:
\begin{align*}
f_{i}\acts b_{S}&= 
\begin{cases}
b_{(S\setminus\{i\})\cup\{i+1\}} & \text{if }S_{i,i+1}=\{ i \}, 
\\
b_{(S\setminus\{-(i+1)\})\cup\{-i\}} & \text{if }S_{i,i+1}=\{-(i+1)\},
\\
b_{(S\setminus\{-(i+1)\})\cup\{-i\}} & \text{if }S_{i,i+1}=\{i,-(i+1)\},
\\
[2]_{\intform}\cdot b_{(S\setminus\{i\})\cup\{i+1\}} & \text{if }S_{i,i+1}=\{i,-i\},
\\
b_{(S\setminus\{-(i+1)\})\cup\{-i\}} & \text{if }S_{i,i+1}=\{i+1,-(i+1)\},
\\
b_{(S\setminus\{i\})\cup\{i+1\}} & \text{if }S_{i,i+1}=\{i,-(i+1),-i\},
\\
b_{(S\setminus\{-(i+1)\})\cup\{-i\}} & \text{if }S_{i,i+1}=\{i,i+1,-(i+1)\},
\end{cases}
\\
f_{n}\cdot b_{S}&=
\begin{cases}
b_{(S\setminus\{n\})\cup\{-n\}} & \text{if }S_{n}=\{n\}.
\end{cases}
\end{align*}
\end{Lemma}

\begin{proof}
This is an exercise of how to work with rainbow diagrams. 
We therefore prove only one of the above equations, namely the one with $[2]_{\intform}$ 
which is the most difficult.

In this case there are five local rainbow diagrams, where the dots are in the $i$th column:
\begin{gather*}
\begin{ytableau}
\none[1] & \none[1] \\
*(spinach!50)\bullet & \phantom{a} \\
*(spinach!50)\bullet & \phantom{a} \\
\none[\phantom{a}] & \none[\phantom{a}]
\end{ytableau}
\,,\quad
\begin{ytableau}
\none[1] & \none[2] \\
*(spinach!50)\bullet & \phantom{a} \\
*(spinach!50)\bullet & \phantom{a} \\
\none[\phantom{a}] & \none[\phantom{a}]
\end{ytableau}
\,,\quad
\begin{ytableau}
\none[1] & \none[\phantom{a}] \\
*(spinach!50)\bullet & \phantom{a} \\
*(spinach!50)\bullet & \phantom{a} \\
\none[\phantom{a}] & \none[\phantom{a}]
\end{ytableau}
\,,\quad
\begin{ytableau}
\none[\phantom{a}] & \none[1] \\
*(spinach!50)\bullet & \phantom{a} \\
*(spinach!50)\bullet & \phantom{a} \\
\none[\phantom{a}] & \none[\phantom{a}]
\end{ytableau}
\,,\quad
\begin{ytableau}
\none[\phantom{a}] & \none[\phantom{a}] \\
*(spinach!50)\bullet & \phantom{a} \\
*(spinach!50)\bullet & \phantom{a} \\
\none[\phantom{a}] & \none[\phantom{a}]
\end{ytableau}
\,.
\end{gather*}
The left case can be checked as in \autoref{E:RainbowCalcs}, and none of the 
other cases can appear by the construction of rainbow diagrams, cf. \autoref{L:Rainbow}.
\end{proof}

\begin{Remark}
Note the missing ``otherwise'' when comparing \autoref{L:IntExtAction} 
to \autoref{L:ExtAction}. Indeed, there are 
other, rather complicated, coefficients such as 
for $f_{3}\acts b_{\{1,2\}}$ when $n=4$.
\end{Remark}

Next, we partially describe the action of $U_{\intform}(\mathfrak{sl}_{2})$ 
on the canonical basis.

\subsection{Crystal basis for fundamental Weyl modules}

The crystal combinatorics for fundamental weights in classical types is nicely summarized 
in \cite[Chapter 6]{BuSc-crystal-bases}, which is also our main source for the below.

\begin{Example}
Before we give the abstract definitions, we start with an example that the reader 
can keep in mind while reading the below.

Roughly, the setting is as follows:
\begin{enumerate}

\item For $\mathfrak{sp}_{2n}$ the ordered underlying set of fillings is 
$\{1<\dots<n<-n<\dots<-1\}$. Note that, as often in the theory of crystals, 
$-i$ is denoted $\overline{i}$ in the examples below.

\item The crystal graph for $\varpi_{k}$ is modeled on Young tableaux with 
one column and $k$ rows (a \emph{column tableaux}). 
That is, the vertices of this graph are Young tableaux 
with entries from $\{1<2<\dots<n<-n<\dots<-2<-1\}$ that
strictly increase when reading columns from top to bottom.

\item The crystal operator $\ff_{i}$, illustrated as $i$-colored edges, 
changes $i$ to $i+1$, or $-(i+1)$ to $-i$ for $i\in\{1,\dots,n-1\}$, 
and $\ff_{n}$ changes $n$ to $-n$. The crystal operator $\ee_{i}$ 
does the opposite and is usually omitted from the pictures.

\end{enumerate}
For $n=3$, that is $\mathfrak{sp}_{6}$, the three
fundamental crystals are as follows.
\begin{gather*}
\varpi_{1}\colon
\scalebox{0.75}{$\begin{tikzpicture}[>=latex,line join=bevel,anchorbase,yscale=0.75]
\node (node_{0}) at (8.5bp,301.5bp) [draw,draw=none] {${\def\lr#1{\multicolumn{1}{|@{\hspace{.6ex}}c@{\hspace{.6ex}}|}{\raisebox{-.3ex}{$#1$}}}\raisebox{-.6ex}{$\begin{array}[b]{*{1}c}\cline{1-1}\lr{2}\\\cline{1-1}\end{array}$}}$};
\node (node_5) at (8.5bp,228.5bp) [draw,draw=none] {${\def\lr#1{\multicolumn{1}{|@{\hspace{.6ex}}c@{\hspace{.6ex}}|}{\raisebox{-.3ex}{$#1$}}}\raisebox{-.6ex}{$\begin{array}[b]{*{1}c}\cline{1-1}\lr{3}\\\cline{1-1}\end{array}$}}$};
\node (node_{1}) at (8.5bp,82.5bp) [draw,draw=none] {${\def\lr#1{\multicolumn{1}{|@{\hspace{.6ex}}c@{\hspace{.6ex}}|}{\raisebox{-.3ex}{$#1$}}}\raisebox{-.6ex}{$\begin{array}[b]{*{1}c}\cline{1-1}\lr{\overline{2}}\\\cline{1-1}\end{array}$}}$};
\node (node_2) at (8.5bp,9.5bp) [draw,draw=none] {${\def\lr#1{\multicolumn{1}{|@{\hspace{.6ex}}c@{\hspace{.6ex}}|}{\raisebox{-.3ex}{$#1$}}}\raisebox{-.6ex}{$\begin{array}[b]{*{1}c}\cline{1-1}\lr{\overline{1}}\\\cline{1-1}\end{array}$}}$};
\node (node_3) at (8.5bp,374.5bp) [draw,draw=none] {${\def\lr#1{\multicolumn{1}{|@{\hspace{.6ex}}c@{\hspace{.6ex}}|}{\raisebox{-.3ex}{$#1$}}}\raisebox{-.6ex}{$\begin{array}[b]{*{1}c}\cline{1-1}\lr{1}\\\cline{1-1}\end{array}$}}$};
\node (node_4) at (8.5bp,155.5bp) [draw,draw=none] {${\def\lr#1{\multicolumn{1}{|@{\hspace{.6ex}}c@{\hspace{.6ex}}|}{\raisebox{-.3ex}{$#1$}}}\raisebox{-.6ex}{$\begin{array}[b]{*{1}c}\cline{1-1}\lr{\overline{3}}\\\cline{1-1}\end{array}$}}$};
\draw [red,->] (node_{0}) ..controls (8.5bp,281.04bp) and (8.5bp,262.45bp)  .. (node_5);
\definecolor{strokecol}{rgb}{0.0,0.0,0.0};
\pgfsetstrokecolor{strokecol}
\draw (17.0bp,265.0bp) node {$2$};
\draw [blue,->] (node_{1}) ..controls (8.5bp,62.042bp) and (8.5bp,43.449bp)  .. (node_2);
\draw (17.0bp,46.0bp) node {$1$};
\draw [blue,->] (node_3) ..controls (8.5bp,354.04bp) and (8.5bp,335.45bp)  .. (node_{0});
\draw (17.0bp,338.0bp) node {$1$};
\draw [red,->] (node_4) ..controls (8.5bp,135.04bp) and (8.5bp,116.45bp)  .. (node_{1});
\draw (17.0bp,119.0bp) node {$2$};
\draw [green,->] (node_5) ..controls (8.5bp,208.04bp) and (8.5bp,189.45bp)  .. (node_4);
\draw (17.0bp,192.0bp) node {$3$};
\end{tikzpicture}$}
,\quad
\varpi_{2}\colon
\scalebox{0.75}{$\begin{tikzpicture}[>=latex,line join=bevel,anchorbase,yscale=0.75]
\node (node_{0}) at (68.5bp,695.5bp) [draw,draw=none] {${\def\lr#1{\multicolumn{1}{|@{\hspace{.6ex}}c@{\hspace{.6ex}}|}{\raisebox{-.3ex}{$#1$}}}\raisebox{-.6ex}{$\begin{array}[b]{*{1}c}\cline{1-1}\lr{1}\\\cline{1-1}\lr{2}\\\cline{1-1}\end{array}$}}$};
\node (node_3) at (68.5bp,610.5bp) [draw,draw=none] {${\def\lr#1{\multicolumn{1}{|@{\hspace{.6ex}}c@{\hspace{.6ex}}|}{\raisebox{-.3ex}{$#1$}}}\raisebox{-.6ex}{$\begin{array}[b]{*{1}c}\cline{1-1}\lr{1}\\\cline{1-1}\lr{3}\\\cline{1-1}\end{array}$}}$};
\node (node_{1}) at (48.5bp,525.5bp) [draw,draw=none] {${\def\lr#1{\multicolumn{1}{|@{\hspace{.6ex}}c@{\hspace{.6ex}}|}{\raisebox{-.3ex}{$#1$}}}\raisebox{-.6ex}{$\begin{array}[b]{*{1}c}\cline{1-1}\lr{1}\\\cline{1-1}\lr{\overline{3}}\\\cline{1-1}\end{array}$}}$};
\node (node_4) at (8.5bp,440.5bp) [draw,draw=none] {${\def\lr#1{\multicolumn{1}{|@{\hspace{.6ex}}c@{\hspace{.6ex}}|}{\raisebox{-.3ex}{$#1$}}}\raisebox{-.6ex}{$\begin{array}[b]{*{1}c}\cline{1-1}\lr{1}\\\cline{1-1}\lr{\overline{2}}\\\cline{1-1}\end{array}$}}$};
\node (node_{1}1) at (68.5bp,440.5bp) [draw,draw=none] {${\def\lr#1{\multicolumn{1}{|@{\hspace{.6ex}}c@{\hspace{.6ex}}|}{\raisebox{-.3ex}{$#1$}}}\raisebox{-.6ex}{$\begin{array}[b]{*{1}c}\cline{1-1}\lr{2}\\\cline{1-1}\lr{\overline{3}}\\\cline{1-1}\end{array}$}}$};
\node (node_2) at (68.5bp,15.5bp) [draw,draw=none] {${\def\lr#1{\multicolumn{1}{|@{\hspace{.6ex}}c@{\hspace{.6ex}}|}{\raisebox{-.3ex}{$#1$}}}\raisebox{-.6ex}{$\begin{array}[b]{*{1}c}\cline{1-1}\lr{\overline{2}}\\\cline{1-1}\lr{\overline{1}}\\\cline{1-1}\end{array}$}}$};
\node (node_{1}2) at (88.5bp,525.5bp) [draw,draw=none] {${\def\lr#1{\multicolumn{1}{|@{\hspace{.6ex}}c@{\hspace{.6ex}}|}{\raisebox{-.3ex}{$#1$}}}\raisebox{-.6ex}{$\begin{array}[b]{*{1}c}\cline{1-1}\lr{2}\\\cline{1-1}\lr{3}\\\cline{1-1}\end{array}$}}$};
\node (node_{1}3) at (8.5bp,355.5bp) [draw,draw=none] {${\def\lr#1{\multicolumn{1}{|@{\hspace{.6ex}}c@{\hspace{.6ex}}|}{\raisebox{-.3ex}{$#1$}}}\raisebox{-.6ex}{$\begin{array}[b]{*{1}c}\cline{1-1}\lr{2}\\\cline{1-1}\lr{\overline{2}}\\\cline{1-1}\end{array}$}}$};
\node (node_5) at (88.5bp,185.5bp) [draw,draw=none] {${\def\lr#1{\multicolumn{1}{|@{\hspace{.6ex}}c@{\hspace{.6ex}}|}{\raisebox{-.3ex}{$#1$}}}\raisebox{-.6ex}{$\begin{array}[b]{*{1}c}\cline{1-1}\lr{\overline{3}}\\\cline{1-1}\lr{\overline{2}}\\\cline{1-1}\end{array}$}}$};
\node (node_6) at (68.5bp,100.5bp) [draw,draw=none] {${\def\lr#1{\multicolumn{1}{|@{\hspace{.6ex}}c@{\hspace{.6ex}}|}{\raisebox{-.3ex}{$#1$}}}\raisebox{-.6ex}{$\begin{array}[b]{*{1}c}\cline{1-1}\lr{\overline{3}}\\\cline{1-1}\lr{\overline{1}}\\\cline{1-1}\end{array}$}}$};
\node (node_7) at (68.5bp,355.5bp) [draw,draw=none] {${\def\lr#1{\multicolumn{1}{|@{\hspace{.6ex}}c@{\hspace{.6ex}}|}{\raisebox{-.3ex}{$#1$}}}\raisebox{-.6ex}{$\begin{array}[b]{*{1}c}\cline{1-1}\lr{3}\\\cline{1-1}\lr{\overline{3}}\\\cline{1-1}\end{array}$}}$};
\node (node_8) at (68.5bp,270.5bp) [draw,draw=none] {${\def\lr#1{\multicolumn{1}{|@{\hspace{.6ex}}c@{\hspace{.6ex}}|}{\raisebox{-.3ex}{$#1$}}}\raisebox{-.6ex}{$\begin{array}[b]{*{1}c}\cline{1-1}\lr{3}\\\cline{1-1}\lr{\overline{2}}\\\cline{1-1}\end{array}$}}$};
\node (node_9) at (48.5bp,185.5bp) [draw,draw=none] {${\def\lr#1{\multicolumn{1}{|@{\hspace{.6ex}}c@{\hspace{.6ex}}|}{\raisebox{-.3ex}{$#1$}}}\raisebox{-.6ex}{$\begin{array}[b]{*{1}c}\cline{1-1}\lr{3}\\\cline{1-1}\lr{\overline{1}}\\\cline{1-1}\end{array}$}}$};
\node (node_{1}0) at (8.5bp,270.5bp) [draw,draw=none] {${\def\lr#1{\multicolumn{1}{|@{\hspace{.6ex}}c@{\hspace{.6ex}}|}{\raisebox{-.3ex}{$#1$}}}\raisebox{-.6ex}{$\begin{array}[b]{*{1}c}\cline{1-1}\lr{2}\\\cline{1-1}\lr{\overline{1}}\\\cline{1-1}\end{array}$}}$};
\draw [red,->] (node_{0}) ..controls (68.5bp,667.62bp) and (68.5bp,650.39bp)  .. (node_3);
\definecolor{strokecol}{rgb}{0.0,0.0,0.0};
\pgfsetstrokecolor{strokecol}
\draw (77.0bp,653.0bp) node {$2$};
\draw [red,->] (node_{1}) ..controls (35.428bp,507.14bp) and (30.255bp,499.42bp)  .. (26.5bp,492.0bp) .. controls (22.315bp,483.74bp) and (18.663bp,474.26bp)  .. (node_4);
\draw (35.0bp,483.0bp) node {$2$};
\draw [blue,->] (node_{1}) ..controls (46.186bp,499.41bp) and (46.221bp,485.65bp)  .. (49.5bp,474.0bp) .. controls (50.695bp,469.75bp) and (52.581bp,465.49bp)  .. (node_{1}1);
\draw (58.0bp,483.0bp) node {$1$};
\draw [green,->] (node_3) ..controls (55.933bp,592.16bp) and (51.625bp,584.55bp)  .. (49.5bp,577.0bp) .. controls (47.22bp,568.9bp) and (46.509bp,559.78bp)  .. (node_{1});
\draw (58.0bp,568.0bp) node {$3$};
\draw [blue,->] (node_3) ..controls (74.977bp,582.62bp) and (79.128bp,565.39bp)  .. (node_{1}2);
\draw (88.0bp,568.0bp) node {$1$};
\draw [blue,->] (node_4) ..controls (8.5bp,412.62bp) and (8.5bp,395.39bp)  .. (node_{1}3);
\draw (17.0bp,398.0bp) node {$1$};
\draw [blue,->] (node_5) ..controls (82.023bp,157.62bp) and (77.872bp,140.39bp)  .. (node_6);
\draw (88.0bp,143.0bp) node {$1$};
\draw [red,->] (node_6) ..controls (68.5bp,72.622bp) and (68.5bp,55.392bp)  .. (node_2);
\draw (77.0bp,58.0bp) node {$2$};
\draw [red,->] (node_7) ..controls (68.5bp,327.62bp) and (68.5bp,310.39bp)  .. (node_8);
\draw (77.0bp,313.0bp) node {$2$};
\draw [green,->] (node_8) ..controls (74.977bp,242.62bp) and (79.128bp,225.39bp)  .. (node_5);
\draw (88.0bp,228.0bp) node {$3$};
\draw [blue,->] (node_8) ..controls (59.475bp,249.37bp) and (57.088bp,242.99bp)  .. (55.5bp,237.0bp) .. controls (53.324bp,228.79bp) and (51.796bp,219.64bp)  .. (node_9);
\draw (64.0bp,228.0bp) node {$1$};
\draw [green,->] (node_9) ..controls (50.581bp,159.56bp) and (52.37bp,145.81bp)  .. (55.5bp,134.0bp) .. controls (56.22bp,131.29bp) and (57.103bp,128.49bp)  .. (node_6);
\draw (64.0bp,143.0bp) node {$3$};
\draw [red,->] (node_{1}0) ..controls (19.434bp,244.69bp) and (25.662bp,230.94bp)  .. (31.5bp,219.0bp) .. controls (32.882bp,216.17bp) and (34.363bp,213.23bp)  .. (node_9);
\draw (40.0bp,228.0bp) node {$2$};
\draw [red,->] (node_{1}1) ..controls (68.5bp,412.62bp) and (68.5bp,395.39bp)  .. (node_7);
\draw (77.0bp,398.0bp) node {$2$};
\draw [green,->] (node_{1}2) ..controls (82.023bp,497.62bp) and (77.872bp,480.39bp)  .. (node_{1}1);
\draw (88.0bp,483.0bp) node {$3$};
\draw [blue,->] (node_{1}3) ..controls (8.5bp,327.62bp) and (8.5bp,310.39bp)  .. (node_{1}0);
\draw (17.0bp,313.0bp) node {$1$};
\end{tikzpicture}$}
,\quad
\varpi_{3}\colon
\scalebox{0.75}{$\begin{tikzpicture}[>=latex,line join=bevel,anchorbase,yscale=0.75]
\node (node_{0}) at (68.5bp,700.5bp) [draw,draw=none] {${\def\lr#1{\multicolumn{1}{|@{\hspace{.6ex}}c@{\hspace{.6ex}}|}{\raisebox{-.3ex}{$#1$}}}\raisebox{-.6ex}{$\begin{array}[b]{*{1}c}\cline{1-1}\lr{1}\\\cline{1-1}\lr{3}\\\cline{1-1}\lr{\overline{3}}\\\cline{1-1}\end{array}$}}$};
\node (node_3) at (48.5bp,603.5bp) [draw,draw=none] {${\def\lr#1{\multicolumn{1}{|@{\hspace{.6ex}}c@{\hspace{.6ex}}|}{\raisebox{-.3ex}{$#1$}}}\raisebox{-.6ex}{$\begin{array}[b]{*{1}c}\cline{1-1}\lr{1}\\\cline{1-1}\lr{3}\\\cline{1-1}\lr{\overline{2}}\\\cline{1-1}\end{array}$}}$};
\node (node_7) at (88.5bp,603.5bp) [draw,draw=none] {${\def\lr#1{\multicolumn{1}{|@{\hspace{.6ex}}c@{\hspace{.6ex}}|}{\raisebox{-.3ex}{$#1$}}}\raisebox{-.6ex}{$\begin{array}[b]{*{1}c}\cline{1-1}\lr{2}\\\cline{1-1}\lr{3}\\\cline{1-1}\lr{\overline{3}}\\\cline{1-1}\end{array}$}}$};
\node (node_{1}) at (68.5bp,506.5bp) [draw,draw=none] {${\def\lr#1{\multicolumn{1}{|@{\hspace{.6ex}}c@{\hspace{.6ex}}|}{\raisebox{-.3ex}{$#1$}}}\raisebox{-.6ex}{$\begin{array}[b]{*{1}c}\cline{1-1}\lr{2}\\\cline{1-1}\lr{3}\\\cline{1-1}\lr{\overline{2}}\\\cline{1-1}\end{array}$}}$};
\node (node_2) at (88.5bp,409.5bp) [draw,draw=none] {${\def\lr#1{\multicolumn{1}{|@{\hspace{.6ex}}c@{\hspace{.6ex}}|}{\raisebox{-.3ex}{$#1$}}}\raisebox{-.6ex}{$\begin{array}[b]{*{1}c}\cline{1-1}\lr{2}\\\cline{1-1}\lr{3}\\\cline{1-1}\lr{\overline{1}}\\\cline{1-1}\end{array}$}}$};
\node (node_8) at (28.5bp,409.5bp) [draw,draw=none] {${\def\lr#1{\multicolumn{1}{|@{\hspace{.6ex}}c@{\hspace{.6ex}}|}{\raisebox{-.3ex}{$#1$}}}\raisebox{-.6ex}{$\begin{array}[b]{*{1}c}\cline{1-1}\lr{2}\\\cline{1-1}\lr{\overline{3}}\\\cline{1-1}\lr{\overline{2}}\\\cline{1-1}\end{array}$}}$};
\node (node_9) at (48.5bp,312.5bp) [draw,draw=none] {${\def\lr#1{\multicolumn{1}{|@{\hspace{.6ex}}c@{\hspace{.6ex}}|}{\raisebox{-.3ex}{$#1$}}}\raisebox{-.6ex}{$\begin{array}[b]{*{1}c}\cline{1-1}\lr{2}\\\cline{1-1}\lr{\overline{3}}\\\cline{1-1}\lr{\overline{1}}\\\cline{1-1}\end{array}$}}$};
\node (node_{1}2) at (8.5bp,506.5bp) [draw,draw=none] {${\def\lr#1{\multicolumn{1}{|@{\hspace{.6ex}}c@{\hspace{.6ex}}|}{\raisebox{-.3ex}{$#1$}}}\raisebox{-.6ex}{$\begin{array}[b]{*{1}c}\cline{1-1}\lr{1}\\\cline{1-1}\lr{\overline{3}}\\\cline{1-1}\lr{\overline{2}}\\\cline{1-1}\end{array}$}}$};
\node (node_4) at (8.5bp,312.5bp) [draw,draw=none] {${\def\lr#1{\multicolumn{1}{|@{\hspace{.6ex}}c@{\hspace{.6ex}}|}{\raisebox{-.3ex}{$#1$}}}\raisebox{-.6ex}{$\begin{array}[b]{*{1}c}\cline{1-1}\lr{3}\\\cline{1-1}\lr{\overline{3}}\\\cline{1-1}\lr{\overline{2}}\\\cline{1-1}\end{array}$}}$};
\node (node_5) at (28.5bp,215.5bp) [draw,draw=none] {${\def\lr#1{\multicolumn{1}{|@{\hspace{.6ex}}c@{\hspace{.6ex}}|}{\raisebox{-.3ex}{$#1$}}}\raisebox{-.6ex}{$\begin{array}[b]{*{1}c}\cline{1-1}\lr{3}\\\cline{1-1}\lr{\overline{3}}\\\cline{1-1}\lr{\overline{1}}\\\cline{1-1}\end{array}$}}$};
\node (node_{1}1) at (28.5bp,118.5bp) [draw,draw=none] {${\def\lr#1{\multicolumn{1}{|@{\hspace{.6ex}}c@{\hspace{.6ex}}|}{\raisebox{-.3ex}{$#1$}}}\raisebox{-.6ex}{$\begin{array}[b]{*{1}c}\cline{1-1}\lr{3}\\\cline{1-1}\lr{\overline{2}}\\\cline{1-1}\lr{\overline{1}}\\\cline{1-1}\end{array}$}}$};
\node (node_6) at (68.5bp,894.5bp) [draw,draw=none] {${\def\lr#1{\multicolumn{1}{|@{\hspace{.6ex}}c@{\hspace{.6ex}}|}{\raisebox{-.3ex}{$#1$}}}\raisebox{-.6ex}{$\begin{array}[b]{*{1}c}\cline{1-1}\lr{1}\\\cline{1-1}\lr{2}\\\cline{1-1}\lr{3}\\\cline{1-1}\end{array}$}}$};
\node (node_{1}3) at (68.5bp,797.5bp) [draw,draw=none] {${\def\lr#1{\multicolumn{1}{|@{\hspace{.6ex}}c@{\hspace{.6ex}}|}{\raisebox{-.3ex}{$#1$}}}\raisebox{-.6ex}{$\begin{array}[b]{*{1}c}\cline{1-1}\lr{1}\\\cline{1-1}\lr{2}\\\cline{1-1}\lr{\overline{3}}\\\cline{1-1}\end{array}$}}$};
\node (node_{1}0) at (28.5bp,21.5bp) [draw,draw=none] {${\def\lr#1{\multicolumn{1}{|@{\hspace{.6ex}}c@{\hspace{.6ex}}|}{\raisebox{-.3ex}{$#1$}}}\raisebox{-.6ex}{$\begin{array}[b]{*{1}c}\cline{1-1}\lr{\overline{3}}\\\cline{1-1}\lr{\overline{2}}\\\cline{1-1}\lr{\overline{1}}\\\cline{1-1}\end{array}$}}$};
\draw [red,->] (node_{0}) ..controls (55.916bp,678.69bp) and (51.609bp,669.69bp)  .. (49.5bp,661.0bp) .. controls (47.515bp,652.82bp) and (46.749bp,643.78bp)  .. (node_3);
\definecolor{strokecol}{rgb}{0.0,0.0,0.0};
\pgfsetstrokecolor{strokecol}
\draw (58.0bp,652.0bp) node {$2$};
\draw [blue,->] (node_{0}) ..controls (75.539bp,666.06bp) and (79.059bp,649.34bp)  .. (node_7);
\draw (88.0bp,652.0bp) node {$1$};
\draw [blue,->] (node_{1}) ..controls (75.539bp,472.06bp) and (79.059bp,455.34bp)  .. (node_2);
\draw (88.0bp,458.0bp) node {$1$};
\draw [green,->] (node_{1}) ..controls (54.332bp,471.85bp) and (47.147bp,454.79bp)  .. (node_8);
\draw (61.0bp,458.0bp) node {$3$};
\draw [green,->] (node_2) ..controls (74.332bp,374.85bp) and (67.147bp,357.79bp)  .. (node_9);
\draw (80.0bp,361.0bp) node {$3$};
\draw [blue,->] (node_3) ..controls (46.318bp,571.01bp) and (46.677bp,557.63bp)  .. (49.5bp,546.0bp) .. controls (50.785bp,540.7bp) and (52.887bp,535.29bp)  .. (node_{1});
\draw (58.0bp,555.0bp) node {$1$};
\draw [green,->] (node_3) ..controls (35.45bp,581.75bp) and (30.276bp,572.6bp)  .. (26.5bp,564.0bp) .. controls (22.841bp,555.67bp) and (19.562bp,546.35bp)  .. (node_{1}2);
\draw (35.0bp,555.0bp) node {$3$};
\draw [blue,->] (node_4) ..controls (6.318bp,280.01bp) and (6.6774bp,266.63bp)  .. (9.5bp,255.0bp) .. controls (10.785bp,249.7bp) and (12.887bp,244.29bp)  .. (node_5);
\draw (18.0bp,264.0bp) node {$1$};
\draw [red,->] (node_5) ..controls (28.5bp,181.19bp) and (28.5bp,164.67bp)  .. (node_{1}1);
\draw (37.0bp,167.0bp) node {$2$};
\draw [green,->] (node_6) ..controls (68.5bp,860.19bp) and (68.5bp,843.67bp)  .. (node_{1}3);
\draw (77.0bp,846.0bp) node {$3$};
\draw [red,->] (node_7) ..controls (81.461bp,569.06bp) and (77.941bp,552.34bp)  .. (node_{1});
\draw (88.0bp,555.0bp) node {$2$};
\draw [red,->] (node_8) ..controls (15.916bp,387.69bp) and (11.609bp,378.69bp)  .. (9.5bp,370.0bp) .. controls (7.5154bp,361.82bp) and (6.7485bp,352.78bp)  .. (node_4);
\draw (18.0bp,361.0bp) node {$2$};
\draw [blue,->] (node_8) ..controls (35.539bp,375.06bp) and (39.059bp,358.34bp)  .. (node_9);
\draw (49.0bp,361.0bp) node {$1$};
\draw [red,->] (node_9) ..controls (41.461bp,278.06bp) and (37.941bp,261.34bp)  .. (node_5);
\draw (49.0bp,264.0bp) node {$2$};
\draw [green,->] (node_{1}1) ..controls (28.5bp,84.188bp) and (28.5bp,67.668bp)  .. (node_{1}0);
\draw (37.0bp,70.0bp) node {$3$};
\draw [blue,->] (node_{1}2) ..controls (6.318bp,474.01bp) and (6.6774bp,460.63bp)  .. (9.5bp,449.0bp) .. controls (10.785bp,443.7bp) and (12.887bp,438.29bp)  .. (node_8);
\draw (18.0bp,458.0bp) node {$1$};
\draw [red,->] (node_{1}3) ..controls (68.5bp,763.19bp) and (68.5bp,746.67bp)  .. (node_{0});
\draw (77.0bp,749.0bp) node {$2$};
\end{tikzpicture}$}
.
\end{gather*}
All of these crystals were produced using SageMath, see \url{https://doc.sagemath.org/html/en/thematic_{T}utorials/lie/crystals.html} (last checked early 2023).
\end{Example}

The following can be found in \cite[Section 6.7 and Proposition 6.7]{BuSc-crystal-bases}.

\begin{Definition}
The \emph{column tableaux model} of the $\mathfrak{sp}_{2n}$ 
crystal with highest weight $\varpi_{k}$ is the set
\begin{gather*}
\scalebox{0.97}{$\mathcal{CT}_{\varpi_{k}}:=
\{(i_{1},\dots,i_{k})\in[1,-1]^{k}|i_{1}<i_{2}<\dots< i_{k}
\text{, and if }i_{p}=\changed{x>0},i_{q}=-x\text{ then }q-p\leq n-x\}.$}
\end{gather*}
The elements are called \emph{standard column Young diagrams} and we will write them as 
words $w$ with the leftmost entry corresponding to the top entry in Young diagram notation.

We further equip this set with the \emph{weight function}
\begin{gather*}
\wt (i_{1},\dots,i_{k})=\sum_{j=1}^{k}\wt i_{j},
\quad\text{where }\wt i=\epsilon_{i}\text{ and }\wt -i=-\epsilon_{i}.
\end{gather*}

We also equip this set with 
the \emph{crystal operators} $\ff_{i}$ and $\ee_{i}$, for $i\in\{1,\dots,n\}$, 
that act as follows. Let $i\in\{1,\dots,n-1\}$, then:
\begin{align*}
\ff_{i}\acts w
&=\begin{cases}
\text{change $i$ to $i+1$ if possible},
\\
\text{else change $-(i+1)$ to $-i$ if possible},
\\
0\text{ else},
\end{cases}
\\
\ff_{n}\acts w
&=\begin{cases}
\text{change $n$ to $-n$ if possible},
\\
0\text{ else},
\end{cases}
\end{align*}
where ``if possible'' means that the result is a standard column Young diagram.
Similarly for $\ee_{i}$, where $i\in\{1,\dots,n-1\}$, and $\ee_{n}$.
\end{Definition}

\begin{Lemma}\label{L:tableauxcyrstalishighweight}
The crystal $\mathcal{CT}_{\varpi_{k}}$ is a highest weight crystal for the 
$\mathfrak{sp}_{2n}$-modules $L(\varpi_{k})$. In particular, the crystal graph is connected.
\end{Lemma}

\begin{proof}
Well-known, see e.g. \cite[Chapter 6]{BuSc-crystal-bases}.
\end{proof}

\begin{Lemma}\label{L:crystalwelldefinedness}
Let $S=\{s_{1}<\dots<s_{k}\}$ be such that 
$|\pm L(S)|=|S_{0}|$. Then \changed{$(s_{1},\dots,s_{k})\in\mathcal{CT}_{\varpi_{k}}$}.
\end{Lemma}

\begin{proof}
View $S\subset[1,-1]$ as a dot diagram \changed{with}, as before, 
columns enumerated $1,2,\dots,n$ when reading left to right. 
The condition that $|\pm L(S)|=|S_{0}|$ 
is equivalent to the condition that for all $x\in S_{0}$, \changed{$x>0$}, 
$|S_{\geq x}|\leq n-x+1$, i.e. the number of dots to the 
right of any \changed{fully dotted column} is less than or equal to the 
index of the \changed{fully dotted column}, ensuring the presence of an 
empty column to match the $x$ column with. If $s_{p}=x$ and 
$s_{q}=-x$, then $q-p+1=|S_{\geq x}|\leq n-x+1$. Thus, 
\changed{$(s_{1},\dots,s_{k})\in\mathcal{CT}_{\varpi_{k}}$}. 
\end{proof}

The following is an alternative model of the 
$\mathfrak{sp}_{2n}$ crystal with highest weight $\varpi_{k}$, and closer 
to our previous notation:

\begin{Definition}\label{D:subsetcrystal}
We define the $\mathfrak{sp}_{2n}$ crystal with highest weight $\varpi_{k}$ to be the set: 
\begin{gather*}
\mathcal{C}_{\varpi_{k}}:=
\{S\subset [1,-1]|\,|S|=k,|\pm L(S)|=|S_{0}|\}.
\end{gather*}

The weight function we use is: 
\begin{gather*}
\wt S=\sum_{s\in S}\wt s,
\quad\text{where }\wt i=\epsilon_{i}\text{ and }\wt -i=-\epsilon_{i},
\end{gather*}
similarly as in \autoref{D:Standard}.

We also equip this set with 
the crystal operators $\ff_{i}$ and $\ee_{i}$, for $i\in\{1,\dots,n\}$, 
that act as follows. Let $i\in\{1,\dots,n-1\}$, then:
\begin{align*}
\ff_{i}\acts S&= 
\begin{cases}
(S\setminus\{i\})\cup\{i+1\} & \text{if }S_{i,i+1}=\{i\},
\\
(S\setminus\{-(i+1)\})\cup\{-i\} & \text{if }S_{i,i+1}=\{-(i+1)\},
\\
(S\setminus\{-(i+1)\})\cup\{-i\} & \text{if }S_{i,i+1}=\{i,-(i+1)\},
\\
(S\setminus\{i\})\cup\{i+1\} & \text{if }S_{i,i+1}=\{i,-i\},
\\
(S\setminus\{-(i+1)\})\cup\{-i\} & \text{if }S_{i,i+1}=\{i+1,-(i+1)\},
\\
(S\setminus\{i\})\cup\{i+1\} & \text{if }S_{i,i+1}=\{i,-(i+1),-i\},
\\
(S\setminus\{-(i+1)\})\cup\{-i\} & \text{if }S_{i,i+1}=\{i,(i+1),-(i+1)\},
\\
0 & \text{otherwise},
\end{cases}
\\
\ff_{n}\acts S&= 
\begin{cases}
(S\setminus\{n\})\cup\{-n\} & \text{if } S_{n}=\{n\},
\\
0 & \text{otherwise}.
\end{cases}
\end{align*}
Similarly for $\ee_{i}$, where $i\in\{1,\dots,n-1\}$, and $\ee_{n}$. 
\end{Definition}

\begin{Example}
The action of the crystal operators on dot diagrams is quite nice: the operator $\ff_{i}$, say for 
$i\in\{1,\dots,n-1\}$, shifts a dot in the $i$th column rightwards, if possible, a dot 
in the $(i+1)$th column leftwards, if possible, or annihilates the dot diagram. This happens in the indicated order, i.e. one first looks for dots in the $i$th column etc. For example:
\begin{gather*}
\ff_{2}\acts\,
\begin{ytableau}
*(spinach!50)\bullet & *(spinach!50)\bullet & \phantom{a} \\
\phantom{a} & \phantom{a} & *(spinach!50)\bullet
\end{ytableau}
=
\begin{ytableau}
*(spinach!50)\bullet & \phantom{a} & *(spinach!50)\bullet \\
\phantom{a} & \phantom{a} & *(spinach!50)\bullet
\end{ytableau}
\,,\quad
\ff_{2}\acts\,
\begin{ytableau}
\phantom{a} & *(spinach!50)\bullet & *(spinach!50)\bullet \\
\phantom{a} & \phantom{a} & *(spinach!50)\bullet
\end{ytableau}
=
\begin{ytableau}
\phantom{a} & *(spinach!50)\bullet & *(spinach!50)\bullet \\
\phantom{a} & *(spinach!50)\bullet & \phantom{a}
\end{ytableau}
\,,
\end{gather*}
are two examples for $n=k=3$.
\end{Example}

\begin{Lemma}\label{L:fundcrystaliso}
The map $\Phi\colon\mathcal{C}_{\varpi_{k}}\to 
\mathcal{CT}_{\varpi_{k}}$ sending $S=\{s_{1}<\dots<s_{k}\}$ to 
$(s_{1},\dots,s_{k})$ is an isomorphism of crystals. 
\end{Lemma}

\begin{proof}
That $\Phi$ is a well-defined bijection follows from 
\autoref{L:crystalwelldefinedness}. Comparing the formulas 
for actions of the crystal operators above it is easy to see 
that, setting $\Phi(0):=0$, $\Phi$ intertwines the 
action of the crystal operators.
\end{proof}

The following two lemmas follow from the theory of crystals via 
\autoref{L:fundcrystaliso}. We, for completeness, included two short arguments.

\begin{Lemma}\label{L:crystalsequenceexists}
For each $S\in\mathcal{C}_{\varpi_{k}}$, there \changed{exists $r\in\N$ and a sequence of elements $j_{1},\dots,j_{r}$ in $\{1,\dots,n\}$}, such that 
\begin{gather*}
\ff_{j_{1}}\dots\ff_{j_{r}}\acts\{1,\dots,k\}=S.
\end{gather*}
\end{Lemma}

\begin{proof}
Thanks to \autoref{L:fundcrystaliso}, the existence of the sequence $j_{1}, \dots, j_{r}$ is a consequence of the crystal graph of $\mathcal{CT}_{\varpi_{k}}$ being connected, which we established in \autoref{L:tableauxcyrstalishighweight}.
\end{proof}

\begin{Definition}
Let $S$ be as in \autoref{L:crystalsequenceexists}. 
Then we say the \emph{length} $\ell(S)$ of $S$ is $r$. 
\end{Definition}

\begin{Lemma}
The length of $S$ is well-defined. 
\end{Lemma}

\begin{proof}
Let $l_{1},\dots,l_{s}$ be a sequence in $\{1,\dots,n\}$ such that
\begin{gather*}
\ff_{l_{1}}\dots\ff_{l_{s}}\acts\{1,\dots,k\} 
=S= 
\ff_{j_{1}}\dots\ff_{j_{r}}\acts\{1,\dots,k\}.
\end{gather*}
Applying the function $\wt$ we find that
\begin{gather*}
-\sum_{i=1}^{r}\alpha_{j_{i}}+\epsilon_{1}+\dots
+\epsilon_{k}=\wt S=-\sum_{i=1}^{s}\alpha_{l_{i}}+\epsilon_{1}+\dots+\epsilon_{k}
\Rightarrow
\sum_{i=1}^{r}\alpha_{j_{i}}=\sum_{i=1}^{s}\alpha_{l_{i}}.
\end{gather*}
Since the simple roots form a basis, it follows that there is a bijection
$\{j_{1},\dots,j_{r}\}\to\{l_{1},\dots,l_{s}\}$. In particular, $r=s$.
\end{proof}

\subsection{Filtrations for \texorpdfstring{$\mathfrak{sl}_{2}$}{sl2}}\label{SS:DualFiltration}

Recall from the $\intform$-version of \autoref{SS:SLAction} that 
$\Lambda_{\intform}$ is, as an $U_{\intform}(\mathfrak{sl}_2)$-module, isomorphic to
$X_{\intform}^{\otimes n}$ where $X_{\intform}={^{1}\intform}_{\ast}\oplus{^{-1}W}_{\intform}\oplus{^{1}\intform}_{\ast}$. In particular, the operators $E^{(k)}$ and $F^{(k)}$ preserve the $\intform$ span of the vectors $v_{S}$. 

\begin{Lemma}\label{L:divided-EF-formulas}
Let $S\subset[1,-1]$ and recall $\textbf{w}$ from \autoref{N:Subsets}. We have:
\begin{gather*}
E^{(k)}\acts v_{S}=q^{\binom{k}{2}}
\sum_{\substack{\{\pm i_{1},\dots,\pm i_{k}\}\subset S_{0}^{c} \\ i_{1}<\dots<i_{k}}} (-q)^{\sum_{t=1}^{k}\textbf{w}_{>i_{t}}(S)}\cdot v_{S\cup\{\pm i_{1},\dots,\pm i_{k}\}},
\\
F^{(k)}\acts v_{S}=q^{\binom{k}{2}}
\sum_{\substack{\{\pm i_{1},\dots,\pm i_{k}\}\subset S_{0} 
\\ i_{1}<\dots<i_{k}}} 
(-q)^{\sum_{t=1}^{k}\textbf{w}_{<i_{t}}(S)}\cdot v_{S\setminus\{\pm i_{1},\dots,\pm i_{k}\}}.
\end{gather*}
\end{Lemma}

\begin{proof}
We prove the claim for $E^{(k)}$. A similar argument works for $F^{(k)}$. 

Let $\{\pm i_{1},\dots,\pm i_{k}\}\subset S_{0}^{c}$, with \changed{$0\leq i_{1}<\dots<i_{k}$}. Note that 
\begin{gather*}
\textbf{w}_{>j}(S\cup\{\pm i_{1},\dots,\pm i_{k}\}) 
=\textbf{w}_{>j}(S)+2\cdot|\{i_{p}\,|\,i_{p}>j\}|.
\end{gather*}
For \changed{$\sigma$ in the symmetric group \changed{$\mathfrak{S}_{k}$} on $\{1,\dots,k\}$}, define 
\begin{gather*}
\textbf{w}_{>\sigma}(S):=\textbf{w}_{>i_{\sigma(1)}}(S)+\textbf{w}_{>i_{\sigma(2)}}
(S\cup\{\pm i_{\sigma(1)}\} 
+\dots+\textbf{w}_{>i_{\sigma(k)}}(S\cup\{\pm i_{\sigma(1)},\dots,\pm i_{\sigma(k-1)}\}).
\end{gather*}
Then, if we write $\ell(\sigma)$ for the length of the permutation $\sigma$, i.e. the number of inversions in $\sigma$, then $\textbf{w}_{>\sigma}(S)=\textbf{w}_{>\id}(S)+2\cdot\ell(\sigma)$.

Repeated application of \autoref{L:EFops} shows that 
we have
\begin{align*}
\sum_{\sigma\in \changed{\mathfrak{S}_{k}}}(-q)^{\textbf{w}_{>\sigma}(S)} 
=(-q)^{\textbf{w}_{>\id}(S)}\sum_{\sigma\in \changed{\mathfrak{S}_{k}}}(-q)^{2\ell(\sigma)}
=(-q)^{\textbf{w}_{>\id}(S)}q^{\binom{k}{2}}[k]_{\intform}!
\end{align*}
as the the coefficient of 
$v_{S\cup\{\pm i_{1},\dots,\pm i_{k}\}}$ in $E^{k}\acts v_{S}$.
This implies:
\begin{align*}
E^{(k)}\acts v_{S} 
&=q^{\binom{k}{2}}
\sum_{\substack{\{\pm i_{1},\dots,\pm i_{k}\}\subset S_{0}^{c} 
\\ i_{1}<\dots<i_{k}}}
(-q)^{\sum_{t=1}^{k}\wt_{>i_{t}}(S\cup\{\pm i_{1},\dots,\pm i_{t-1}\})}
\cdot v_{S\cup\{\pm i_{1},\dots,\pm i_{k}\}}
\\
&=q^{\binom{k}{2}}
\sum_{\substack{\{\pm i_{1},\dots,\pm i_{k}\}\subset S_{0}^{c} 
\\ i_{1}<\dots<i_{k}}}(-q)^{\sum_{t=1}^{k}\wt_{> i_{t}}(S)}\cdot 
v_{S\cup\{\pm i_{1},\dots,\pm i_{k}\}}.
\end{align*}
This proves the statement.
\end{proof}

The following can be, for example, found in \cite[Section 5.2]{Ja-lectures-qgroups}. 
With a bit \changed{of} care, cf. \cite[Chapter 5]{Lu-quantumgroups-book},
these operators can be defined and used for type-1 and type-(-1) 
$U_{q}(\mathfrak{sl}_{2})$-modules.

\begin{Definition}\label{D:quantum-Weyl-group}
Let $M$ be a finite dimensional 
type-$\epsilon$ $U_{q}(\mathfrak{sl}_{2})$-module. 
We define \emph{Lusztig's quantum Weyl group element} 
$T$ as the operator $T\colon M\to M$, which acts 
on $v\in M[k]=\{v\in M\,|\,K\acts v=\epsilon q^{k}\cdot v\}$ as
\begin{gather*}
T\acts v=
\sum_{\substack{a,b,c\geq 0\\-a+b-c=k}}
(-1)^{b}q^{b-ac}\cdot E^{(a)}F^{(b)}E^{(c)}\acts v.
\end{gather*}
(Note that these sums contain only finitely many terms, 
since $M$ is finite dimensional so $E$ and $F$ act nilpotently.)
\end{Definition}

\begin{Lemma}\label{L:triple-exponent-iso}
Let $M$ be a finite dimensional type-$\epsilon$ $U_{q}(\mathfrak{sl}_{2})$-module. Then 
\begin{gather*}
T|_{M[k]}\colon M[k]\xrightarrow{\cong}M[-k]
\end{gather*}
is an $\field$-linear isomorphism. The inverse of $T$ is the operator 
$T^{-1}\colon M\to M$, which acts on $w\in M[k]$ as
\begin{gather*}
T^{-1}\acts w=
\sum_{\substack{a,b,c\geq 0 \\ a-b+c=m}}
(-1)^{b}q^{ac-b}\cdot F^{(a)}E^{(b)}F^{(c)}\acts w.
\end{gather*}
\end{Lemma}

\begin{proof}
See \cite[Section 8.4]{Ja-lectures-qgroups} or \cite[Proposition 5.2.7]{Lu-quantumgroups-book}.
\end{proof}

\begin{Lemma}\label{L:intquantumWeyliso}
Suppose that $M$ is a finite dimensional 
type-$\epsilon$ $U_{q}(\mathfrak{sl}_{2})$-module 
and that $M_{\intform}\subset M$ is preserved 
by $U_{\intform}(\mathfrak{sl}_{2})$ and 
satisfies $M_{\intform}[k]=M_{\intform}\cap M[k]$ for all $k\in\Z$. Then 
\begin{gather*}
T\changed{|}_{M_{\intform}[k]}\colon M_{\intform}[k]\xrightarrow{\cong}M_{\intform}[-k]
\end{gather*}
is an $\intform$-linear isomorphism.
\end{Lemma}

\begin{proof}
Since $E^{(i)}$ and $F^{(i)}$ preserve $M_{\intform}$, 
it follows that $T$ preserves $M_{\intform}$. Then, 
noting that $T$ maps $M[k]$ to $M[-k]$ and 
$M_{\intform}[\pm k]= M_{\intform}\cap M[\pm k]$, 
we see $T$ maps $M_{\intform}[k]$ to $M_{\intform}[-k]$. 
Similarly, we find that $T^{-1}$ maps $M_{\intform}[-k]$ to $M_{\intform}[k]$. 
\end{proof}

\begin{Proposition}\label{P:sl2tilting}
The $U_{\intform}(\mathfrak{sl}_{2})$-module $\Lambda_{\intform}$ is a tilting module. In particular, $\Lambda_{\intform}$ has a Weyl module filtration with Weyl character
\begin{gather*}
[\Lambda_{\intform}]= 
\sum_{k=0}^{n}\mathrm{rk}_{\intform}\Delta_{\intform}(\varpi_{k})[{^{\epsilon}\Delta_{\intform}(n-k)}],
\end{gather*}
where $\epsilon=1$, if $k$ is even, and $\epsilon=-1$, if $k$ is odd.
\end{Proposition}

\begin{proof}
Because we have that ${^{1}\intform}_{\ast}$ and ${^{-1}W}_{\intform}$ 
are tilting $U_{\intform}(\mathfrak{sl}_{2})$-modules (this is well-known 
and easy to check by hand), the claim that $\Lambda_{\intform}$ is a tilting module, such that $\Lambda^{even}$ is a type-1 module and $\Lambda^{odd}$ is a type-(-1) module, follows 
directly from $\Lambda_{\intform}\cong X_{\intform}^{\otimes n}$ and the fact that 
tensor products of tilting modules are tilting, see \autoref{L:TensorTilting}. To compute the Weyl character of $\Lambda_{\intform}$ it suffices to determine the character of $\Lambda_{q}$, which is known by \autoref{T:QuantumHowe}.
\end{proof}

\subsection{Embedding fundamental Weyl modules in exterior powers}

\begin{Notation}
Let $B^{k}_{\Delta(\varpi_{k})}=\{b_{S}\,|\,|S|=k\text{ and }|\pm L(S)|=|S_{0}|\}$.
\end{Notation}

\begin{Lemma}\label{L:generatingker}
For each $b_{S}\in B_{\Delta(\varpi_{k})}^{k}$ there is an element $u_{S}\in U_{\intform}(\mathfrak{sp}_{2n})$ such that $u_{S}\acts v_{\{1,\dots,k\}}=b_{S}$, that is
\begin{gather*}
u_{S}\acts
\begin{ytableau}
*(spinach!50)\bullet & \dots & *(spinach!50)\bullet & \phantom{a} & \dots\\
\phantom{a} & \phantom{a} & \phantom{a} & \phantom{a} & \dots
\end{ytableau}
=
b_{S},
\end{gather*}
where the rightmost dot is in the $k$th column.
\end{Lemma}

\begin{proof}
Let $b_{S}\in B_{\Delta(\varpi_{k})}^{k}$, i.e. 
$S\in\mathcal{C}_{\varpi_{k}}$, such that $\ell(S)=r$. 

If $r=0$, then $S=\{1,\dots,k\}$, $1\in U_{\intform}(\mathfrak{sp}_{2n})$, 
and $1\acts v_{\{1,\dots,k\}}=b_{\{1,\dots,k\}}$. 
Suppose inductively that the result is true for all 
$b_{T}\in B_{\Delta(\varpi_{k})}^{k}$ such that $\ell(T)<r$.

Let $r>0$. There is a sequence of elements $j_{1},\dots, j_{r}$ in 
$\{1,\dots,n\}$ such that 
\begin{gather*}
\ff_{j_{1}}\dots\ff_{j_{r}}\acts\{1,\dots,k\} = S. 
\end{gather*}
Note that from \autoref{D:subsetcrystal} we have two options.
The first option is that $j_{1}\in\{1,\dots,n-1\}$ and 
$(\ff_{j_2}\dots\ff_{j_{r}}\acts\{1,\dots,k\})_{j_{1},j_{1}+1}$ 
is one of the following:
\begin{gather*}
\{j_{1}\},\{-(j_{1}+1)\},
\{j_{1}, -(j_{1}+1)\},\{ \pm j_{1}\},
\{\pm(j_{1}+1)\},\{\pm j_{1},-(j_{1}+1)\},
\{j_{1}, \pm(j_{1}+1)\}.
\end{gather*}
Or, as a second option, $j_{1}=n$ and 
$(\ff_{j_2}\dots\ff_{j_{r}}\acts\{1,\dots,k\})_{j_{1}}=\{n\}$.

By induction, there is an element $u_{\ff_{j_2}\dots
\ff_{j_{r}}\acts\{1,\dots,k\}}\in U_{\intform}(\mathfrak{sp}_{2n})$ such that
\begin{gather*}
u_{\ff_{j_2}\dots\ff_{j_{r}}\acts\{1,\dots,k\}}
\acts v_{\{1,\dots,k\}}=b_{\ff_{j_2}\dots\ff_{j_N}\acts\{1,\dots,k\}}.
\end{gather*}
Excluding the $\{\pm j_{1}\}$ case, we can compare the formulas in \autoref{D:subsetcrystal} and \autoref{L:IntExtAction} to deduce that
\begin{gather*}
F_{j_{1}}b_{\ff_{j_2}\dots\ff_{j_{r}}\acts\{1,\dots,k\}}
=b_{\ff_{j_{1}}\ff_{j_2}\dots\ff_{j_{r}}\acts\{1,\dots,k\}}.
\end{gather*}
Thus, we have $F_{j_{1}}u_{\ff_{j_2}\dots\ff_{j_{r}}\acts\{1,\dots,k\}}
\in U_{\intform}(\mathfrak{sp}_{2n})$ such that 
\begin{gather*}
F_{j_{1}}u_{\ff_{j_2}\dots\ff_{j_{r}}\acts\{1,\dots,k\}}\acts v_{\{1,\dots,k\}}=b_{s}.
\end{gather*}

It remains to consider the case where $j_{1}\in\{1,\dots,n-1\}$ 
and $(\ff_{j_2}\dots\ff_{j_{r}}\acts\{1,\dots,k\})_{j_{1}, j_{1}+1} 
=\{\pm j_{1}\}$. It follows that
\begin{gather*}
(\ff_{j_2}\dots\ff_{j_{r}}\acts\{1,\dots,k\})\setminus\{-j_{1}\} 
\cup\{-(j_{1}+1)\}\in\mathcal{C}_{\varpi_{k}},
\end{gather*}
and therefore there is a sequence $l_{1},\dots,l_{s}$ in $\{1,\dots,n\}$ such that 
\begin{gather*}
\ff_{l_{1}}\dots\ff_{l_{s}}\acts\{1,\dots,k\} 
=\big((\ff_{j_2}\dots\ff_{j_{r}}\acts\{1,\dots,k\})\setminus\{-j_{1}\}\big)\cup\{-(j_{1}+1)\}.
\end{gather*}
Then from \autoref{D:subsetcrystal} we find
\begin{gather*}
\ff_{j_{1}}\ff_{j_{1}}\ff_{l_{1}}\dots\ff_{l_{s}}\acts\{1,\dots, k\} = S. 
\end{gather*}
In particular, we have $2+s=\ell(S)=r$, so $s<r$. By induction, there is an element $u_{\ff_{l_{1}}\dots\ff_{l_{s}}\acts\{1,\dots,k\}}\in U_{\intform}(\mathfrak{sp}_{2n})$ such that $u_{\ff_{l_{1}}\dots\ff_{l_{s}}\acts\{1,\dots,k\}}\acts v_{\{1,\dots,k\}} = b_{\ff_{l_{1}}\dots\ff_{l_{s}}\acts\{1,\dots,k\}}$. Comparing the formulas in \autoref{D:subsetcrystal} and \autoref{L:IntExtAction}, we find
\begin{gather*}
f_{j_{1}}f_{j_{1}}u_{\ff_{l_{1}}\dots\ff_{l_{s}}\acts\{1,\dots,k\}}\acts v_{\{1,\dots,k\}}= [2]_{\intform}b_{\ff_{j_{1}}\ff_{j_{1}}\ff_{l_{1}}\dots\ff_{l_{s}}\acts\{1,\dots,k\}}. 
\end{gather*}
Thus, we have $f_{j_{1}}^{(2)}u_{\ff_{l_{1}}\dots\ff_{l_{s}}\acts\{1,\dots,k\}}\in U_{\intform}(\mathfrak{sp}_{2n})$ and
\begin{gather*}
f_{j_{1}}^{(2)}u_{\ff_{l_{1}}\dots\ff_{l_{s}}\acts\{1,\dots,k\}}\acts v_{\{1,\dots,k\}}=b_{s}. 
\end{gather*}
This completes the proof.
\end{proof}

Recall that $F$ denotes the 
Chevalley generator from $U_{\intform}(\mathfrak{sl}_{2})$.

\begin{Lemma}\label{L:caninker}
If $b_{S}\in B^{k}_{\Delta(\varpi_{k})}$, then $F\acts b_{S}=0$. 
\end{Lemma}

\begin{proof}
From \autoref{L:EFops}, we have $F\acts v_{\{1,\dots,k\}}=0$. 
Let $b_{S}\in B^{k}_{\Delta(\varpi_{k})}$. By 
\autoref{L:generatingker}
there is $u_{S}$ such that $b_{S}=u_{S}\acts v_{\{1,\dots,k\}}$. It follows that 
\begin{gather*}
F\acts b_{S}= 
F\acts(u_{S}\acts v_{\{1,\dots,k\}})= 
u_{S}\acts (F\acts v_{\{1,\dots,k\}})= 
u_{S}\acts 0=0,
\end{gather*}
which uses \autoref{P:IntActionExt}.(c).
\end{proof}

\begin{Lemma}\label{L:WeylequalskeroverK}
We have $\field B_{\Delta(\varpi_{k})}^{k}= 
\ker F|_{\Lambda_{q}^{k}}$ as $U_{q}(\mathfrak{sp}_{2n})$-modules.
\end{Lemma}

\begin{proof}
\autoref{L:caninker} and \autoref{C:sskernelisisimple} implies that
\begin{gather*}
\field B_{\Delta(\varpi_{k})}^{k}\subset\ker F|_{\Lambda_{q}^{k}} 
\xrightarrow{\cong}L_{q}(\varpi_{k}).
\end{gather*}
Moreover, from \autoref{D:subsetcrystal} and \autoref{L:fundcrystaliso}, one can observe that the character of $\field B_{\Delta(\varpi_{k})}^{k}$ is the same as the character of $L_{q}(\varpi_{k})$. Hence, the inclusion $\field B_{\Delta(\varpi_{k})}^{k}\to L_{q}(\varpi_{k})$ is an isomorphism.
\end{proof}

Let $X=X(v\to b)$ be the matrix with rows indexed by 
$\{T\subset[1,-1]\,|\,|T|=k\}$ and columns indexed 
by $\{S\subset[1,-1]\,|\,|S|=k,|R(S)|=|S_{0}|\}$ such 
that the $(T,S)$ entry is the coefficient of $v_{T}$ in $b_{S}$.
The matrix $X$ is a type of \emph{base change matrix} 
between the standard and the canonical basis.

\begin{Lemma}\label{L:leftinverseinA}
There is an $\intform$-matrix $Y=Y(b\to v)$ with the opposite indexing 
when compared to $X$, such that $YX=\id$.
\end{Lemma}

\begin{proof}
If we partially order the subsets of $[1,-1]$ by how far fully dotted columns
are to the right, then 
$X$ is a unitriangular matrix with values in $\intform$. The lemma follows.
\end{proof}

\begin{Remark}
The base change argument in \autoref{L:leftinverseinA} is 
often used in the theory of crystals,
see for example 
\cite[Proof of Lemma 4.5]{MiMi-affine-sln-crystal} for 
an early reference from where we got the idea from.
\end{Remark}

\begin{Lemma}\label{L:canbasisspanker}
We have $\intform B_{\Delta(\varpi_{k})}^{k}= 
\ker F|_{\Lambda_{q}^{k}}$ as $U_{\intform}(\mathfrak{sp}_{2n})$-modules.
\end{Lemma}

\begin{proof}
Because of \autoref{L:caninker}, it suffices to show that $\ker F|_{\Lambda_{\intform}^{k}}$ is contained in $\intform B_{\Delta(\varpi_{k})}^{k}$. Let 
$x\in\ker F|_{\Lambda_{\intform}^{k}}$. Then in particular,
\begin{gather*}
x=\sum_{T\subset[1,-1],|T|=k}a_{T}\cdot v_{T},\quad a_{T}\in\intform.
\end{gather*}
Viewing $x\in \ker F|_{\Lambda_{q}^{k}}=\field B_{\Delta(\varpi_{k})}$, with the equality coming from 
\autoref{L:WeylequalskeroverK}, we can write
\begin{gather*}
x=\sum_{\substack{S\subset[1,-1],|S|=k\\ |L(S)|=|S_{0}|}}k_{S}\cdot b_{S},\quad k_{S}\in\field.
\end{gather*}
It follows from \autoref{L:leftinverseinA} that
\begin{gather*}
k_{S}=\sum_{T\subset[1,-1],|T|=k}Y_{S,T}a_{T}\in\intform,
\end{gather*}
where $Y$ is the matrix from \autoref{L:leftinverseinA}. 
Hence, we get $x\in\intform B_{\Delta(\varpi_{k})}^{k}$.
\end{proof}

Let $v_{\varpi_{k}}^{+}$ denote the highest weight vector of 
$\Delta_{\intform}(\varpi_{k})$ (this is unique up to some choice 
of scalars in $\intform$; the choice 
does not play any role for us)

\begin{Proposition}\label{P:kerisWeyl}
There is an isomorphism of $U_{\intform}(\mathfrak{sp}_{2n})$-modules
\begin{gather*}
\Delta_{\intform}(\varpi_{k})\to
\ker F|_{\Lambda_{\intform}^{k}}
\end{gather*}
such that $v_{\varpi_{k}}^{+}\mapsto v_{\{1,\dots,k\}}$.
\end{Proposition}

\begin{proof}
Write $L_{q}(\varpi_{k})=\field\otimes\Delta_{\intform}(\varpi_{k})$, and continue to denote the highest weight vector by $v_{\varpi_{k}}^{+}$. By \autoref{C:sskernelisisimple} we have
\begin{gather*}
v_{\{1,\dots,k\}}\in 
\ker F|_{\Lambda^{k}_{\intform}}\subset 
\ker F|_{\Lambda^{k}_{q}}\xrightarrow{\cong}L_{q}(\varpi_{k})
\end{gather*}
such that $v_{\{1,\dots,k\}}\mapsto v_{\varpi_{k}}^{+}$.

Thanks to \autoref{P:IntActionExt}, \autoref{L:kernelispreserved} implies that $\ker F|_{\Lambda_{\intform}^{k}}$ is preserved by $U_{\intform}(\mathfrak{sp}_{2n})$. \autoref{L:canbasisspanker} and \autoref{L:generatingker} imply that $\ker F|_{\Lambda_{\intform}^{k}}$ is generated over $U_{\intform}(\mathfrak{sp}_{2n})$ by $v_{\{1,\dots,k\}}$. The result follows from the description of Weyl modules in \autoref{L:Weylmodule-omnibuslemma}.
\end{proof}

\subsection{Filtrations for $\mathfrak{sp}_{2n}$}

We now aim to show that the exterior algebra $\Lambda_{\intform}$ is 
a tilting $U_{\intform}(\mathfrak{sp}_{2n})$-module.

\begin{Lemma}\label{L:kthPower}
Let $0\leq k\leq n$ and let $i\geq 0$ be minimal such that 
$k-2i< 0$. Then $\Lambda_{\intform}^{k}\cong\ker F^{(i)}|_{\Lambda_{\intform}^{k}}$ 
as $U_{\intform}(\mathfrak{sp}_{2n})$-modules.
\end{Lemma}

\begin{proof}
This follows since $F$ removes a fully dotted column from dot diagrams, see \autoref{D:EFOps}.
\end{proof}

Note that \autoref{L:kthPower} implies that
$\Lambda_{\intform}^{k}$ has the 
following filtration by $U_{\intform}(\mathfrak{sp}_{2n})$-submodules:
\begin{gather*}
0=\ker F^{(0)}|_{\Lambda_{\intform}^{k}} 
\subset 
\ker F^{(1)}|_{\Lambda_{\intform}^{k}}
\subset\dots\subset 
\ker F^{(i-1)}|_{\Lambda_{\intform}^{k}}\subset
\ker F^{(i)}=\Lambda_{\intform}^{k}.
\end{gather*}
Our goal is to show that this is a Weyl filtration.

\begin{Notation}\label{N:almosthwvectors}
Let $0\leq k\leq n$ and let $i\geq 0$ such that $k-2i\geq 0$. Then we write
\begin{gather*}
v_{k,k-2i}
\leftrightsquigarrow
\begin{ytableau}
*(spinach!50)\bullet & \dots & *(spinach!50)\bullet & \phantom{a} & \dots & \phantom{a} & *(spinach!50)\bullet & \dots & *(spinach!50)\bullet
\\
\phantom{a} & \phantom{a} & \phantom{a} & \phantom{a} & \dots & \phantom{a} & *(spinach!50)\bullet & \dots & *(spinach!50)\bullet
\end{ytableau}
\,,
\end{gather*}
with the single dots in columns $1$ to $k-2i$ and the \changed{fully dotted columns} $(n-i+1)$ to $n$.
Or in formulas, $v_{k,k-2i}:=v_{\{1,\dots, k-2i\}\cup\{\pm (n-i+1),\dots,\pm n\}}$.
\end{Notation}

\begin{Remark}
A consequence of \autoref{P:kerisWeyl} is that $v_{k,k}=v_{\{1,\dots,k\}}$ is a 
highest weight vector generating a $U_{\intform}(\mathfrak{sp}_{2n})$-submodule 
isomorphic to $\Delta_{\intform}(\varpi_{k})$. 
However, for $i>0$, the vector $v_{k,k-2i}$ is not a singular vector.
\end{Remark}

\begin{Remark}
If $x\in\ker F^{(i+1)}|_{\Lambda_{\intform}^{k}}$, then 
$F^{(i)}\acts x\in\Lambda^{k-2i}_{\intform}$ and 
$F\acts(F^{(i)}\acts x)=[i+1]_{\intform}\cdot 
F^{(i+1)}\acts x=0$. Combining this observation with 
\autoref{P:IntActionExt}, we see that $F^{(i)}$ induces a map
$F^{(i)}\acts\placeholder\in\Hom_{U_{\intform}(\mathfrak{sp}_{2n})}(\ker F^{(i+1)}|_{\Lambda_{\intform}^{k}},\ker F|_{\Lambda_{\intform}^{k-2i}})$.
\end{Remark}

\begin{Lemma}\label{L:Fi-mapspartialhwtohw}
Let $0\leq k\leq n$ and let $i\geq 0$ such that $k-2i\geq 0$. Then
\begin{gather*}
F^{(i)}\acts v_{k,k-2i}=\xi\cdot v_{k-2i,k-2i},
\end{gather*}
where $\xi=\pm q^{j}$ for some $j\in\Z$. In particular, $v_{k,k-2i}\in\ker F^{(i+1)}|_{\Lambda_{\intform}^{k}}$.
\end{Lemma}

\begin{proof}
It follows from \autoref{L:divided-EF-formulas} that
$F^{(i)}(v_{k,k-2i})=q^{\binom{i}{2}}(-q)^{i(n-k+i)+\binom{i}{2}}\cdot v_{k-2i,k-2i}$.
\end{proof}

\begin{Lemma}\label{L:Fi-is-surjective}
There is a surjection
\begin{gather*}
F^{(i)}\acts\placeholder\colon\ker F^{(i+1)}|_{\Lambda_{\intform}^{k}}
\twoheadrightarrow\ker F|_{\Lambda_{\intform}^{k-2i}}.
\end{gather*}
\end{Lemma}

\begin{proof}
By \autoref{L:Fi-mapspartialhwtohw} the vector 
$v_{k,k-2i}\in\ker F^{(i+1)}|_{\Lambda_{\intform}^{k}}$ 
and $F^{(i)}\acts v_{k,k-2i}=\xi\cdot v_{k-2i,k-2i}$ for $\xi=\pm q^{j}$ for some $j\in\Z$. Let $x\in \ker F|_{\Lambda_{\intform}^{k-2i}}$. 
It follows from \autoref{P:kerisWeyl} that there is $u\in U_{\intform}(\mathfrak{sp}_{2n})$ such that $u\acts v_{k-2i, k-2i} = x$. Noting that $\xi^{-1}\in \intform$, we can apply \autoref{L:kernelispreserved} to see that $\xi^{-1}\cdot u\acts v_{k, k-2i}\in \ker F^{(i+1)}|_{\Lambda_{\intform}^k}$. Since
\[
F^{(i)}\acts (\xi^{-1}\cdot u\acts v_{k, k-2i}) = \xi^{-1}\cdot u\acts (F^{(i)}\acts v_{k, k-2i}) = x,
\]
the lemma follows.
\end{proof}

We now have the desired Weyl filtration.

\begin{Proposition}\label{P:spfilt}
Let $k\in\{1,\dots,n\}$. We have an isomorphism of 
$U_{\intform}(\mathfrak{sp}_{2n})$-modules
\begin{gather*}
\ker F^{(i+1)}|_{\Lambda_{\intform}^{k}}/\ker F^{(i)}|_{\Lambda_{\intform}^{k}} 
\cong\Delta_{\intform}(\varpi_{k-2i}).
\end{gather*}
In particular, the $U_{\intform}(\mathfrak{sp}_{2n})$-module $\Lambda_{\intform}^{k}$ has a 
filtration by Weyl modules with Weyl character
\begin{gather*}
[\Lambda_{\intform}^{k}] 
=\sum_{k-2i\geq 0}[\Delta_{\intform}(\varpi_{k-2i})].
\end{gather*}
(The final equation holds in the Grothendieck ring.)
\end{Proposition}

\begin{proof}
By \autoref{L:Fi-is-surjective} and \autoref{P:kerisWeyl}, we have
\begin{gather*}
\ker F^{(i+1)}|_{\Lambda_{\intform}^{k}}/\ker F^{(i)}|_{\Lambda_{\intform}^{k}}
\cong\ker F|_{\Lambda_{\intform}^{k-2i}}\cong\Delta_{\intform}(\varpi_{k-2i}).
\end{gather*}
Suppose $j$ is such that $k-2j< 0$ and $k-2j+2\geq 0$. Then 
$\Lambda^{k}\subset\ker F^{(j)}$ and the filtration
\begin{gather*}
0=\ker F^{(0)}|_{\Lambda_{\intform}^{k}} 
\subset 
\ker F^{(1)}|_{\Lambda_{\intform}^{k}}
\subset\dots\subset 
\ker F^{(j-1)}|_{\Lambda_{\intform}^{k}}\subset
\ker F^{(j)}=\Lambda_{\intform}^{k}
\end{gather*}
has subquotients isomorphic to 
\begin{gather*}
\Delta_{\intform}(\varpi_{k}),\Delta_{\intform}(\varpi_{k-2}),\dots,\Delta_{\intform}(\varpi_{k-2j+2}).
\end{gather*}
Thus, $\Lambda_{\intform}^{k}$ has Weyl character as claimed in the proposition.
\end{proof}

\begin{Proposition}\label{P:sptilting}
The $U_{\intform}(\mathfrak{sp}_{2n})$-module $\Lambda_{\intform}$ is a tilting module. 
\end{Proposition}

\begin{proof}
Let $0\le k\le n$. \autoref{P:spfilt} implies that $\Lambda_{\intform}^k$ has a Weyl filtration. Thanks to \autoref{L:twist-sigma-equivalence}, \autoref{R:twist-intform-iso}, and \autoref{P:SelfDual}, we find that $\Lambda_{\intform}^k\cong ({^{\omega}\Lambda_{\intform}^k})^*$, so $\Lambda^k$ also has a dual Weyl filtration. 

Since $U_{\intform}(\mathfrak{sl}_{2})$ acts on $\Lambda_{\intform}$, with weight spaces $\Lambda_{\intform}^{i}$, such that $0\leq i\leq 2n$, \autoref{L:intquantumWeyliso} implies that $T|_{\Lambda_{\intform}^{k}}$ is an isomorphism $\Lambda_{\intform}^{k}\xrightarrow{\cong}\Lambda_{\intform}^{2n-k}$. Since $T$ is an $\intform$-linear combination of elements in $U_{\intform}(\mathfrak{sl}_{2})$, \autoref{P:IntActionExt} implies that 
\begin{gather*}
T\in\End_{U_{\intform}(\mathfrak{sp}_{2n})}(\Lambda_{\intform}).
\end{gather*}
Thus, $\Lambda_{\intform}^{i}$ is a tilting module for $0\leq i\leq 2n$, and since a direct sum of tilting modules is a tilting module $\Lambda_{\intform}=\bigoplus_{i=0}^{2n}\Lambda_{\intform}^{i}$ is a tilting module.
\end{proof}

\section{Howe duality -- part II}\label{S:IntHowe}

We now prove \autoref{T:IntroMain}. We start below by collecting 
well-known material which we then apply to get new results.

\subsection{Some generalities, e.g. based modules and Weyl filtrations}

We fix a finite dimensional simple Lie algebra $\mathfrak{g}$. There are associated weights $X$, roots $\Phi$, positive roots $\Phi_{+}$, dominant weights $X_{+}$, and 
the Weyl group $W$ which acts on $X$ and has longest element $w_{0}$.

\cite[Chapter 24]{Lu-quantumgroups-book} defines a non-unital $\field$-algebra, 
the \emph{idempotented form} $\dot{U}_{q}(\mathfrak{g})$, with mutually orthogonal idempotents 
$\pr_{\lambda}$ for $\lambda\in X$, such that 
\begin{gather*}
\dot{U}_{q}(\mathfrak{g}) 
=\bigoplus_{\lambda,\mu\in X}\pr_{\mu}\dot{U}_{q}(\mathfrak{g})\pr_{\lambda},
\end{gather*}
and a divided powers form $\dot{U}_{\intform}(\mathfrak{g})$ which is generated as an $\intform$-algebra by $E_{\alpha}^{(m)}\pr_{\lambda}$ and $F_{\alpha}^{(m)}\pr_{\lambda}$. 

\cite[Chapter 25]{Lu-quantumgroups-book} describes a way to extend the canonical basis for 
the positive part $U_{q}^{+}(\mathfrak{g})$ to a basis $\dcanB$ of $\dot{U}_{q}(\mathfrak{g})$ such that
\begin{gather*}
\intform\dcanB=\dot{U}_{\intform}(\mathfrak{g}).
\end{gather*}
Moreover, if $b\in\dcanB$, then $b\in\dot{U}_{q}(\mathfrak{g})\pr_{\lambda}$ 
for a unique $\lambda\in X$.

The following can be found in \cite{Lu-quantumgroups-book}, and is nicely summarized in 
\cite[Section 1]{Ka-Based-Filt} (we however drop the condition on being finite dimensional 
from those papers):

\begin{Definition}
\changed{Let} $(M,B,\phi_{M})$ be a triple consisting of a (finite dimensional) type-1 
$\dot{U}_{q}(\mathfrak{g})$-module $M$, a $\ring$-basis $B$ of $M$ and 
an involution $\phi_{M}\colon M\to M$.
The triple is a \emph{based 
$\dot{U}_{\intform}(\mathfrak{g})$-module} if 
$B=\bigcup_{v\in X}(B\cap M_{v})$ (weight spaces), $M_{\intform}=\intform B$ is 
$\dot{U}_{\intform}(\mathfrak{g})$-stable, $\phi_{M}(b)=b$ for $b\in B$, 
$\phi_{M}(u\acts m)=\phi(u)\acts\phi_{M}(m)$ with $\phi$ swapping $K_{i}$ and $K^{-1}_{i}$, 
and the pair $(L(M)=\Q(q)_{(q)}B,\overline{B})$ is a crystal base of $M$ for $\overline{B}$ the image of $B$ in $L(M)/qL(M)$.
\end{Definition} 
 
\begin{Example}
The left regular $\dot{U}_{q}(\mathfrak{g})$-module $\dot{U}_{q}(\mathfrak{g})$ is a based module with basis $\dcanB$. \changed{The involution used in this example is as in \cite[Section 11.9]{Ja-lectures-qgroups} (essentially, send $e_{i}$ to $e_{i}$, $f_{i}$ to $f_{i}$, and invert $k_{i}$ and $q$).}
\end{Example}

\begin{Example}
Weyl modules are prototypical examples of based modules. That is,
each Weyl module $\Delta_{\intform}(\lambda)$, 
where $\lambda\in X_{+}$, comes with a distinguished 
highest weight vector $v_{\lambda}^{+}$ and canonical 
basis $\canB_{\Delta(\lambda)}$. \changed{The involution is again as in \cite[Section 11.9]{Ja-lectures-qgroups}.} There is a unique canonical basis element 
$v_{\lambda}^{-}\in\Delta(\lambda)$ which lies in the $w_{0}(\lambda)$ 
weight space.
\end{Example}

For based modules we have the notion of tensor product (that works similarly 
as for crystals) \cite[Section 27.3]{Lu-quantumgroups-book}. In particular, for each $\lambda,\mu\in X_{+}$, there is a canonical basis $\canB_{\Delta(\lambda)\otimes\Delta(\mu)}$ of $\Delta_{\intform}(\lambda)\otimes\Delta_{\intform}(\mu)$. \changed{The involution is described in \cite[Theorem 24.3.3]{Lu-quantumgroups-book}.}

\changed{\begin{Notation}
For $\lambda\in X_{+}$, we write
$M[\lambda]$ to denote the sum of submodules of $M$ isomorphic to $\Delta(\lambda)$, and write $M[\ge \lambda]:=\oplus_{\mu\ge  \lambda}M[\mu]$.  
\end{Notation}}

\changed{\begin{Definition}
For $b\in B$, there is a unique $\mu \in X_+$ which is maximal (with respect to the dominance order on $X_+$) such that $b\in M[\ge \mu]$ \cite[Section 27.2.1]{Lu-quantumgroups-book}. This defines a map $B\rightarrow X_+$. Denote by $B[\mu]$ the fiber of this map over $\mu\in X_+$. 
\end{Definition}}

Recall the following from \cite[Section 29.1.1]{Lu-quantumgroups-book}:
Let $b\in\dcanB\cap\dot{U}_q(\mathfrak{g})\pr_{\zeta}$. 
Choose $\lambda\in X_{+}$ such that $(\alpha_{i}^{\vee},\lambda)$ is 
``large enough'', see \cite[Section 25.2]{Lu-quantumgroups-book}, for all $i$. Then $b\cdot v_{-w_0(\lambda)}^{-}\otimes v_{\lambda+ \xi}^{+}$ is in the canonical basis $\canB_{\Delta(-w_0(\lambda))\otimes\Delta(\lambda+\xi)}$, 
and by \cite[Section 27.2.1]{Lu-quantumgroups-book} there is a unique $\mu\in X_{+}$ such that 
\begin{gather*}
b\cdot v_{-w_0(\lambda)}^{-}\otimes v_{\lambda+\xi}^{+}\in 
\canB_{\Delta(-w_0(\lambda))\otimes\Delta(\lambda+\xi)}[\mu].
\end{gather*}
The weight $\mu\in X_{+}$ does not depend on the choice of $\lambda$ and we say $b\in\dcanB[\mu]$. We have $\dcanB=\coprod_{\lambda\in X_{+}}\dcanB[\lambda]$.

\begin{Notation} 
If $S\subset X_{+}$, then write $\dcanB[S]=\coprod_{\lambda\in S}\dcanB[\lambda]$ and 
$\dot{U}_{\intform}(\mathfrak{g})[S]:=\intform\dcanB[S]$.
\end{Notation}

\begin{Lemma}\label{L:basedideal}
Suppose that $S\subset X_{+}$ is such that:
\begin{enumerate}[itemsep=0.15cm,label=\emph{\upshape(\roman*)}]

\item If $\lambda\in S$ and $\lambda^{\prime}\in X_{+}$ is such that 
$\lambda^{\prime}>\lambda$, then $\lambda\in S$.

\item $X_{+}\setminus S$ is finite.

\end{enumerate}
Then $\dot{U}_{\intform}(\mathfrak{g})[S]\subset\dot{U}_{\intform}(\mathfrak{g})$ is a 
two-sided ideal. 
\end{Lemma}

\begin{proof}
This is \cite[Section 29.2.1]{Lu-quantumgroups-book}.
\end{proof}

Recall the notion of a homomorphism of based modules 
from \cite[Section 27.1.3]{Lu-quantumgroups-book}: such a homomorphism 
is a $U_{q}(\mathfrak{g})$-equivariant map $f\colon M\to M^{\prime}$ such that 
$B\subset B^{\prime}\cup\{0\}$.

\begin{Lemma}\label{L:U-mod-S-is-based}
Suppose that $S\subset X_{+}$ satisfies the hypothesis of \autoref{L:basedideal}. Then the 
quotient map 
\begin{gather*}
\dot{U}_{\intform}(\mathfrak{g})\to
\dot{U}_{\intform}(\mathfrak{g})/\dot{U}_{\intform}(\mathfrak{g})[S]
\end{gather*}
is a homomorphism of based modules.
\end{Lemma}

\begin{proof}
Since $\dot{U}_{\intform}(\mathfrak{g})[S]$ is a based ideal, the quotient, equipped with basis $\dcanB[X_{+}\setminus S]$, is a based module by \cite[Section 27.1.4]{Lu-quantumgroups-book}. 
\end{proof}

We recall the following key observation which is the main point in 
\cite[Corollary 1.5]{Ka-Based-Filt}:

\begin{Lemma}\label{L:based-is-weylffiltered}
Suppose that $M$ is a (finite dimensional) based 
$\dot{U}_{\intform}(\mathfrak{g})$-module. Then $M$ has a filtration by Weyl modules.
\end{Lemma}

\begin{proof}
This is a consequence of repeatedly applying \cite[Proposition 27.1.7]{Lu-quantumgroups-book}.
\end{proof}

\subsection{Quantum Schur algebras}\label{SS:QSchur}

We now recall the setting in \cite{Do-presenting-qSchur}, and explain 
how we can use this to quantize an important argument from \cite{AdRy-tilting-howe-positive-char}.

\begin{Definition}
Let $\pi\subset X_{+}$ be a set of dominant weights such that if $\lambda\in\pi$ 
and $\mu\leq\lambda$, then $\mu\in\pi$. We say 
that $\pi$ is a \emph{saturated} set of weights and write 
\begin{gather*}
W\acts\pi:=\{w\acts\lambda|\lambda\in\pi\}.
\end{gather*}
Write $\pi^{c}:=X_{+}\setminus\pi$.
\end{Definition}

Note that $\pi^{c}$ satisfies the hypothesis 
in \autoref{L:basedideal}.

\begin{Definition}\label{D:Schur}
We define the \emph{(generalized quantum) Schur algebra} 
$S_{\intform}^{\pi}(\mathfrak{g})$ as the quotient of 
$\dot{U}_{\intform}(\mathfrak{g})$ by the ideal generated by $\pr_{\chi}$ 
such that $\chi\notin W\changed{\acts}\pi$.
\end{Definition}

\autoref{D:Schur} appeared in \cite[Definition 1.5]{Do-presenting-qSchur}.
It is easy to see that the $\intform$-algebra $S_{\intform}^{\pi}(\mathfrak{g})$ 
is unital with unit $\pr_{\pi}:=\sum_{\lambda\in W\changed{\acts} \pi}\pr_{\lambda}$.

\begin{Proposition}\label{P:presentqschur}
There is an explicit isomorphism $\dot{U}_{\intform}(\mathfrak{g})/\dot{U}_{\intform}(\mathfrak{g})[\pi^{c}]\to S_{\intform}^{\pi}(\mathfrak{g})$ and $\dcanB[\pi]$ descends to an $\intform$-basis of $S_{\intform}^{\pi}(\mathfrak{g})$.
\end{Proposition}

\begin{proof}
See \cite[Theorem 4.2]{Do-presenting-qSchur}.
\end{proof}

\begin{Remark}
Let $\lambda\in\pi$. Then $\pr_{\pi}$ acts as the 
identity on $\Delta_{\intform}(\lambda)$ and
$\nabla_{\intform}(\lambda)$. Moreover, if $\ring$ 
is a field over $\intform$, then $\pr_{\pi}\acts T_{\ring}(\lambda)=T_{\ring}(\lambda)$.
\end{Remark}

\begin{Lemma}\label{L:Schurmod-is-highestweight}
Let $\ring$ be a field over $\intform$. The 
category $S_{\ring}^{\pi}(\mathfrak{g})\text{-mod}$ is a 
highest weight category, with poset $\pi$, and with 
standard modules, costandard modules, and tilting modules
given by $\Delta_{\ring}(\lambda)$, $\nabla_{\ring}(\lambda)$ and 
$T_{\ring}(\lambda)$ with indexing poset $\pi$.
\end{Lemma}

\begin{proof}
See \cite[Theorem 5.4]{Do-presenting-qSchur}
\end{proof}

The following is a standard consequence of \autoref{L:Schurmod-is-highestweight}:

\begin{Lemma}\label{C:dimhom-tilt}
If $T$ and $T^{\prime}$ are tilting 
$S_{\ring}^{\pi}(\mathfrak{g})$-modules with Weyl characters 
\begin{gather*}
[T]=\sum_{\lambda\in\pi}m_{\lambda}\cdot[\Delta_{\ring}(\lambda)] 
,\quad
[T^{\prime}]=\sum_{\lambda\in\pi}m^{\prime}_{\lambda}\cdot[\Delta_{\ring}(\lambda)],
\end{gather*}
in the Grothendieck ring, then we have
\begin{gather*}
\dim_{\ring}\Hom_{S_{\ring}^{\pi}(\mathfrak{g})}(T,T^{\prime}) 
=\sum_{\lambda\in\pi}m_{\lambda}m^{\prime}_{\lambda}.
\end{gather*}
\end{Lemma}

\begin{proof}
The dual Weyl character is the same as the Weyl character, so 
\begin{gather*}
[T^{\prime}]=\sum_{\lambda\in\pi}m^{\prime}_{\lambda}[\nabla_{\ring}(\lambda)].
\end{gather*}
\autoref{L:Schurmod-is-highestweight} implies 
the usual Ext-vanishing, cf. \autoref{L:Weylmodule-omnibuslemma}.(c),
and the result then follows from a standard long exact sequence argument.
\end{proof}

\begin{Lemma}\label{L:dimofschur}
The Schur algebra $S_{\intform}^{\pi}(\mathfrak{g})$ is free over $\intform$ of rank 
\begin{gather*}
\mathrm{rk}_{\intform}S_{\intform}^{\pi}(\mathfrak{g})=
\sum_{\lambda\in \pi}\left(\dim L_{q}(\lambda)\right)^{2}.
\end{gather*}
\end{Lemma}

\begin{proof}
By \cite[Section 29.2.1]{Lu-quantumgroups-book} we have
\begin{gather*}
\dim_{\field}\dot{U}_{q}(\mathfrak{g})/\dot{U}_{q}(\mathfrak{g})[\pi^{c}]=
\sum_{\lambda\in\pi}\big(\dim_{\field}L_{q}(\lambda)\big)^{2},
\end{gather*}
so $|\dcanB[\pi]|=\sum_{\lambda\in\pi}\left(\dim L_{q}(\lambda)\right)^{2}$. 
\autoref{P:presentqschur} implies that $S_{\intform}^{\pi}(\mathfrak{g})$ 
has a basis naturally in bijection with $\dcanB[\pi]$.
\end{proof}

\begin{Lemma}
If $\ring$ is a field over $\intform$, then
\begin{gather*}
\dim_{\ring}S_{\ring}^{\pi}(\mathfrak{g})= 
\sum_{\lambda\in\pi}\big(\dim\Delta_{\ring}(\lambda)\big)^{2}.
\end{gather*}
\end{Lemma}

\begin{proof}
Since Weyl modules are free over $\intform$, it follows that 
\begin{gather*}
\dim_{\ring}\Delta_{\ring}(\lambda)
=\dim_{\field}\Delta_{q}(\lambda)
=\dim_{\field}L_{q}(\lambda).
\end{gather*}
Therefore, the claim follows from \autoref{L:dimofschur}. 
\end{proof}

\begin{Lemma}\label{L:intschurhasfilt}
The left regular $S_{\intform}^{\pi}(\mathfrak{g})$-module $S_{\intform}^{\pi}(\mathfrak{g})$ has a filtration by standard modules. 
\end{Lemma}

\begin{proof}
\autoref{L:U-mod-S-is-based} and \autoref{L:based-is-weylffiltered} imply that $\dot{U}_{\intform}(\mathfrak{g})/\dot{U}_{\intform}(\mathfrak{g})[\pi^{c}]$ is Weyl filtered as a module over $\dot{U}_{\intform}(\mathfrak{g})$. The claim then follows from \autoref{P:presentqschur}.
\end{proof}

\begin{Lemma}\label{L:schurhasfilt}
The left regular $S_{\ring}^{\pi}(\mathfrak{g})$-module $S_{\ring}^{\pi}(\mathfrak{g})$ has a filtration by standard modules. 
\end{Lemma}

\begin{proof}
Directly from \autoref{L:intschurhasfilt}.
\end{proof}

\begin{Lemma}\label{L:standardembeds}
Let $X$ be an $S_{\ring}^{\pi}(\mathfrak{g})$-module.
If $X$ has a finite standard filtration, then $X$ embeds in a tilting module.
\end{Lemma}

\begin{proof}
Each standard module $\Delta$ embeds in an indecomposable tilting module $T_{\Delta}$. If 
\begin{gather*}
0\rightarrow K\xrightarrow{i}X\xrightarrow{p}\Delta\rightarrow 0
\end{gather*}
and $K$ has a standard filtration, then we can assume inductively that there are embeddings into tilting modules:
\begin{gather*}
0\rightarrow K\xrightarrow{a}T_{K} 
\quad\text{and}\quad 
0\rightarrow\Delta\xrightarrow{b}T_{\Delta}.
\end{gather*}
Since $\Ext^{1}(\Delta,T_{K})=0$, there is a surjection
\begin{gather*}
\Hom(X,T_{K})\xrightarrow{i^{*}}\Hom(K,T_{K})\rightarrow 0.
\end{gather*}
Let $a^{\prime}\colon X\rightarrow T_{K}$ be such that $a^{\prime}\circ i=a$ and define
\begin{gather*}
X\to T_{K}\oplus T_{\Delta},\quad x\mapsto\big(a^{\prime}(x),b\circ p(x)\big).
\end{gather*}
If $x\mapsto 0$, then injectivity of $b$ implies $x\in\ker p=i(K)$. Injectivity of $a$ 
then implies $x=i(k)=0$. Therefore, $X$ embeds in the tilting module $T_{K}\oplus T_{\Delta}$. 
\end{proof}

\begin{Lemma}\label{L:schurembeds}
The left regular $S_{\ring}^{\pi}(\mathfrak{g})$-module 
$S_{\ring}^{\pi}(\mathfrak{g})$ embeds in a tilting module. 
\end{Lemma}

\begin{proof}
Combing \autoref{L:schurhasfilt} with \autoref{L:standardembeds}.
\end{proof}

The following is crucial:

\begin{Proposition}\label{L:fulltiltisfaithful}
If $T\in S_{\ring}^{\pi}(\mathfrak{g})\modules$ is a full 
tilting module, then $T$ is a faithful module.
\end{Proposition}

\begin{proof}
Let $T$ be a full tilting module. The left regular module $S_{\ring}^{\pi}(\mathfrak{g})$ is faithful and embeds in a tilting module. Hence, it embeds in a full tilting module $T^{\prime}$ and it follows that $T^{\prime}$ is faithful. Since 
\begin{gather*}
T=\bigoplus_{\lambda\in\pi}T^{\oplus m_{\lambda}}(\lambda)
\quad\text{and}\quad 
T^{\prime}=\bigoplus_{\lambda\in\pi}T^{\oplus m^{\prime}_{\lambda}}(\lambda)
\end{gather*}
such that \changed{$m_{\lambda},m_{\lambda}^{\prime}\in\N$}, it follows that 
there is some $n\in\N$ such that $T^{\prime}$ is a summand of 
$(T)^{\oplus n}$. Since $T^{\prime}$ is faithful, 
$T^{\oplus n}$ is also faithful. A module $M$ is faithful 
if and only if $M^{\oplus n}$ is faithful for some $n$, 
so $T$ is faithful as well.  
\end{proof}

\subsection{Integral Howe duality}

Recall from \autoref{P:IntActionExt} that we have commuting actions 
of $U_{\intform}(\mathfrak{sp}_{2n})$ and $U_{\intform}(\mathfrak{sl}_{2})$ on 
$\Lambda_{\intform}$. We will use these actions for the remainder of this section.

\begin{Lemma}\label{L:Surjective}
An $\intform$-linear map between two finitely generated free $\intform$-modules is an isomorphism 
if it is an isomorphism for all specializations.
\end{Lemma}

\begin{proof}
Let $f\colon\intform^{m}\rightarrow\intform^{n}$ be an $\intform$ module map. Since $\field\otimes f\colon\field^{m}\rightarrow\field^{n}$ is an isomorphism, it follows that $m=n$. In particular, $f$ can be represented as a square matrix. 

Suppose now for the sake of contradiction that $\det f$ is not invertible in $\intform$. Therefore, the ideal $(\det f)\subset\intform$ is contained in some maximal ideal $M$ and 
$\intform/M\otimes\det f=0$. Since $f$ is an isomorphism for all specializations, 
it follows that $\intform/M\otimes f$ is an isomorphism. Therefore 
$\intform/M\otimes\det f=\det(\intform/M\otimes f)\neq 0$, 
a contradiction. 

Since $\det f$ is 
therefore invertible in $\intform$, it follows that $f$ is invertible. 
\end{proof}

\begin{Proposition}\label{P:IntHowe}
The homomorphisms induced from the actions
\begin{gather*}
\phi_{\intform}\colon
U_{\intform}(\mathfrak{sp}_{2n})\to\End_{{U}_{\intform}(\mathfrak{sl}_{2})}(\Lambda_{\intform})
,\quad
\psi_{\intform}\colon
{U}_{\intform}(\mathfrak{sl}_{2})\to\End_{U_{\intform}(\mathfrak{sp}_{2n})}(\Lambda_{\intform})
\end{gather*}
are surjective. 
\end{Proposition}

\begin{proof}
The set $\pi:=\{\varpi_{0},\varpi_{1},\dots,\varpi_{n}\}\subset X_{C_{n}}$ 
is a saturated set of dominant weights, and the set 
of weights occurring in $\Lambda_{\intform}$ is equal 
to $W_{C_{n}}\cdot\pi$. Therefore, the map $\phi_{\intform}$ 
induces an $\intform$-algebra homomorphism 
$\overline{\phi}_{\intform}\colon S_{\intform}^{\pi}(\mathfrak{sp}_{2n})
\to\End_{{U}_{\intform}(\mathfrak{sl}_{2})}(\Lambda_{\intform})$. 

Since $\Lambda_{\ring}=\oplus_{k=0}^{2n}\Lambda_{\ring}^{k}$ and 
$T_{\ring}(\varpi_{k})$ is a summand of 
$\Lambda_{\ring}(\varpi_{k})$ for $k\in\{0,1,\dots,n\}$, 
it follows that $\Lambda_{\ring}$ is a full tilting module 
for $S_{\ring}^{\pi}(\mathfrak{sp}_{2n})$. By 
\autoref{L:fulltiltisfaithful} $\overline{\phi}_{\ring}$ is injective. 

Using \autoref{C:dimhom-tilt} and \autoref{P:sl2tilting} we find 
\begin{gather*}
\dim_{\ring}\Hom_{{U}_{\ring}(\mathfrak{sl}_{2})}(\Lambda_{\ring},\Lambda_{\ring})
=\sum_{k=0}^{n}\big(\dim_{\ring}\Delta_{\ring}(\varpi_{k})\big)^{2}.
\end{gather*}
It then follows from \autoref{L:dimofschur} that 
\begin{gather*}
\dim S_{\ring}^{\pi}(\mathfrak{sp}_{2n})= 
\dim_{\ring}\End_{{U}_{\ring}(\mathfrak{sl}_{2})}(\Lambda_{\ring}),
\end{gather*}
so $\overline{\phi}_{\ring}$ is surjective. Hence, $\overline{\phi}_{\intform}$ is surjective 
by \autoref{L:Surjective}, so $\phi_{\intform}$ is surjective.

The set of $\mathfrak{sl}_{2}$ weights $\tau:=\{0,1,\dots,n\}\subset X_{A_{1}}$ 
is a saturated set of dominant weights, and the set of $\mathfrak{sl}_{2}$ weights 
occurring in $\Lambda_{\intform}$ is equal to $W_{A_{1}}\cdot\tau$. So there is an induced map $\overline{\psi}_{\intform}:S_{\intform}^{\tau}(\mathfrak{sl}_2)\rightarrow \End_{U_{\intform}(\mathfrak{sp}_{2n})}(\Lambda_{\intform})$. Also, $\Lambda_{\ring}$ is a full 
tilting module for $S_{\ring}^{\tau}(\mathfrak{sl}_{2})$, 
so $\overline{\psi}_{\ring}$ is injective. 
\autoref{P:sptilting} implies $\Lambda_{\ring}$ is a 
$U_{\ring}(\mathfrak{sp}_{2n})$ tilting module, with 
$\Delta$ character given by \autoref{P:spfilt}, so a 
similar argument to the one above shows that
\begin{gather*}
\dim_{\ring}S_{\ring}^{\tau}(\mathfrak{sl}_{2})= 
\sum_{k=0}^{n}\big(\dim_{\ring}\Delta_{\ring}(k)\big)^{2}= 
\dim_{\ring}\End_{U_{\ring}(\mathfrak{sp}_{2n})}(\Lambda_{\ring}),
\end{gather*}
and therefore that $\overline{\psi}_{\ring}$ is surjective. This implies that 
$\overline{\psi}_{\intform}$ is surjective, again by \autoref{L:Surjective}. Thus, $\psi_{\intform}$ is surjective.
\end{proof}

\begin{Proposition}\label{P:IntHoweTwo}
The $\textbf{U}_{\intform}$-module $\Lambda_{\intform}$ is a
Howe tilting module satisfying
\autoref{Eq:MainOne} and \autoref{Eq:MainTwo}.
\end{Proposition}

\begin{proof}
Note that \autoref{Eq:MainOne} implies \autoref{Eq:MainTwo}, and recall that
the tilting module properties of $\Lambda_{\intform}$ were shown in
\autoref{P:sl2tilting}, \autoref{P:sptilting} and \autoref{SS:QSchur}.
Thus, it remains to show \autoref{Eq:MainOne}. This is however a formal 
consequence of Ringel duality, see 
\cite[Appendix]{Do-q-schur}, and that Ringel duality can be used in our setting is 
\autoref{P:IntHowe}.
\end{proof}

\begin{proof}[Proof of \autoref{T:IntroMain}]
This theorem follows now from \autoref{P:IntActionExt}, \autoref{P:IntHowe} and \autoref{P:IntHoweTwo}. 
(Note that we formulated \autoref{P:IntActionExt} in terms of a left and a right action, hence 
the appearance of the opposite algebra when one compares 
that theorem to \autoref{P:IntHoweTwo}, see also \autoref{R:IntroMain}.(a).)
\end{proof}

\begin{Remark}
For the reader wondering whether, as $S_{\intform}^{\pi}(\mathfrak{sp}(V))$-modules,
we have $\pr_{\pi}\cdot V^{\otimes k}\cong\Lambda^{k}$: this is not the case. For example,
$V^{\otimes 3}$ has Weyl filtration containing three copies of $V$ and one copy of $\Delta(\varpi_{3})$, while $\Lambda^{3}$ has Weyl filtration $\Delta(\varpi_{3})$ and one copy of $V$.
\end{Remark}

\section{Further results I: The dual action via differential operators}\label{SS:DualAction}

In the final three sections, including this one, we collect results that are not relevant for \autoref{T:IntroMain} but of independent interest.

Our construction of the $\mathfrak{sl}_{2}$ in \autoref{SS:SLAction} is quite different 
from Howe's approach in \cite{Ho-perspectives-invariant-theory} who uses 
a form of harmonic analysis, that is, differential operators. 
We now present the quantum version of this approach, following ideas in \cite{Su-harmonic-analysis-son}.

\begin{Remark}
The results in this section, up to some 
details regarding idempotented forms, give the same results as 
\autoref{SS:SLAction}. We therefore decided to leave all proofs to the reader.
\end{Remark}

Using \autoref{P:SelfDual}, we will write $\Lambda_{\intform}^{\ast}$ for $\Lambda_{\intform}$ 
with the generating set being $\{d_{i}|i\in[1,-1]\}$.
In the following definition, which defines an $\intform$-algebra containing $\Lambda_{\intform}$ 
and $\Lambda_{\intform}^{\ast}$,
we write $\dv_{i}\leftrightsquigarrow v_{i}$ for the elements of $\Lambda_{\intform}$
and $\dd_{j}\leftrightsquigarrow d_{j}$ for 
the elements of $\Lambda_{\intform}^{\ast}$ (we decided to change notation 
to stress the differential operator calculus).

\begin{Definition}\label{D:DiffOp}
Let $\mathcal{D}_{\intform}$ be the $\intform$-algebra generated by elements 
\begin{gather*}
\dv_{1},\dots,\dv_{n},\dv_{-n},\dots,\dv_{-1},
\dd_{1},\dots,\dd_{n},\dd_{-n},\dots,\dd_{-1},
\end{gather*}
subject only to the following relations.
The relations 
of $\Lambda_{\intform}$ for the $\dv_{i}$, the relations of 
$\Lambda_{\intform}$ for the $\dd_{j}$ and:
\begin{gather*}
\dd_{i}\dv_{i^{\prime}}=
-q^{2}\cdot\dv_{i^{\prime}}\dd_{i},i\in[1,-1],
\\
\dd_{i}\dv_{i}=-\dv_{i}\dd_{i}+q(q-q^{-1})\sum_{k\in[1,i-1]}\dv_{k}\dd_{k}+1,i\in[1,n],
\\
\dd_{-i}\dv_{-i}=-\dv_{-i}\dd_{-i}
+q^{2(n-i+1)+1}(q-q^{-1})\cdot\dv_{i}\dd_{i}+1
+q(q-q^{-1})\sum_{k\in[1,-(i+1)]}\dd_{k}\dv_{k},i\in[1,n],
\end{gather*}
\changed{and for the final three relations we always have $i,j\in[1,-1],i\notin\{j,j^{\prime}\}$:}
\begin{gather*}
\dd_{i}\dv_{j}=-q\cdot\dv_{j}\dd_{i},i^{\prime}>j,
\\
\dd_{i}\dv_{j}=-q\cdot\dv_{j}\dd_{i}
+q(-q)^{|j|-|i^{\prime}|}(q-q^{-1})\cdot\dv_{i^{\prime}}\dd_{j^{\prime}}, 
i^{\prime}<j,\text{sgn}(i)\neq\text{sgn}(j),
\\
\dd_{i}\dv_{j}=-q\cdot\dv_{j}\dd_{i}+
q^{2}(-q)^{|j|-|i^{\prime}|}(q-q^{-1})\cdot\dv_{i^{\prime}}\dd_{j^{\prime}}, 
i^{\prime}<j,\text{sgn}(i)=\text{sgn}(j).
\end{gather*}
Here, for $i\in[1,-1]$, let $i^{\prime}=2n+1-i$.
Moreover, $\text{sgn}(i)=1$ and $|i|=i$ for $i\in[1,n]$, 
and $\text{sgn}(i)=-1$ and $|i|=2n+i+1$ for $i\in[-n,-1]$.
\end{Definition}

The relation $\dd_{i}\dv_{i}=-\dv_{i}\dd_{i}+q(q-q^{-1})\sum_{k\in[1,i-1]}\dv_{k}\dd_{k}+1$
and its variations are sometimes called \emph{$q$-Heisenberg commutation relations}. The name comes 
from the classical version of this relation being $\dd_{i}\dv_{i}=-\dv_{i}\dd_{i}+1$.

\begin{Remark}
The $\intform$-algebra $\mathcal{D}_{\intform}$ from \autoref{D:DiffOp} 
is the symplectic analog of \cite[(2.26)]{Su-harmonic-analysis-son}, 
i.e. the \emph{differential operator algebra} on $\Lambda_{\intform}$.
The $\intform$-algebra $\mathcal{D}_{\intform}$ has similar properties 
as the orthogonal analog in \cite{Su-harmonic-analysis-son}. 
For example, using the same notation as for $\Lambda_{\intform}$, 
the first result one can show is that
\begin{gather*}
\{\dd_{T}\dv_{S}|S,T\subset[1,-1]\}
\end{gather*}
is an $\intform$-basis of $\mathcal{D}_{\intform}$ so that 
$\mathcal{D}_{\intform}\cong\Lambda_{\intform}
\otimes_{\intform}\Lambda_{\intform}^{\ast}$ as free $\intform$-modules.
For brevity, we decided to not include any further discussion of $\mathcal{D}_{\intform}$ 
and its properties.
\end{Remark}

Similarly as in \cite{Su-harmonic-analysis-son} one can show that 
$\mathcal{D}_{\intform}$ gives rise to an analog of the oscillator 
$U_{q}(\mathfrak{sl}_{2})$-action on $\Lambda_{\intform}$. 
The rest of this section is an algebraic description of this action.

To this end, we will work with the idempotented quantum group for $\mathfrak{sl}_{2}$ which need the 
following weight operators:

\begin{Definition}\label{D:WeightOps}
For $i\in[1,n]$ define the \emph{($\mathfrak{sl}_{2}$) weight operators} $\pr_{i}$ on $\Lambda_{q}$ by
$\pr_{i}\acts v_{S}:=\delta_{i,-n+|S|}\cdot v_{S}$.
\end{Definition}

Note that $\pr_{i}\pr_{j}=\delta_{ij}\pr_{i}$.

\begin{Remark}
The nondegenerate alternating two form from \autoref{SS:HowePartOne}, determines an nonzero vector $\omega\in\Lambda_{1}^{2}$. The operator $E$ is then multiplication by $\omega$, while $F$ is contraction with $\omega$. 

The elements
\begin{gather*}
\omega_{q}:=\sum_{i=1}^{n}(-q)^{i}\dv_{i}\dv_{-i}
,\quad
\omega_{q}^{\vee}:=\sum_{i=1}^{n}(-q)^{-i}\dd_{i}\dd_{-i}.
\end{gather*}
both span trivial $U_{\intform}(\mathfrak{sp}_{2n})$-submodules of $\mathcal{D}_{\intform}$. If we view $\omega_{q}$ as an element in $\Lambda_{\intform}$, then this is the $q$-analog of $\omega$. There is a natural action of $\mathcal{D}_{\intform}$ on $\Lambda_{\intform}$, inducing operators in $\End_{U_{\intform}(\mathfrak{sp}_{2n})}(\Lambda_{\intform})$: multiplication by $\omega_{q}$ and contraction with $\omega_{q}^{\vee}$. The following definition records the action of these operators explicitly.
\end{Remark}

\begin{Definition}\label{D:EFOps}
For $S\subset[1,-1]$ and $i\in[1,n]$ we write 
$S_{i\rightarrow}:=S\cap[i+1,-(i+1)]$ and $S^{c}=[1,-1]\setminus S$.
Define operators $e$ and $f$ on $\Lambda_{q}$ by
\begin{gather*}
e\acts v_{S}:=\sum_{\{i,-i\}\subset S^{c}}(-q)^{i}(-q)^{|S_{i\rightarrow}|}\cdot v_{S\cup\{i,-i\}}
,\quad
f\acts v_{S}:=\sum_{\{i,-i\}\subset S}(-q)^{i}(-q)^{|S_{i\rightarrow}|}\cdot v_{S\setminus\{i,-i\}}
\end{gather*}
We call $e$ and $f$ the \emph{($\mathfrak{sl}_{2}$) Chevalley operators}.
\end{Definition}

Let us explain the action of the Chevalley generators on dot diagrams in words:
\begin{enumerate}[itemsep=0.15cm,label=\emph{\upshape(\roman*)}]

\item The only nonzero action of $f$ is on fully dotted columns, and $f$ acts summing over all such columns.

\item $f$ removes the two dots in fully dotted columns.

\item The $q$-factor that one gets is obtained from the column number $i$, giving $(-q)^{i}$, and the number of dots to the right, contributing $(-q)^{|S_{i\rightarrow}|}$. Here the set $S_{i\rightarrow}$ counts nodes to the right, for example, as follows:
\begin{gather*}
\scalebox{0.95}{$S_{1\rightarrow}=\{2,3,-3,-2\}\leftrightsquigarrow
\begin{ytableau}
i & *(blue!50) & *(blue!50) \\
\phantom{a} & *(blue!50) & *(blue!50)
\end{ytableau}
\,,\,
S_{2\rightarrow}=\{3,-3\}\leftrightsquigarrow
\begin{ytableau}
\phantom{a} & i & *(blue!50) \\
\phantom{a} & \phantom{a} & *(blue!50)
\end{ytableau}
\,,\,
S_{3\rightarrow}=\emptyset\leftrightsquigarrow
\begin{ytableau}
\phantom{a} & \phantom{a} & i \\
\phantom{a} & \phantom{a} & \phantom{a}
\end{ytableau}
\,,$}
\end{gather*}
where the shaded nodes are in $S_{i\rightarrow}$.

\item The action of $e$ is similar, but $e$ acts on completely empty columns and adds dots instead of 
removing them.

\end{enumerate}
Here is an explicit example:

\begin{Example}\label{E:ActionF}
The action of the operator $f$ on dot diagrams is exemplified by:
\begin{gather*}
f\acts
\begin{ytableau}
*(spinach!50)\bullet & *(spinach!50)\bullet & \phantom{a} \\
*(spinach!50)\bullet & *(spinach!50)\bullet & \phantom{a}
\end{ytableau}
=
(-q)^{1}(-q)^{2}\cdot
\begin{ytableau}
\phantom{a} & *(spinach!50)\bullet & \phantom{a} \\
\phantom{a} & *(spinach!50)\bullet & \phantom{a}
\end{ytableau}
+(-q)^{2}(-q)^{0}\cdot
\begin{ytableau}
*(spinach!50)\bullet & \phantom{a} & \phantom{a} \\
*(spinach!50)\bullet & \phantom{a} & \phantom{a}
\end{ytableau}
\,.
\end{gather*}
Note that this is different from the action of the \changed{Chevalley generators $f_{i}$ of the symplectic side} which fix the number of dots.
\end{Example}

\begin{Lemma}\label{L:Commute}
The operators $e$ and $f$ commute with the operators $e_{i}$ and $f_{i}$ for $\{1,\dots,n\}$.
\end{Lemma}

\begin{proof}
Omitted.
\end{proof}

Following \cite{BeLuMaPh-quantum-group},
the \emph{idempotented quantum group} $\dot{U}_{q}(\mathfrak{sl}_{2})$ 
is obtained by the $\field$-algebra extension of $U_{q}(\mathfrak{sl}_{2})$ by orthogonal 
idempotents $\pr_{k}$ for $k\in\Z$, the 
\emph{(dot $\mathfrak{sl}_{2}$) weight generators}, such that, 
\begin{gather}\label{Eq:DotSL2}
\pr_{k+2}E=E\pr_{k}
,\quad
\pr_{k-2}F=F\pr_{k}
,\quad
\pr_{k}K=K\pr_{k}=q^{k}\pr_{k},
\end{gather}
and then idempotent truncation, namely
$\dot{U}_{q}(\mathfrak{sl}_{2})=\bigoplus_{k,j\in\Z}\pr_{k}U_{q}(\mathfrak{sl}_{2})\pr_{j}$.

In $\dot{U}_{q}(\mathfrak{sl}_{2})$ we can use \autoref{Eq:DotSL2} to 
write elements of $\dot{U}_{q}(\mathfrak{sl}_{2})$ as strings 
of $E$ and $F$ with one weight generator on the right, and we will do this in 
the following. Note that in $\dot{U}_{q}(\mathfrak{sl}_{2})$ we have that 
$\pr_{k+2}E=E\pr_{k}$ replaces $KE=q^{2}EK$,
$\pr_{k-2}F=F\pr_{k}$ replaces $KF=q^{-2}FK$ and 
and $EF-FE=\frac{K-K^{-1}}{q-q^{-1}}$ gets replaced by 
\begin{gather*}
EF\pr_{k}-FE\pr_{k}=[k]_{q}\pr_{k}.
\end{gather*}

We almost have an action of $\dot{U}_{q}(\mathfrak{sl}_{2})$ on $\Lambda_{q}$: 

\begin{Example}
Continuing \autoref{E:ActionF} by applying $e$ gives:
\begin{gather*}
\scalebox{0.95}{$e\acts\bigg(f\acts
\begin{ytableau}
*(spinach!50)\bullet & *(spinach!50)\bullet & \phantom{a} \\
*(spinach!50)\bullet & *(spinach!50)\bullet & \phantom{a}
\end{ytableau}\,\bigg)
=
(-q)^{6}\cdot
\begin{ytableau}
\phantom{a} & *(spinach!50)\bullet & *(spinach!50)\bullet \\
\phantom{a} & *(spinach!50)\bullet & *(spinach!50)\bullet
\end{ytableau}
+(-q)^{5}\cdot
\begin{ytableau}
*(spinach!50)\bullet & \phantom{a} & *(spinach!50)\bullet \\
*(spinach!50)\bullet & \phantom{a} & *(spinach!50)\bullet
\end{ytableau}
+\big((-q)^{6}+(-q)^{4}\big)\cdot
\begin{ytableau}
*(spinach!50)\bullet & *(spinach!50)\bullet & \phantom{a} \\
*(spinach!50)\bullet & *(spinach!50)\bullet & \phantom{a}
\end{ytableau}
\,.$}
\end{gather*}
On the other hand, we also have
\begin{gather*}
f\acts\bigg(e\acts
\begin{ytableau}
*(spinach!50)\bullet & *(spinach!50)\bullet & \phantom{a} \\
*(spinach!50)\bullet & *(spinach!50)\bullet & \phantom{a}
\end{ytableau}\,\bigg)
=(-q)^{3}\cdot f\acts
\begin{ytableau}
*(spinach!50)\bullet & *(spinach!50)\bullet & *(spinach!50)\bullet \\
*(spinach!50)\bullet & *(spinach!50)\bullet & *(spinach!50)\bullet
\end{ytableau}
\\
=
(-q)^{8}\cdot
\begin{ytableau}
\phantom{a} & *(spinach!50)\bullet & *(spinach!50)\bullet \\
\phantom{a} & *(spinach!50)\bullet & *(spinach!50)\bullet
\end{ytableau}
+(-q)^{7}\cdot
\begin{ytableau}
*(spinach!50)\bullet & \phantom{a} & *(spinach!50)\bullet \\
*(spinach!50)\bullet & \phantom{a} & *(spinach!50)\bullet
\end{ytableau}
+(-q)^{6}\cdot
\begin{ytableau}
*(spinach!50)\bullet & *(spinach!50)\bullet & \phantom{a} \\
*(spinach!50)\bullet & *(spinach!50)\bullet & \phantom{a}
\end{ytableau}
\,.
\end{gather*}
In other words, the action of $(ef-fe)$ is a bit off compared to 
what we would want for quantum $\mathfrak{sl}_{2}$ \changed{(namely a scalar multiple of the dot diagram one starts with)}, which is the reason for the rescaling 
in \autoref{D:EFTwo} below.
\end{Example}

\begin{Definition}\label{D:EFTwo}
Define the \emph{dot ($\mathfrak{sl}_{2}$) Chevalley operators} $E$ and $F$
\begin{gather*}
E\pr_{-n+|S|}:=q^{-|S|-1}\cdot e\pr_{-n+|S|}
,\quad
F\pr_{-n+|S|}:= q^{-n}\cdot f\pr_{-n+|S|},
\end{gather*}
with $e$ and $f$ being as in \autoref{D:EFOps}, and with $\pr_{i}$ as in \autoref{D:WeightOps}.
\end{Definition}

Next, the \changed{``correct''} action:

\begin{Proposition}\label{P:SL2Action}
There is an $\dot{U}_{q}(\mathfrak{sl}_{2})$-action on $\Lambda_{q}$ by assigning 
the generators of $\dot{U}_{q}(\mathfrak{sl}_{2})$ to the operators with the same name in 
\autoref{D:EFTwo} and \autoref{D:WeightOps}.

In particular, we have
\begin{gather*}
(EF-FE)\pr_{-n+|S|}\acts v_{S}=[-n +|S|]_{q}\acts v_{S}.
\end{gather*}
\end{Proposition}

\begin{proof}
Omitted.
\end{proof}

\begin{Proposition}\label{P:QBimodule}
The two actions, the $U_{q}(\mathfrak{sp}_{2n})$-action and the $\dot{U}_{q}(\mathfrak{sl}_{2})$-action, on $\Lambda_{q}$ commute.
\end{Proposition}

\begin{proof}
Omitted.
\end{proof}

Let $\dot{U}_{q}=U_{q}(\mathfrak{sp}_{2n})\otimes\dot{U}_{q}(\mathfrak{sl}_{2})$
and note that \autoref{P:QBimodule} gives a $\dot{U}_{q}$-action on $\Lambda_{q}$.

\begin{Proposition}\label{P:ActionsAgree}
The $\dot{U}_{q}$-action and the $U_{q}$-action on 
$\Lambda_{q}$, the one in this section and the one from 
\autoref{SS:SLAction}, give the same image.
\end{Proposition}

\begin{proof}
Omitted.
\end{proof}

\begin{Remark}
After defining the appropriate integral form of the action of $\dot{U}_{q}$,
there is an analogous statement for the integral actions. Then
\autoref{P:ActionsAgree} says that, for the purpose of this paper, 
both actions do the same job, and it depends on the problem at hand which 
action is the more suitable one.
\end{Remark}

As explained in \cite{Su-harmonic-analysis-son} 
for the orthogonal case, pushing this further leads to the symplectic version of harmonic analysis on quantum spheres.

\section{Further results II: On fundamental tilting modules}\label{SS:Tilting}

For $U_{\intform}(\mathfrak{sl}_{2})$ it is well-known when the Weyl modules 
are simple, see for example \cite[Section 3.4]{Do-q-schur} or \cite[Proposition 3.3]{SuTuWeZh-mixed-tilting}. We now partially address the case of $U_{\intform}(\mathfrak{sp}_{2n})$.

Recall from \autoref{P:spfilt} and \autoref{P:sptilting} that
\begin{gather*}
\mathrm{rk}_{\intform}\Hom_{U_{\intform}(\mathfrak{sp}_{2n})}
\big(\Delta_{\intform}(\varpi_{k}),\Lambda_{\intform}^{k}\big) =1=
\mathrm{rk}_{\intform}\Hom_{U_{\intform}(\mathfrak{sp}_{2n})}
\big(\Lambda_{\intform}^{k},\nabla_{\intform}(\varpi_{\intform})\big).
\end{gather*}
We fix bases for these homomorphism spaces as follows: 
Fix highest weight vectors $v_{\varpi_{k}}^{+}\in\Delta_{\intform}(\varpi_{k})$ 
and $v_{\varpi_{k}}^{-}\in\nabla_{\intform}(\varpi_{k})$.
Let $\iota_{k}\colon\Delta_{\intform}(\varpi_{k})\to\Lambda_{\intform}^{k}$ be the unique homomorphism such that $v_{\varpi_{k}}^{+}\mapsto v_{\{1,\dots,k\}}$ and let 
$\pi_{k}\colon\Lambda_{\intform}^{k}\to\nabla_{\intform}(\varpi_{k})$ 
be the unique homomorphism such that $v_{\{1,\dots,k\}}\mapsto v_{\varpi_{k}}^{-}$. 
Note that the usual Ext-vanishing \autoref{L:Weylmodule-omnibuslemma}.(c) gives
\begin{gather*}
\mathrm{rk}_{\intform}\Hom_{U_{\intform}(\mathfrak{sp}_{2n})}
\big(\Delta_{\intform}(\varpi_{k}),\nabla_{\intform}(\varpi_{k})\big)=1
\end{gather*}
\changed{and a basis for this homomorphism space is 
given by $hs_{k}:=\pi_{k}\circ\iota_{k}$.}

\begin{Remark}
The construction above is standard, see e.g. \cite[Section 3]{AnStTu-cellular-tilting} 
for essentially the same type of basis.
\end{Remark}

\begin{Lemma}\label{L:PrePairing}
We have
\begin{gather*}
\Hom_{U_{\intform}(\mathfrak{sp}_{2n})}(\Lambda_{\intform}^{k},\Lambda_{\intform}^{k+2i})\circ\iota_{k} 
=\intform\{(E^{(i)}\acts\placeholder)\circ\iota_{k}\},
\\
\pi_{k}\circ\Hom_{U_{\intform}(\mathfrak{sp}_{2n})}(\Lambda_{\intform}^{k+2i},\Lambda_{\intform}^{k}) 
=\intform\{\pi_{k}\circ(F^{(i)}\acts\placeholder)\}.
\end{gather*}
\end{Lemma}

\begin{proof}
By \autoref{P:spfilt}, \autoref{P:sptilting} and 
\autoref{L:Weylmodule-omnibuslemma}.(c) we know the $\intform$-rank on both sides match, 
hence, for $F^{(i)}$ is this \autoref{L:Fi-mapspartialhwtohw}. The first claim follows similarly.
\end{proof}

\begin{Definition}
Consider the pairing 
\begin{gather*}
\beta_{k}^{k+2i}\colon
\Hom_{U_{\intform}(\mathfrak{sp}_{2n})}(\Lambda_{q}^{k+2i},\Lambda^{k})
\times
\Hom_{U_{\intform}(\mathfrak{sp}_{2n})}(\Lambda_{q}^{k},\Lambda_{q}^{k+2i})
\to
\Hom_{U_{\intform}(\mathfrak{sp}_{2n})}
\big(\Delta_{\intform}(\varpi_{k}),\nabla_{\intform}(\varpi_{k})\big)
\end{gather*}
defined by $(x,y)\mapsto\pi_{k}\circ(x\circ y)\circ\iota_{k}$. 
\end{Definition}

Recall that $\ring$ denotes a specialization of $\intform$.

\begin{Proposition}\label{P:BinomPairing}
The rank of the pairing $(\beta_{k}^{k+2i})_{\ring}$ is equal to the 
multiplicity of $T_{\ring}(\varpi_{k})$ as a summand of $\Lambda_{\ring}^{k+2i}$. 
Moreover,
\begin{gather*}
\mathrm{rk}_{\ring}(\beta_{k}^{k+2i})_{\ring}
=
\begin{cases}
1 & \text{if }\qbinn{n-k}{i}_{\ring}\neq 0,
\\
0 & \text{else}.
\end{cases}
\end{gather*}
\end{Proposition}

\begin{proof}
We have the well-known relation
\begin{gather*}
F^{(a)}E^{(b)}\pr_{c}
=
\sum_{d=0}^{\min\{a,b\}}
\qbinn{a-b-c}{d}_{\intform}
E^{(b-d)}F^{(a-d)}\pr_{c}
\end{gather*}
that holds in the idempotented form $\dot{U}_{q}(\mathfrak{sl}_{2})$ 
which we will use for $a=b=i$ and $j=n-k$. That is, by \autoref{L:PrePairing} 
and the weight space decomposition for the $U_{q}(\mathfrak{sl}_{2})$-action on $\Lambda_{\intform}$,
we need to know the value of $F^{(i)}E^{(i)}\pr_{n-k}$. This value is $\qbinn{n-k}{i}_{\intform}$ 
by the above formula as only the summands for $d=\min\{a,b\}$ is nonzero.
This shows that $\mathrm{rk}_{\ring}(\beta_{k}^{k+2i})_{\ring}$ is as claimed. 
The first part of the statement follows from general theory similarly to 
\cite[Section 4C]{AnStTu-cellular-tilting}.
\end{proof}

Let $p=\mathrm{char}(\ring)$ 
(with, as usual, $\mathrm{char}(\ring)=0$ interpreted as $p=\infty$) 
and let $\ell\in\Z_{\geq 1}$ be minimal such that $[\ell]_{\ring}=0$ with $\ell=\infty$ 
if no such minimal number exists.
Write $p^{(i)}=p^{i-1}\ell$ with $p^{(0)}=1$. The \emph{$(p,\ell)$-adic expansion} of 
$m\in\N$ is the tuple $[\dots,a_{k},\dots,a_{0}]_{p,\ell}$ with $a_{i}\in\{0,\dots,p-1\}$ 
for $i\neq 0$ and $a_{0}\in\{0,\dots,\ell-1\}$ such that $m=\sum_{k\in\N}a_{k}p^{(k)}$ 
(we often omit the infinitely many zeros for $k\gg 0$).
We compare $(p,\ell)$-adic expansions componentwise.

As a consequence, we get:

\begin{Proposition}\label{P:isasummand}
Let $n-k=[\dots,a_{0}]_{p,\ell}$ and $i=[\dots,b_{0}]_{p,\ell}$ be the 
$(p,\ell)$-adic expansions.

\begin{enumerate}

\item $T_{\ring}(\varpi_{k})$ appears as direct summands of $\Lambda_{\ring}^{k+2i}$ 
if and only if $[\dots,a_{0}]_{p,\ell}\geq[\dots,b_{0}]_{p,\ell}$.

\item If $\ell\geq n+1$ 
($n+1$ is the dual Coxeter number of $\mathfrak{sp}_{2n}$), then $T_{\ring}(\varpi_{k})$ 
is always a direct summand of $\Lambda_{\ring}^{k+2i}$.

\end{enumerate}
\end{Proposition}

\begin{proof}
\textit{(a).} By \autoref{P:BinomPairing} and the famous quantum Lucas' theorem. For the latter see 
e.g. \cite[(1.2.4)]{Ol-gen-powers} for an early account.

\textit{(b).} In this case all 
involved $(p,\ell)$-adic expansion 
have only the zeroth digit, and this then follows from (a).
\end{proof}

\begin{Example}\label{E:CalcsWithTilting}
\changed{If $\ring=\mathbb{F}_{7}$ and $\xi=2\in\ring$, then $(p=7,\ell=3)$. Let us take $74=[3,3,2]_{7,3}$, 
$28=[1,2,1]_{7,3}$ and $25=[1,1,1]_{7,3}$. We get $\qbinn{74}{28}_{\ring}=0$ and 
$\qbinn{74}{25}_{\ring}\neq 0$. Thus, if $n=78$ and $k=5$ so that $78-5+1=74$, then $T_{\ring}(\varpi_{5})$ 
is a direct summand of $\Lambda_{\ring}^{55=5+2\cdot 25}$ but not of $\Lambda_{\ring}^{61=5+2\cdot 28}$, which we can directly read of from $[3,3,2]_{7,3}\not\geq[1,2,1]_{7,3}$ but 
$[3,3,2]_{7,3}\geq[1,1,1]_{7,3}$.}
\end{Example}

\begin{Theorem}\label{T:tiltingcharacter}
The following equations inductively determine the Weyl characters of each fundamental tilting 
$U_{\ring}(\mathfrak{sp}_{2n})$-module. First, the initial conditions are
$[T_{\ring}(\varpi_{0})]=[\Delta_{\ring}(\varpi_{0})]$ and 
$[T_{\ring}(\varpi_{1})]=[\Delta_{\ring}(\varpi_{1})]$.
Then, for $k>1$:
\begin{gather*}
[T_{\ring}(\varpi_{k})] 
=\sum_{i\geq 0,k-2i\geq 0}[\Delta_{\ring}(\varpi_{k-2i})]
-\sum_{\substack{i\geq 1,k-2i\geq 0\\ 
\qbinn{n-k}{i}_{\ring}\neq 0}}
[T_{\ring}(\varpi_{k-2i})].
\end{gather*}
\end{Theorem}

\begin{proof}
For all $k\in\{0,\dots,n\}$, \autoref{P:isasummand} implies that
\begin{gather*}
\Lambda_{\ring}^{k}\cong
\bigoplus_{\substack{i\geq 0,k-2i\geq 0
\\ 
\qbinn{n-k}{i}_{\ring}\neq 0}}
T_{\ring}(\varpi_{k-2i}) 
= 
T_{\ring}(\varpi_{k})\oplus\bigoplus_{\substack{i\geq 1,k-2i \geq 0
\\ \qbinn{n-k}{i}_{\ring}\neq 0}}T_{\ring}(\varpi_{k-2i}).
\end{gather*}
In particular, 
\begin{gather*}
\Delta(\varpi_{0})\cong\Lambda_{\ring}^{0}\cong T_{\ring}(\varpi_{0}),\quad
\Delta(\varpi_{1})\cong\Lambda_{\ring}^{1}\cong T_{\ring}(\varpi_{1}).
\end{gather*}
The desired result follows \ochanged{by} looking at the equations determined by these isomorphisms in the Grothendieck ring, noting that the Weyl character of $\Lambda_{\ring}^{k}$ is determined by \ochanged{the results in} \autoref{P:spfilt}.
\end{proof}

\begin{Example}
In the setting of \autoref{E:CalcsWithTilting}, let $n=78$. The matrix
\begin{center}
\scalebox{0.9}{\begin{tabular}{C||C|C|C|C|C|C|C|C|C|C}
& T_{\ring}(\varpi_{0}) & T_{\ring}(\varpi_{1}) & T_{\ring}(\varpi_{2}) & T_{\ring}(\varpi_{3}) & T_{\ring}(\varpi_{4}) & T_{\ring}(\varpi_{5}) & T_{\ring}(\varpi_{6}) & T_{\ring}(\varpi_{7}) & T_{\ring}(\varpi_{8}) & T_{\ring}(\varpi_{9}) \\
\hline\hline
\Delta_{\ring}(\varpi_{0})	& 1 &  &  &  &  &  & &  &  &  \\
\hline
\Delta_{\ring}(\varpi_{1})	&  & 1 &  & 1 &  &  &  &  &  &  \\
\hline
\Delta_{\ring}(\varpi_{2})	&  &  & 1 &  &  &  & 1 &  &  &  \\
\hline
\Delta_{\ring}(\varpi_{3})	&  &  &  & 1 &  &  & &  &  &  \\
\hline
\Delta_{\ring}(\varpi_{4})	&  &  &  &  & 1 &  & 1 &  &  &  \\
\hline
\Delta_{\ring}(\varpi_{5})	&  &  &  &  &  & 1 &  &  &  & 1 \\
\hline
\Delta_{\ring}(\varpi_{6})	&  &  &  &  &  &  & 1 &  &  &  \\
\hline
\Delta_{\ring}(\varpi_{7})	&  &  &  &  &  &  &  & 1 &  & 1 \\
\hline
\Delta_{\ring}(\varpi_{8})	&  &  &  &  &  &  &  &  & 1 &  \\
\hline
\Delta_{\ring}(\varpi_{9})	&  &  &  &  &  &  &  &  &  & 1 \\
\end{tabular}}
\end{center}
is the matrix of tilting:Weyl character for $k\in\{0,\dots,9\}$. 
The whole matrix (left $(p=7,\ell=3)$ and on the right $(3,2)$ for comparison) is
\begin{gather*}
\begin{tikzpicture}[anchorbase,scale=1]
\node at (0,0) {\includegraphics[height=6cm]{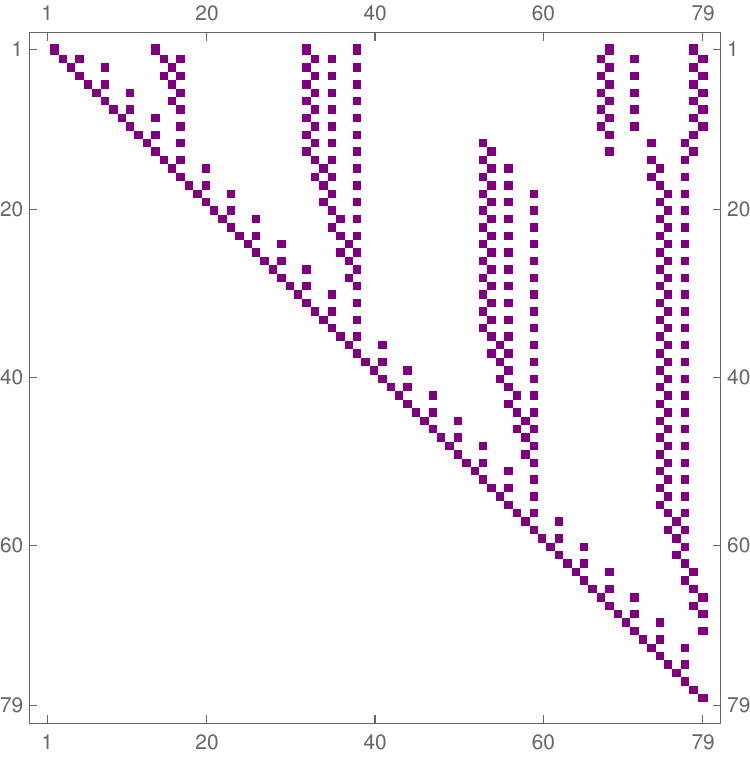}};
\node at (0,-2) {$(p=7,\ell=3)$};
\end{tikzpicture}
,\quad
\begin{tikzpicture}[anchorbase,scale=1]
\node at (0,0) {\includegraphics[height=6cm]{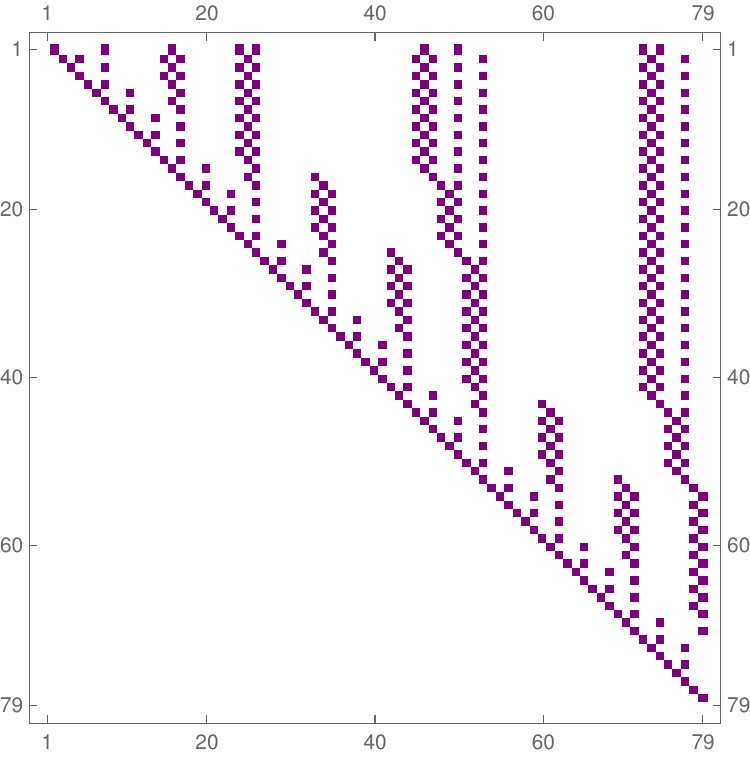}};
\node at (0,-2) {$(p=3,\ell=2)$};
\end{tikzpicture}
\end{gather*}
where entries $1$ are shaded boxes. Not displayed entries are zero in both cases.
\end{Example}

\begin{Remark}
Since $[T_{\ring}(\varpi_{k})]=[\Delta_{\ring}(\varpi_{k})]$ 
if and only if $\Delta_{\ring}(\varpi_{k})\cong L_{\ring}(\varpi_{k})$ 
by \autoref{L:AllFour}, \autoref{T:tiltingcharacter} also completely determines whether or not a Weyl module with fundamental highest weight is simple. When $q=1$ and the characteristic of $\ring$ is positive, this question has 
been answered in \cite{PrSu-fund-weyl-sp}. It is however not immediate
how to recover their result from \autoref{T:tiltingcharacter}, or how 
to deduce the $q$ analog of their result.
\end{Remark}

\begin{Proposition}
If $\ell\geq n+1$, then each fundamental Weyl $U_{\ring}(\mathfrak{sp}_{2n})$-module is simple. 
\end{Proposition}

\begin{proof}
By \autoref{T:tiltingcharacter} using \autoref{L:AllFour}.
\end{proof}

The following describes tilting:Weyl multiplicities on the symplectic 
side in terms of dual Weyl:simple on the side of $\mathfrak{sl}_{2}$, and also vice versa.

\begin{Proposition}
We have $\big(T(\varpi_{k}):\Delta(\varpi_{l})\big)=[\nabla(n-l):L(n-k)]$
and $\big(T(n-k):\Delta(n-l)\big)=[\nabla(\varpi_{l}):L(\varpi_{k})]$.
\end{Proposition}

The point is that the numbers $[\nabla(n-l):L(n-k)]$ and $\big(T(n-k):\Delta(n-l)\big)$
are well-known, see e.g.
\cite[Proposition 3.3]{SuTuWeZh-mixed-tilting}, so that we can get also 
the tilting:Weyl multiplicities.

\begin{proof}
The classical version of this is explained in \cite[Proposition 6.1.3]{McNi-charp-howe} 
(which is a consequence of abstract Ringel duality similarly as \autoref{Eq:MainOne}). 
Using \autoref{T:IntroMain}, the arguments 
given therein can be copied without problems.
\end{proof}

\begin{Example}\label{E:CalcsWithTiltings}
Continuing \autoref{E:CalcsWithTilting}, let $n=78$ and $k=5$. Then
$78-5+1=[3,3,2]_{7,3}$ and $\big(T(n-k):\Delta(n-l)\big)=0$ 
or $\big(T(n-k):\Delta(n-l)\big)=1$. The latter happens if and only if 
$n-l+1\in\{[3,3,2]_{7,3},[3,3,-2]_{7,3},[3,-3,2]_{7,3},[3,-3,-2]_{7,3}\}=\{74,70,56,52\}$.
This follows for example from
\cite[Proposition 3.3]{SuTuWeZh-mixed-tilting}.
Thus, we get $[\nabla(\varpi_{l}):L(\varpi_{k})]=1$ for $l\in\{5,9,23,27\}$.
\end{Example}

Lastly, we record a question 
we arrived at computing examples using \autoref{T:tiltingcharacter}.
In this question $q\!\dim_{\ring}$ denotes the \emph{quantum dimension}.

\begin{Question}
Let $n=[\dots,a_{0}]_{p,\ell}$ and $i=[\dots,b_{0}]_{p,\ell}$. Suppose that $\ell>2$. Do we have 
\begin{gather*}
\big(q\!\dim_{\ring}T_{\ring}(\varpi_{i})\neq 0\big) 
\Leftrightarrow
\big([\dots,a_{0}]_{p,\ell}\geq[\dots,b_{0}]_{p,\ell}\big)?
\end{gather*}
\end{Question}

\section{Further results III: Canonical bases}\label{SS:Canonical}

In this final section we justify that we called 
the $b_{S}$ from \autoref{D:Canonical}
canonical vectors.

\begin{Lemma}\label{L:generation}
The vectors $\{v_{k,k-2i}|0\leq k\leq 2n,k-2i\in\N\}$, described in \autoref{N:almosthwvectors}, generate $\Lambda_{\intform}$ as a $U_{\intform}(\mathfrak{sp}_{2n})$-module. Similarly for fixed $k$.
\end{Lemma}

\begin{proof}
We prove the claim for each fixed $k$, form which the general claim follows. The case $k=0,1$ is immediate. Fix $k\geq 2$ and suppose the claim is true for $k-2$. 

Let $v\in \Lambda_{\intform}^{k}$. Using 
\autoref{L:Fi-mapspartialhwtohw}, we find for 
each $i\geq 1$ such that $k-2i\geq 0$, $\xi_{i}$, a 
unit in $\intform$, such that $F(v_{k,k-2i})=\xi_{i}\cdot v_{k-2,k-2i}$. 
Since $F(v)\in \Lambda^{k-2}_{\intform}$, our hypothesis 
implies there are $u_{k-2,k-2i}\in U_{\intform}(\mathfrak{sp}_{2n})$ 
such that $F(v)= \sum_{i} u_{k-2, k-2i}\acts v_{k-2, k-2i}$. 
Since $F$ commutes with $U_{\intform}(\mathfrak{sp}_{2n})$, the vector
\begin{gather*}
v-\sum_{i}\xi_{i}^{-1}\cdot u_{k-2,k-2i}\acts v_{k,k-2i}\in\Lambda_{\intform}^{k}
\end{gather*}
is in the kernel of $F$. \autoref{P:kerisWeyl} implies that $F|_{\Lambda_{\intform}^{k}}$ is generated over $U_{\intform}(\mathfrak{sp}_{2n})$ by $v_{k,k}$. The result follows. 
\end{proof}

The following is our version of
\cite[Corollary 2.13]{Br-super-kl-glmn}:

\begin{Lemma}\label{L:brundan-lemma}
Let $\theta\colon\Lambda_{\intform}\to\Lambda_{\intform}$ be $U_{\intform}(\mathfrak{sp}_{2n})$-equivariant. 
Suppose that $\theta(v_{k,k-2i})=v_{k,k-2i}$ for 
all $k\in\{0,\dots,2n\}$ and $i\in\N$ 
such that $k-2i\geq 0$. Then $\theta=\id$. 
\end{Lemma}

\begin{proof}
Let $x\in\Lambda_{\intform}^{k}$. By \autoref{L:generation} there 
are $u_{k,k-2i}\in U_{\intform}(\mathfrak{sp}_{2n})$ such 
that $x=\sum_{i}u_{k,k-2i}\cdot v_{k,k-2i}$ with $u_{k,k-2i}\in\intform$. Thus,
\begin{gather*}
\theta(x) = \sum_{i}\theta(u_{k,k-2i}\cdot v_{k,k-2i})
=\sum_{i}u_{k,k-2i}\cdot\theta(v_{k,k-2i})
=\sum_{i}u_{k,k-2i}\cdot v_{k,k-2i}=x,
\end{gather*}
where we used that $\theta(v_{k,k-2i})=v_{k,k-2i}$.
\end{proof}

\begin{Definition}
Let $M$, $N$ be $\intform$-modules. A $\Z$-module homomorphism 
$a\colon M\to N$ such that $a(q\cdot m)=q^{-1}\cdot a(m)$ 
will be referred to as \emph{$\intform$-antilinear}.
\end{Definition}

\begin{Lemma}\label{L:anitlinear-uniqueness}
Suppose that $a$ and $b$ are two $\intform$-antilinear maps $a,b\colon
\Lambda_{\intform}\to\Lambda_{\intform}$ such that $\alpha\in\{a,b\}$:
\begin{enumerate}[itemsep=0.15cm,label=\emph{\upshape(\roman*)}]

\item Commutes with the Chevalley generators $e_{i}$ and $f_{i}$.

\item Satisfies $\alpha(k_{i}\acts x)=k_{i}^{-1}\acts\alpha(x)$ for all $x\in\Lambda_{\intform}$.

\item Fixes the vectors $\{v_{k,k-2i}|0\leq k\leq 2n,k-2i\geq 0\}$.

\end{enumerate}
Then $a=b$ and $a^{2}=\id$. 
\end{Lemma}

\begin{proof}
The composition of two $\intform$-antilinear maps is $\intform$-linear. 
Thus, \autoref{L:brundan-lemma} applies to $a^{2}$, showing that $a$ is an involution, 
and $a\circ b$, showing that $a\circ b=\id$. Algebra autopilot gives 
\begin{gather*}
b=\id\circ b=a\circ a\circ b=a\circ\id=a.
\end{gather*}
and the statement is proven.
\end{proof}

\begin{Definition}
Let $S=\{x_{1},x_{2},\dots,x_{|S|}\}\subset[1,-1]$ and 
suppose $x_{1}<x_{2}<\dots<x_{|S|}$. We will write
\begin{gather*}
v_{S}^{rev}:=v_{x_{|S|}}\cdots v_{x_{2}}v_{x_{1}}. 
\end{gather*}
for the \emph{reversed standard basis vectors}.
\end{Definition}

\begin{Lemma}
For all $k\in\{0,1,\dots,2n\}$, the sets
\begin{gather*}
\{
v_{S}^{rev}|S\subset[1,-1]
\}
,\quad
\{
v_{S}^{rev}|S\subset[1,-1],|S|=k
\}
\end{gather*}
form an $\intform$-basis for $\Lambda_{\intform}$ and $\Lambda_{\intform}^{k}$, respectively.
\end{Lemma}

\begin{proof}
Similarly as the proof of \autoref{P:ExtFlat} and omitted.
\end{proof}

\begin{Lemma}\label{L:opExtAct}
For $i\in\{1,\dots,n-1\}$ and $i=n$ we have:
\begin{align*}
f_{i}\acts v_{S}^{rev}&= 
\begin{cases}
v_{(S\setminus\{i\})\cup\{i+1\}}^{rev} & \text{if }S_{i,i+1}=\{i\},
\\
v_{(S\setminus\{-(i+1)\})\cup\{-i\}}^{rev} & \text{if }S_{i,i+1}=\{-(i+1)\}, 
\\
q^{-1}\cdot v_{(S\setminus\{-(i+1)\}\cup)\{-i\}}^{rev} 
+v_{(S\setminus\{i\}\cup)\{i+1\}}^{rev} & \text{if }S_{i,i+1}=\{i,-(i+1)\},
\\
v_{(S\setminus\{i\}\cup)\{i+1\}}^{rev} & \text{if }S_{i,i+1}=\{i,-i\},
\\
q\cdot v_{(S\setminus\{-(i+1)\}\cup)\{-i\}}^{rev} & \text{if }S_{i,i+1}=\{i+1,-(i+1)\},
\\
v_{(S\setminus\{i\}\cup)\{i+1\}}^{rev} & \text{if }S_{i,i+1}=\{i,-(i+1),-i\},
\\
v_{(S\setminus\{-(i+1)\}\cup)\{-i\}}^{rev} & \text{if }S_{i,i+1}=\{i,(i+1),-(i+1)\},
\\
0 & \text{otherwise}.
\end{cases}
\\
f_{n}\acts v_{S}^{rev}&= 
\begin{cases} 
v_{(S\setminus\{n\}\cup)\{-n\}}^{rev}& \text{if }S_{n}=\{n\},
\\
0 & \text{otherwise}.
\end{cases}
\end{align*}
Similarly, for $e_{i}\acts v_{S}^{rev}$ and $e_{n}\acts v_{S}^{rev}$.
\end{Lemma}

\begin{proof}
As for \autoref{L:ExtAction}.
\end{proof}

\begin{Definition}\label{D:involution}
We call the $\intform$-antilinear map defined by
\begin{gather*}
v_{S}\mapsto\overline{v_{S}} 
:=(-q^{-1})^{\binom{|S|}{2}}(q^{-1})^{\frac{|S\cap -S|}{2}}\cdot v_{S}^{rev}
\end{gather*}
the \emph{bar involution}.
\end{Definition}

Note that for $g(q)\in\intform$, we have $\overline{g(q)\cdot v_{S}}=g(q^{-1})\cdot\overline{v_{S}}$.
We are about to show that this map is an involution, justifying the choice of terminology.

\begin{Lemma}\label{L:bar-commutes-Chevalley}
We have $\overline{e_{i}\acts v_{S}}=e_{i}\acts\overline{v_{S}}$, 
$\overline{f_{i}\acts v_{S}}=f_{i}\acts\overline{v_{S}}$ and $\overline{k_{i}\acts v_{S}}= k_{i}^{-1}\acts\overline{v_{S}}$ for all $i\in\{1,\dots,n\}$.
\end{Lemma}

\begin{proof}
The proof is a direct calculation using \autoref{L:ExtAction}, \autoref{D:involution}, and \autoref{L:opExtAct}. We give an example of one of the calculations, leaving the rest to the reader. Suppose that $S_{i,i+1}=\{i,-(i+1)\}$ and write $T:=(S\setminus\{-(i+1)\})\cup\{-i\}$ and $U:=(S\setminus\{i\})\cup\{i+1\}$. Then 
\begin{align*}
\overline{f_{i}\acts v_{S}} &=\overline{v_{T}+q^{-1}\cdot v_{U}}
\\
&=(-q^{-1})^{\binom{|T|}{2}}(q^{-1})^{\frac{|T\cap -T|}{2}}\cdot v^{rev}_{T}+q\cdot (-q^{-1})^{\binom{|U|}{2}}(q^{-1})^{\frac{|U\cap -U|}{2}}\cdot v^{rev}_{U} 
\\
&=(-q^{-1})^{\binom{|S|}{2}}\cdot q^{-1}\cdot(q^{-1})^{\frac{|S\cap -S|}{2}}\cdot v^{rev}_{T}+q\cdot (-q^{-1})^{\binom{|S|}{2}}\cdot q^{-1}\cdot(q^{-1})^{\frac{|S\cap -S|}{2}}\cdot v^{rev}_{U}
\\
&=(-q^{-1})^{\binom{|S|}{2}} (q^{-1})^{\frac{|S\cap -S|}{2}}\cdot \big(q^{-1}\cdot v^{rev}_{T} + v^{rev}_{U}\big) 
\\
&=(-q^{-1})^{\binom{|S|}{2}}(q^{-1})^{\frac{|S\cap -S|}{2}}\cdot f_{i}\acts v_{S}^{rev} 
\\
&=f_{i}\acts\overline{v_{S}}.
\end{align*}
All other calculations are similar.
\end{proof}

\begin{Remark}\label{R:sometimesviv-icommutes}
First, note that the defining relations for $\Lambda_{\intform}$ imply that for all $k>0$, 
\begin{gather*}
v_{-(x+k)}v_{-(x+1)}\dots v_{-n}v_{n}\dots v_{x+1}=0,
\end{gather*}
and therefore $v_{x}v_{-x}v_{-(x+1)}\dots v_{-n}v_{n}\dots v_{x+1}
=-q^{2}\cdot v_{-x}v_{x}v_{-(x+1)}\dots v_{-n}v_{n}\dots v_{x+1}$.
\end{Remark}

\begin{Lemma}\label{L:bar-fixes-partial-hwvectors}
We have $\overline{v_{k,k-2i}}=v_{k,k-2i}$ for all 
$0\leq k\leq 2n$ and $i\in\N$ with $k-2i\in\N$.
\end{Lemma}

\begin{proof}
The relations for $\Lambda_{\intform}$ and \autoref{R:sometimesviv-icommutes} imply that 
\begin{align*}
v_{k,k-2i}^{rev}
&=(-q)^{(k-1)}v_{1}v_{-(n-i+1)}\dots v_{-n}v_{n}\dots v_{n-i+1}v_{k-2i}\dots v_{2} 
\\
&=(-q)^{(k-1)+(k-2)}v_{1}v_{2}v_{-(n-i+1)}\dots v_{-n}v_{n}\dots v_{n-i+1}v_{k-2i}\dots v_{3} 
\\
&=\dots 
\\
&=(-q)^{(k-1)+(k-2)+\dots+(k-(k-2i))}\cdot v_{1}\dots 
v_{k-2i}v_{-(n-i+1)}\dots v_{-n}v_{n}\dots v_{n-i+1} 
\\
&=(-q)^{(k-1)+(k-2)+\dots+(2i)+(2i-1)}q\cdot 
v_{1}\dots v_{k-2i}v_{n-i+1}v_{-(n-i+1)}v_{-(n-i+2)}\dots v_{-n}v_{n}\dots v_{n-i+2} 
\\
&=(-q)^{(k-1)+\dots+(2i-2)}q^2\cdot 
v_{1}\dots v_{k-2i}v_{n-i+1}v_{n-i+2}v_{-(n-i+1)}v_{-(n-i+2)}\dots v_{-n}v_{n}\dots v_{n-i+3} 
\\
&=\dots 
\\
&=(-q)^{(k-1)+\dots+(2i-i)}q^{i}\cdot 
v_{1}\dots v_{k-2i}v_{n-i+1}\dots v_{n}v_{-(n-i+1)}\dots v_{-n} 
\\
&=(-q)^{(k-1)+\dots+(i)+(i-1)}q^{i}\cdot 
v_{1}\dots v_{k-2i}v_{n-i+1}\dots v_{n}v_{-n}v_{-(n-i+1)}\dots v_{-(n-1)} 
\\
&=(-q)^{(k-1)+\dots+(i-2)}q^{i}\cdot 
v_{1}\dots v_{k-2i}v_{n-i+1}\dots v_{n}v_{-n}v_{-(n-1)}v_{-(n-i+1)}\dots v_{-(n-2)}
\\
&= \dots 
\\
&=(-q)^{(k-1)+\dots +(1)}q^{i}\cdot 
v_{1}\dots v_{k-2i}v_{n-i+1}\dots v_{n}v_{-n}\dots v_{-(n-i+1)} 
\\
&= (-q)^{\binom{k}{2}}q^{i}\cdot v_{k, k-2i}.
\end{align*}
By the definition of the bar involution we find
\begin{gather*}
\overline{v_{k,k-2i}}=(-q^{-1})^{\binom{k}{2}}(q^{-1})^{i}\cdot v_{k,k-2i}^{rev}=v_{k,k-2i}.
\end{gather*}
The proof completes.
\end{proof}

\begin{Proposition}
The bar involution is an $\intform$-antilinear involution.
\end{Proposition}

\begin{proof}
Combining \autoref{L:bar-commutes-Chevalley} and \autoref{L:bar-fixes-partial-hwvectors}, 
this follows from \autoref{L:anitlinear-uniqueness}.
\end{proof}

\begin{Lemma}\label{L:Bar}
The elements $b_{S}$ are invariant under the bar involution. 
\end{Lemma}

\begin{proof}
For brevity, we only sketch a proof.

Let $S\subset[1,-1]$. Suppose that $|S|=k$ and
$|\pm L(S)|+2i=|S_{0}|$. We will show that there is 
$u_{S}\in U_{\intform}(\mathfrak{sp}_{2n})$ of the form 
$u_{S}=f_{j_{1}}\dots f_{j_{t}}$, for some $j_{p}\in\{1,\dots,n\}$, 
such that $u_{S}\acts v_{k,k-2i}=b_{S}$. It will then follow 
from \autoref{L:bar-commutes-Chevalley} and \autoref{L:bar-fixes-partial-hwvectors} that 
\begin{gather*}
\overline{b_{S}}=\overline{u_{S}\acts v_{k, k-2i}}=
u_{S}\acts\overline{v_{k, k-2i}}=u_{S}\acts v_{k, k-2i}=b_{S}. 
\end{gather*}

To this end, let 
\begin{gather*}
\mathcal{C}:=\{S\subset[1,-1]\}, 
\quad\mathcal{C}^{k}:=\{S\subset[1,-1]\,|\,|S|=k\},
\\
\mathcal{C}^{k}_{k-2i}:=
\{S\subset[1,-1]\,|\,|S|=k\text{ and }|\pm L(S)|+2i=|S_{0}|\}.
\end{gather*}
Note that
\begin{gather*}
\mathcal{C}=\coprod_{0\leq k\leq 2n}\mathcal{C}^{k}=
\coprod_{\substack{0\leq k\leq 2n \\k-2i\geq 0}}\mathcal{C}^{k}_{k-2i}.
\end{gather*}
Define operators $\tilde{E},\tilde{F}\colon\mathcal{C}\to\mathcal{C}\cup\{0\}$ as follows:
\begin{gather*}
\tilde{E}\acts S:=
\begin{cases} 
S\cup\{\pm i_{k}\} & S_{0}\setminus \pm L(S)=\{\pm i_{1},\dots,\pm i_{r}\}
\text{ and }i_{1}<\dots<i_{r},
\\
0 & S_{0}=\pm L(S),
\end{cases}
\\
\tilde{F}\acts S:=
\begin{cases} 
S\setminus\{\pm i_{1}\} & S_{0}\setminus \pm L(S)=\{\pm i_{1},\dots,\pm i_{r}\}
\text{ and }i_{1}<\dots<i_{r}, 
\\
0 & S_{0}=\pm L(S).
\end{cases}
\end{gather*}
We write $\mathcal{C}(m):=\{-m,-m+2,\dots,m-2,m\}$ to denote the $\mathfrak{sl}_{2}$ crystal of highest weight $m$, with the evident action of $\tilde{E}$ an $\tilde{F}$. One can then use \autoref{D:TensorToDot} to determine a bijection
\begin{gather*}
\big(\mathcal{C}(0)\coprod\mathcal{C}(1)\coprod\mathcal{C}(0)\big)^{\otimes n}\to\mathcal{C},
\end{gather*}
which is easily seen to be an isomorphism of $\mathfrak{sl}_{2}$ crystals. 

For $S\in\mathcal{C}$ such that 
$S_{0}\setminus \pm L(S)=\{\pm i_{1},\dots,\pm i_{r}\}$, we have 
$S\setminus\{\pm i_{1},\dots,\pm i_{r}\}\in\mathcal{C}_{\varpi_{k}}$. 
Moreover, the assignments
\begin{gather*}
S\mapsto 
(S\setminus\{\pm i_{1},\dots,\pm i_{r}\},-n+|S|)\in
\mathcal{C}_{\varpi_{|S|-2r}}\otimes\mathcal{C}(n-|S|+2r)
\end{gather*}
induce a bijection 
\begin{gather*}
\mathcal{C}\rightarrow\coprod_{0\leq k\leq n}\mathcal{C}_{\varpi_{k}}\otimes\mathcal{C}(n-k).
\end{gather*}
The formulas in \autoref{D:subsetcrystal} define 
operators $\ee_{j},\ff_{j}\colon\mathcal{C}\to\mathcal{C}\cup\{0\}$, and we can then view the above bijection as an isomorphism of $\mathfrak{sp}_{2n}\times\mathfrak{sl}_{2}$ crystals. 

For $S\in\mathcal{C}$ such that $|S|=k$ and $|S_{0}\setminus \pm L(S)|=r$, then from \autoref{L:generatingker} we know there is 
$u_{\tilde{F}^{r}\acts S}=f_{j_{1}}\dots f_{j_{t}}$ 
such that $u_{\tilde{F}^{r}\acts S}\acts v_{k-2r,k-2r}=b_{\tilde{F}^{r}\acts S}$. 
Now, similar to the proof of \autoref{L:generatingker}, one 
compares the descriptions of $f_{j}\acts b_{S}$ and $\ff_{j}\acts S$, 
but now for $S\in\mathcal{C}$, to deduce that 
$u_{\tilde{F}^{r}\acts S}\acts v_{k, k-2r}=b_{S}$, as desired. 
\end{proof}

Given two subsets $X,Y\subset \{1,\dots,n\}$, we write 
$X\leq Y$ if $X=\{x_{1}<x_{2}<\dots<x_{m}\}$, 
$Y=\{y_{1}< y_{2}<\dots<y_{m}\}$, and $x_{k}\leq y_{k}$ for $k\in\{1,\dots,m\}$. 

\begin{Definition}\label{D:PartialOrder}
We define a partial order on 
the standard basis $\{v_{S}|S\subset[1,-1]\}$ by $v_{S}\leq v_{T}$ if and only if:
\begin{enumerate}[itemsep=0.15cm,label=\emph{\upshape(\roman*)}]

\item $S\setminus S_{0}=T\setminus T_{0}$.

\item $S_{0}\leq T_{0}$.

\end{enumerate}
(Note that $v_{S}\leq v_{T}$ is equivalent 
to $\wt v_{S}=\wt v_{T}$ and $S_{0}\leq T_{0}$.)
\end{Definition}

\begin{Lemma}\label{L:UniqueForm}
There is a unique symmetric bilinear form 
\begin{gather*}
(\placeholder,\placeholder)\colon
\Lambda_{\intform}\times\Lambda_{\intform}\to\intform
\end{gather*}
determined by $(v_{S},v_{T})=\delta_{S,T}$.
\end{Lemma}

\begin{proof}
Clear.
\end{proof}

Using \autoref{L:UniqueForm} we define:

\begin{Definition}
We define the sesquilinear form 
\begin{align*}
\langle\placeholder,\placeholder\rangle\colon
\Lambda_{\intform}\times\Lambda_{\intform}\to\intform
,\quad
\langle x,y\rangle:=(x,\overline{y}). 
\end{align*}
(\emph{Sesquilinear} means $q$-linear in first 
argument and $q$-antilinear in the second argument).
\end{Definition}

Note that $\langle\overline{y},\overline{x}\rangle=\langle x,y\rangle$
and $\langle v_{S},\overline{v_{T}}\rangle=\delta_{S,T}$.

\begin{Lemma}
Let $S\subset[1,-1]$. We have
\begin{gather*}
\overline{v_{S}}
\in v_{S}+\sum_{v_{T}\geq v_{S}}\intform\cdot v_{T}. 
\end{gather*}
In particular, if $S$ does not have any fully dotted columns, then $\overline{v_{S}}=v_{S}$. 
\end{Lemma}

\begin{proof}
The proof is similar to the proof of \autoref{L:bar-fixes-partial-hwvectors}. That is, going from $\overline{v_{S}}$ to $v_{S}$
can be done using the defining relations of $\Lambda_{\intform}$ 
to reverse the order. The 
only relation one needs to be careful in the process is
\begin{gather*}
v_{-i}v_{i}=-q^{2}\cdot v_{i}v_{-i}
+(q-q^{-1})
\sum_{k\in[1,n-i]}(-q)^{k+1}\cdot v_{i+k}v_{-(i+k)}.
\end{gather*}
The term in this relation that is multiplied by $(q-q^{-1})$ produces 
fully dotted columns further to the right. Applying this observation repeatedly 
shows the statement by additionally observing that the scalar in front of the leading term 
comes out as claimed.
\end{proof}

\begin{Lemma}
The free $\intform$-module $\Lambda_{\intform}$ 
equipped with standard basis $\{v_{S}|S\subset[1,-1]\}$, 
partial order $\leq$ from 
\autoref{D:PartialOrder}, and bar involution $v_{S}\mapsto\overline{v_{S}}$ 
equip $\Lambda_{\intform}$ satisfies the conditions of a 
\emph{balanced pre-canonical structure} in the terminology of 
\cite[Definitions 1.1 and 1.5]{We-canonical-bases-higher-rep}.
(The balanced part comes from $\langle v_{S},\overline{v_{T}}\rangle=\delta_{S,T}$.)
\end{Lemma}

\begin{proof}
Directly from the definitions.
\end{proof}

By \emph{canonical} we mean \cite[Definition 1.7]{We-canonical-bases-higher-rep}.

\begin{Theorem}
The basis $\{b_{S}|S\subset[1,-1]\}$ is canonical.
\end{Theorem}

\begin{proof}
Because of 
\cite[Lemma 1.8]{We-canonical-bases-higher-rep}, this follows 
from  
\begin{gather*}
b_{S}\in v_{S}+\sum_{T\geq S}q^{-1}\Z[q^{-1}]\cdot v_{T}.
\end{gather*}
and bar invariance of $b_{S}$ as in \autoref{L:Bar}.
\end{proof}

\begin{Remark}
Using \cite[Theorem 1.9]{We-canonical-bases-higher-rep} 
one can now show that $\{b_{S}|S\subset[1,-1]\}$ is 
\textit{the} canonical basis (up to signs). Details are 
omitted.
\end{Remark}

\section*{Declarations}

\subsection*{Ethical Approval} 
Not applicable.

\subsection*{Competing interests} 
No financial or personal competing interest.

\subsection*{Authors' contributions }
All authors contributed equally with respect to every section of the paper.

\subsection*{Funding} 
E.B. was supported by the National Science Foundation's M.S.P.R.F.-2202897, and by the University of Oregon's Lokey Fellowship. D.T. shouldn't be supported but they were still supported by \ochanged{the ARC Future Fellowship FT230100489.}

\subsection*{Availability of data and materials} 
Not applicable.


\begin{thebibliography}{STWZ21}

\bibitem[AR96]{AdRy-tilting-howe-positive-char}
A.M.~Adamovich and G.L.~Rybnikov.
\newblock Tilting modules for classical groups and {H}owe duality in positive
characteristic.
\newblock {\em Transform. Groups}, 1(1-2):1--34, 1996.
\newblock \href {https://doi.org/10.1007/BF02587733}
{\path{doi:10.1007/BF02587733}}.

\bibitem[APW91]{AnPoWe-representation-qalgebras}
H.H.~Andersen, P.~Polo, and K.X.~Wen.
\newblock Representations of quantum algebras.
\newblock {\em Invent. Math.}, 104(1):1--59, 1991.
\newblock \href {https://doi.org/10.1007/BF01245066}
{\path{doi:10.1007/BF01245066}}.

\bibitem[AST18]{AnStTu-cellular-tilting}
H.H.~Andersen, C.~Stroppel, and D.~Tubbenhauer.
\newblock Cellular structures using {$\mathrm{U}_q$}-tilting modules.
\newblock {\em Pacific J. Math.}, 292(1):21--59, 2018.
\newblock URL: \url{https://arxiv.org/abs/1503.00224}, \href
{https://doi.org/10.2140/pjm.2018.292.21}
{\path{doi:10.2140/pjm.2018.292.21}}.

\bibitem[AST17]{AnStTu-semisimple-tilting}
H.H.~Andersen, C.~Stroppel, and D.~Tubbenhauer.
\newblock Semisimplicity of {H}ecke and (walled) {B}rauer algebras.
\newblock {\em J. Aust. Math. Soc.}, 103(1):1--44, 2017.
\newblock URL: \url{https://arxiv.org/abs/1507.07676}, \href
{https://doi.org/10.1017/S1446788716000392}
{\path{doi:10.1017/S1446788716000392}}.

\bibitem[AT17]{AnTu-tilting}
H.H.~Andersen and D.~Tubbenhauer.
\newblock Diagram categories for {$\textbf{U}_q$}-tilting modules at roots of
unity.
\newblock {\em Transform. Groups}, 22(1):29--89, 2017.
\newblock URL: \url{http://arxiv.org/abs/1409.2799}, \href
{https://doi.org/10.1007/s00031-016-9363-z}
{\path{doi:10.1007/s00031-016-9363-z}}.

\bibitem[BLM90]{BeLuMaPh-quantum-group}
A.A.~Beilinson, G.~Lusztig, and R.~MacPherson.
\newblock A geometric setting for the quantum deformation of {${\rm GL}_n$}.
\newblock {\em Duke Math. J.}, 61(2):655--677, 1990.
\newblock \href {https://doi.org/10.1215/S0012-7094-90-06124-1}
{\path{doi:10.1215/S0012-7094-90-06124-1}}.

\bibitem[BZ08]{BeZw-braided-ext-sym}
A.~Berenstein and S.~Zwicknagl.
\newblock Braided symmetric and exterior algebras.
\newblock {\em Trans. Amer. Math. Soc.}, 360(7):3429--3472, 2008.
\newblock URL: \url{https://arxiv.org/abs/math/0504155}, \href
{https://doi.org/10.1090/S0002-9947-08-04373-0}
{\path{doi:10.1090/S0002-9947-08-04373-0}}.

\bibitem[Ber78]{Be-diamond-lemma}
G.M.~Bergman.
\newblock The diamond lemma for ring theory.
\newblock {\em Adv. in Math.}, 29(2):178--218, 1978.
\newblock \href {https://doi.org/10.1016/0001-8708(78)90010-5}
{\path{doi:10.1016/0001-8708(78)90010-5}}.

\bibitem[Bod20]{Bo-c2-tilting}
E.~Bodish.
\newblock Web calculus and tilting modules in type {$C_2$}.
\newblock {\em Quantum Topol.} 13 (2022), no. 3, 407--458.
\newblock URL: \url{https://arxiv.org/abs/2009.13786}, 
\href {https://doi.org/10.4171/qt/166}
{\path{doi:10.4171/qt/166}}.

\bibitem[BERT21]{BoElRoTa-c-webs}
E.~Bodish, B.~Elias, D.E.V.~Rose, and L.~Tatham.
\newblock Type {C} webs.
\newblock 2021.
\newblock URL: \url{https://arxiv.org/abs/2103.14997}.

\bibitem[Bou02]{Bo-chapters-4-6}
N.~Bourbaki.
\newblock {\em Lie groups and {L}ie algebras. {C}hapters 4--6}.
\newblock Elements of Mathematics (Berlin). Springer-Verlag, Berlin, 2002.
\newblock Translated from the 1968 French original by Andrew Pressley.
\newblock \href {https://doi.org/10.1007/978-3-540-89394-3}
{\path{doi:10.1007/978-3-540-89394-3}}.

\bibitem[Bru03]{Br-super-kl-glmn}
J.~Brundan.
\newblock Kazhdan--{L}usztig polynomials and character formulae for the {L}ie
superalgebra {$\mathfrak{gl}(m|n)$}.
\newblock {\em J. Amer. Math. Soc.}, 16(1):185--231, 2003.
\newblock URL: \url{https://arxiv.org/abs/math/0203011}, \href
{https://doi.org/10.1090/S0894-0347-02-00408-3}
{\path{doi:10.1090/S0894-0347-02-00408-3}}.

\bibitem[BS17]{BuSc-crystal-bases}
D.~Bump and A.~Schilling.
\newblock {\em Crystal bases}.
\newblock World Scientific Publishing Co. Pte. Ltd., Hackensack, NJ, 2017.
\newblock Representations and combinatorics.
\newblock \href {https://doi.org/10.1142/9876} {\path{doi:10.1142/9876}}.

\bibitem[CKM14]{CaKaMo-webs-skew-howe}
S.~Cautis, J.~Kamnitzer, and S.~Morrison.
\newblock Webs and quantum skew {H}owe duality.
\newblock {\em Math. Ann.}, 360(1-2):351--390, 2014.
\newblock URL: \url{https://arxiv.org/abs/1210.6437}, \href
{https://doi.org/10.1007/s00208-013-0984-4}
{\path{doi:10.1007/s00208-013-0984-4}}.

\bibitem[dCP76]{CoPr-invariant-theory}
C.~de~Concini and C.~Procesi.
\newblock A characteristic free approach to invariant theory.
\newblock {\em Advances in Math.}, 21(3):330--354, 1976.
\newblock \href {https://doi.org/10.1016/S0001-8708(76)80003-5}
{\path{doi:10.1016/S0001-8708(76)80003-5}}.

\bibitem[Don93]{Do-tilting-alg-groups}
S.~Donkin.
\newblock On tilting modules for algebraic groups.
\newblock {\em Math. Z.}, 212(1):39--60, 1993.
\newblock \href {https://doi.org/10.1007/BF02571640}
{\path{doi:10.1007/BF02571640}}.

\bibitem[Don98]{Do-q-schur}
S.~Donkin.
\newblock {\em The {$q$}-{S}chur algebra}, volume 253 of {\em London
Mathematical Society Lecture Note Series}.
\newblock Cambridge University Press, Cambridge, 1998.
\newblock \href {https://doi.org/10.1017/CBO9780511600708}
{\path{doi:10.1017/CBO9780511600708}}.

\bibitem[Don00]{Do-sp-fund}
R.G.~Donnelly.
\newblock Explicit constructions of the fundamental representations of the
symplectic {L}ie algebras.
\newblock {\em J. Algebra}, 233(1):37--64, 2000.
\newblock \href {https://doi.org/10.1006/jabr.2000.8446}
{\path{doi:10.1006/jabr.2000.8446}}.

\bibitem[Dot03]{Do-presenting-qSchur}
S.~Doty.
\newblock Presenting generalized {$q$}-{S}chur algebras.
\newblock {\em Represent. Theory}, 7:196--213, 2003.
\newblock URL: \url{https://arxiv.org/abs/math/0305208}, \href
{https://doi.org/10.1090/S1088-4165-03-00176-6}
{\path{doi:10.1090/S1088-4165-03-00176-6}}.

\bibitem[DPS98]{DuPaSc-schur-weyl-tilting}
J.~Du, B.~Parshall, and L.~Scott.
\newblock Quantum {W}eyl reciprocity and tilting modules.
\newblock {\em Comm. Math. Phys.}, 195(2):321--352, 1998.
\newblock \href {https://doi.org/10.1007/s002200050392}
{\path{doi:10.1007/s002200050392}}.

\bibitem[ES18]{EhSt-nw-algebras-howe}
M.~Ehrig and C.~Stroppel.
\newblock Nazarov--{W}enzl algebras, coideal subalgebras and categorified skew {H}owe duality.
\newblock {\em Adv. Math.}, 331:58--142, 2018.
\newblock URL: \url{https://arxiv.org/abs/1310.1972}, \href
{https://doi.org/10.1016/j.aim.2018.01.013}
{\path{doi:10.1016/j.aim.2018.01.013}}.

\bibitem[Eli15]{El-ladders-clasps}
B.~Elias.
\newblock Light ladders and clasp conjectures.
\newblock 2015.
\newblock URL: \url{https://arxiv.org/abs/1510.06840}.

\bibitem[Fou05]{Fo-symplectic-simple}
S.~Foulle.
\newblock Characters of the irreducible representations with fundamental
highest weight for the symplectic group in characteristic p.
\newblock 2005.
\newblock See also the authors Ph.D. thesis, Universit{\'e} Claude Bernard Lyon 1.
\newblock URL: \url{https://arxiv.org/abs/math/0512312}.

\bibitem[GM13]{GiMa-modular-dim-symplectic}
P.M.~Gilmer and G.~Masbaum.
\newblock Dimension formulas for some modular representations of the symplectic
group in the natural characteristic.
\newblock {\em J. Pure Appl. Algebra}, 217(1):82--86, 2013.
\newblock URL: \url{https://arxiv.org/abs/1111.0240}, \href
{https://doi.org/10.1016/j.jpaa.2012.04.005}
{\path{doi:10.1016/j.jpaa.2012.04.005}}.

\bibitem[How95]{Ho-perspectives-invariant-theory}
R.~Howe.
\newblock Perspectives on invariant theory: {S}chur duality, multiplicity-free actions and beyond.
\newblock In {\em The {S}chur lectures (1992) ({T}el {A}viv)}, volume~8 of {\em Israel Math. Conf. Proc.}, pages 1--182. Bar-Ilan Univ., Ramat Gan, 1995.

\bibitem[How89]{Ho-remarks-invariant-theory}
R.~Howe.
\newblock Remarks on classical invariant theory.
\newblock {\em Trans. Amer. Math. Soc.}, 313(2):539--570, 1989.
\newblock \href {https://doi.org/10.2307/2001418} {\path{doi:10.2307/2001418}}.

\bibitem[Jan96]{Ja-lectures-qgroups}
J.C.~Jantzen.
\newblock {\em Lectures on quantum groups}, volume~6 of {\em Graduate Studies
in Mathematics}.
\newblock American Mathematical Society, Providence, RI, 1996.

\bibitem[Kan98]{Ka-Based-Filt}
M.~Kaneda.
\newblock Based modules and good filtrations in algebraic groups.
\newblock {\em Hiroshima Math. J.}, 28(2):337--344, 1998.
\newblock URL: \url{http://projecteuclid.org/euclid.hmj/1206126765}.

\bibitem[LTV22]{LaTuVa-verma-howe}
A.~Lacabanne, D.~Tubbenhauer, and P.~Vaz.
\newblock Verma {H}owe duality and {LKB} representations.
\newblock 2022.
\newblock URL: \url{https://arxiv.org/abs/2207.09124}.

\bibitem[LZZ11]{LeZhZh-q-first-fundamental-theorem}
G.I.~Lehrer, H.~Zhang, and R.B.~Zhang.
\newblock A quantum analogue of the first fundamental theorem of classical
invariant theory.
\newblock {\em Comm. Math. Phys.}, 301(1):131--174, 2011.
\newblock URL: \url{https://arxiv.org/abs/0908.1425}, \href
{https://doi.org/10.1007/s00220-010-1143-3}
{\path{doi:10.1007/s00220-010-1143-3}}.

\bibitem[Lus10]{Lu-quantumgroups-book}
G.~Lusztig.
\newblock {\em Introduction to quantum groups}.
\newblock Modern Birkh\"{a}user Classics. Birkh\"{a}user/Springer, New York,
2010.
\newblock Reprint of the 1994 edition.
\newblock \href {https://doi.org/10.1007/978-0-8176-4717-9}
{\path{doi:10.1007/978-0-8176-4717-9}}.

\bibitem[Lus90]{Lu-qgroups-root-of-1}
G.~Lusztig.
\newblock Quantum groups at roots of {$1$}.
\newblock {\em Geom. Dedicata}, 35(1-3):89--113, 1990.
\newblock \href {https://doi.org/10.1007/BF00147341}
{\path{doi:10.1007/BF00147341}}.

\bibitem[McN00]{McNi-charp-howe}
G.J.~McNinch.
\newblock Filtrations and positive characteristic {H}owe duality.
\newblock {\em Math. Z.}, 235(4):651--685, 2000.
\newblock \href {https://doi.org/10.1007/s002090000157}
{\path{doi:10.1007/s002090000157}}.

\bibitem[MM90]{MiMi-affine-sln-crystal}
K.~Misra and T.~Miwa.
\newblock Crystal base for the basic representation of
{$U_{q}(\hat{\mathfrak{sl}}(n))$}.
\newblock {\em Comm. Math. Phys.}, 134(1):79--88, 1990.

\bibitem[NUW96]{NoUmWa-sl2-son-duality}
M.~Noumi, T.~Umeda, and M.~Wakayama.
\newblock Dual pairs, spherical harmonics and a {C}apelli identity in quantum
group theory.
\newblock {\em Compositio Math.}, 104(3):227--277, 1996.

\bibitem[Oli65]{Ol-gen-powers}
G.~Olive.
\newblock Generalized powers.
\newblock {\em Amer. Math. Monthly}, 72:619--627, 1965.
\newblock \href {https://doi.org/10.2307/2313851} {\path{doi:10.2307/2313851}}.

\bibitem[Par94]{Pa-tilting-tensor}
J.~Paradowski.
\newblock Filtrations of modules over the quantum algebra.
\newblock In {\em Algebraic groups and their generalizations: quantum and
infinite-dimensional methods ({U}niversity {P}ark, {PA}, 1991)}, volume~56 of
{\em Proc. Sympos. Pure Math.}, pages 93--108. Amer. Math. Soc., Providence,
RI, 1994.

\bibitem[PS83]{PrSu-fund-weyl-sp}
A.A.~Premet and I.D.~Suprunenko.
\newblock The {W}eyl modules and the irreducible representations of the
symplectic group with the fundamental highest weights.
\newblock {\em Comm. Algebra}, 11(12):1309--1342, 1983.
\newblock \href {https://doi.org/10.1080/00927878308822907}
{\path{doi:10.1080/00927878308822907}}.

\bibitem[Rin91]{Ri-good-filtrations}
C.M.~Ringel.
\newblock The category of modules with good filtrations over a quasi-hereditary
algebra has almost split sequences.
\newblock {\em Math. Z.}, 208(2):209--223, 1991.
\newblock \href {https://doi.org/10.1007/BF02571521}
{\path{doi:10.1007/BF02571521}}.

\bibitem[RT16]{RoTu-symmetric-howe}
D.E.V.~Rose and D.~Tubbenhauer.
\newblock Symmetric webs, {J}ones--{W}enzl recursions, and {$q$}-{H}owe
duality.
\newblock {\em Int. Math. Res. Not. IMRN}, (17):5249--5290, 2016.
\newblock URL: \url{https://arxiv.org/abs/1501.00915}, \href
{https://doi.org/10.1093/imrn/rnv302} {\path{doi:10.1093/imrn/rnv302}}.

\bibitem[RH03]{RyHa-q-kempf}
S.~Ryom-Hansen.
\newblock A {$q$}-analogue of {K}empf's vanishing theorem.
\newblock {\em Mosc. Math. J.}, 3(1):173--187, 260, 2003.
\newblock URL: \url{https://arxiv.org/abs/0905.0236}, \href
{https://doi.org/10.17323/1609-4514-2003-3-1-173-187}
{\path{doi:10.17323/1609-4514-2003-3-1-173-187}}.

\bibitem[ST19]{SaTu-bcd-webs}
A.~Sartori and D.~Tubbenhauer.
\newblock Webs and {$q$}-{H}owe dualities in types {BCD}.
\newblock {\em Trans. Amer. Math. Soc.}, 371(10):7387--7431, 2019.
\newblock URL: \url{https://arxiv.org/abs/1701.02932}, \href
{https://doi.org/10.1090/tran/7583} {\path{doi:10.1090/tran/7583}}.

\bibitem[Sug95]{Su-harmonic-analysis-son}
T.~Sugitani.
\newblock Harmonic analysis on quantum spheres associated with representations
of {$U_{q}(\mathfrak{so}_{N})$} and {$q$}-{J}acobi polynomials.
\newblock {\em Compositio Math.}, 99(3):249--281, 1995.

\bibitem[STWZ23]{SuTuWeZh-mixed-tilting}
L.~Sutton, D.~Tubbenhauer, P.~Wedrich, and J.~Zhu.
\newblock {SL2} tilting modules in the mixed case.
\newblock {\em Selecta Math. (N.S.)}, 29(3):39, 2023.
\newblock URL: \url{https://arxiv.org/abs/2105.07724}, \href
{https://doi.org/10.1007/s00029-023-00835-0}
{\path{doi:10.1007/s00029-023-00835-0}}.

\bibitem[TVW17]{TuVaWe-super-howe}
D.~Tubbenhauer, P.~Vaz, and P.~Wedrich.
\newblock Super {$q$}-{H}owe duality and web categories.
\newblock {\em Algebr. Geom. Topol.}, 17(6):3703--3749, 2017.
\newblock URL: \url{https://arxiv.org/abs/1504.05069}, \href
{https://doi.org/10.2140/agt.2017.17.3703}
{\path{doi:10.2140/agt.2017.17.3703}}.

\bibitem[TW21]{TuWe-quiver-tilting}
D.~Tubbenhauer and P.~Wedrich.
\newblock Quivers for {$\mathrm{SL}_2$} tilting modules.
\newblock {\em Represent. Theory}, 25:440--480, 2021.
\newblock URL: \url{https://arxiv.org/abs/1907.11560}, \href
{https://doi.org/10.1090/ert/569} {\path{doi:10.1090/ert/569}}.

\bibitem[Web15]{We-canonical-bases-higher-rep}
B.~Webster.
\newblock Canonical bases and higher representation theory.
\newblock {\em Compos. Math.}, 151(1):121--166, 2015.
\newblock URL: \url{https://arxiv.org/abs/1209.0051}, \href
{https://doi.org/10.1112/S0010437X1400760X}
{\path{doi:10.1112/S0010437X1400760X}}.

\bibitem[Zwi09]{Zw-r-matrix-poisson}
S.~Zwicknagl.
\newblock {$R$}-matrix {P}oisson algebras and their deformations.
\newblock {\em Adv. Math.}, 220(1):1--58, 2009.
\newblock URL: \url{https://arxiv.org/abs/0706.0351}, \href
{https://doi.org/10.1016/j.aim.2008.08.006}
{\path{doi:10.1016/j.aim.2008.08.006}}.

\end{thebibliography}
\end{document}